\documentclass[12pt]{amsart}
\usepackage{latexsym, amssymb, amscd, rotating}
\renewcommand\a{\alpha}
\renewcommand\b{\beta}
\newcommand\g{\gamma}
\renewcommand\d{\delta}
\newcommand\la{\lambda}
\newcommand\z{\zeta}
\newcommand\e{\eta}
\renewcommand\th{\theta}
\newcommand\io{\iota}
\newcommand\m{\mu}
\newcommand\n{\nu}
\newcommand\s{\sigma}
\newcommand\x{\chi}
\newcommand\f{\phi}
\newcommand\vf{\varphi}
\newcommand\p{\psi}
\renewcommand\t{\tau}
\renewcommand\r{\rho}

\newcommand\Om{\Omega}
\newcommand\w{\omega}
\newcommand\vS{\varSigma}
\newcommand\D{\Delta}
\newcommand\vD{\varDelta}

\newcommand\vL{\varLambda}

\newcommand{\vT}{\varTheta}
\newcommand\vG{\varGamma}
\newcommand\ve{\varepsilon}

\newcommand{\ZZ}{\mathbb Z}

\newcommand\Fq{{\mathbf F}_q}

\newcommand\Ql{\bar{\mathbf Q}_l}
\newcommand\BQ{\mathbf Q}
\newcommand\BP{\mathbf P}
\newcommand\BF{\mathbf F}

\newcommand\BR{\mathbf R}
\newcommand\BZ{\mathbf Z}

\newcommand\Bw{\mathbf w}

\newcommand\CA{\mathcal{A}}

\newcommand\CH{\mathcal{H}}
\newcommand\CI{\mathcal{I}}
\newcommand\CE{\mathcal{E}}
\newcommand\CL{\mathcal{L}}
\newcommand\DD{\mathcal{D}}
\newcommand\CS{\mathcal{S}}

\newcommand\CM{\mathcal{M}}

\newcommand\CO{\mathcal{O}}

\newcommand\CP{\mathcal{P}}
\newcommand\CV{\mathcal{V}}

\newcommand\CF{\mathcal{F}}

\newcommand\CT{ \mathcal{T}}
\newcommand\CZ{ \mathcal{Z}}
\newcommand\CX{ \mathcal{X}}
\newcommand\CY{ \mathcal{Y}}
\newcommand\CW{ \mathcal{W}}

\newcommand\Fu{\mathfrak u}
\newcommand\Fg{\mathfrak g}
\newcommand\Fs{\mathfrak s}
\newcommand\Fl{\mathfrak l}

\newcommand\FS{\mathfrak S}

\newcommand\Fz{\mathfrak z}

\newcommand\iv{^{-1}}
\newcommand\td{\tilde}
\newcommand\wh{\widehat}
\newcommand\wt{\widetilde}
\newcommand\wg{^{\wedge}}

\newcommand\ol{\overline}
\newcommand\hra{\hookrightarrow}
\newcommand\hla{\hookleftarrow}
\newcommand\lra{\leftrightarrow}

\newcommand\ssim{/\!\!\sim}
\newcommand\ch{\operatorname{ch}}
\newcommand\Id{\operatorname{Id}}
\newcommand\IC{\operatorname{IC}}
\newcommand\Ker{\operatorname{Ker}}

\newcommand\End{\operatorname{End}}
\newcommand\ind{\operatorname{ind}}
\newcommand\Ind{\operatorname{Ind}}

\newcommand\Irr{\operatorname{Irr}}
\newcommand\codim{\operatorname{codim}}
\newcommand\supp{\operatorname{supp}\,}
\newcommand\Lie{\operatorname{Lie}}
\newcommand\Tr{\operatorname{Tr}\,}

\newcommand\Ad{\operatorname{Ad}}
\newcommand\ad{\operatorname{ad}}

\newcommand\der{_{\operatorname{der}}}

\newcommand\uni{_{\operatorname{uni}}}
\newcommand\nil{_{\operatorname{nil}}}

\newcommand\lp{\operatorname{\!\langle\!}}
\newcommand\rp{\operatorname{\!\rangle}}

\newcommand\Diag{\operatorname{Diag}}
\newcommand\nat{^{\natural}}

\newcommand\dw{\dot w}

\newcommand\dc{\dot c}

\newcommand\ds{\dot s}
\newcommand\dx{\dot x}
\newcommand\dy{\dot y}
\newcommand\dz{\dot z}

\newcommand{\isom}{\,\raise2pt\hbox{$\underrightarrow{\sim}$}\,}
\numberwithin{equation}{section}

\newtheorem{thm}{Theorem}[section]
\newtheorem{lem}[thm]{Lemma}
\newtheorem{cor}[thm]{Corollary}
\newtheorem{prop}[thm]{Proposition}

\def \para#1{\par\medskip\textbf{#1}
              \addtocounter{thm}{1}}

\def \remark#1{\par\medskip\noindent
                \textbf{Remark #1}
                \addtocounter{thm}{1}}

\begin{document}
\setlength{\baselineskip}{4.9mm}
\setlength{\abovedisplayskip}{4.5mm}
\setlength{\belowdisplayskip}{4.5mm}
\renewcommand{\theenumi}{\roman{enumi}}
\renewcommand{\labelenumi}{(\theenumi)}
\renewcommand{\thefootnote}{\fnsymbol{footnote}}

\parindent=20pt
\pagestyle{myheadings}
\markboth{SHOJI}{LUSZTIG'S CONJECTURE}
\begin{center}
{\bf Lusztig's conjecture for finite special linear groups} \\
\vspace{1cm}
Toshiaki Shoji \\
\vspace{0.5cm}
Graduate School of Mathematics \\
Nagoya University  \\
Chikusa-ku, Nagoya 464-8602,  Japan
\end{center}
\title{}
\maketitle
\begin{abstract}
In this paper, we prove Lusztig's conjecture for 
$G^F = SL_n(\BF_q)$, i.e., we show that characteristic 
functions of character sheaves of $G^F$ coincide with 
almost characters of $G^F$ up to scalar constants, 
assuming that $\ch \BF_q$ is not too small.  
We determine these scalars explicitly.  Our result 
gives a method of computing irreducible characters of 
$G^F$. 
\end{abstract}
\par\bigskip
\addtocounter{section}{-1}
\section{Introduction}
Let $G = SL_n$ defined over a finite field $\BF_q$ with 
the standard Frobenius map $F$, and $G^F = SL_n(\BF_q)$ 
the finite special linear group.   In [S2], 
a parametrization of irreducible characters of $G^F$ was given 
by making use of modified generalized Gelfand-Graev characters. Also the 
almost characters of $G^F$ was defined, and it was shown that 
the Shintani descent of irreducible characters of $G^{F^m}$, 
for sufficiently divisible $m$,   
coincides with almost characters up to scalar.  
However as explained in the remark of the last part of [S2], 
the relationship of our parametrization of irreducible characters
and the parametrization in terms of the Harish-Chandra induction 
was not so clear.  
Now Lusztig's conjecture is formulated in the form that 
almost characters coincide with the characteristic functions 
of character sheaves of $G$ under a suitable parametrization. 
\par
In this paper we prove that Lusztig's conjecture holds for $G^F$,
 assuming that $\ch \BF_q$ is not too small so that Lusztig's 
results ([L7]) for generalized Gelfand-Graev characters are 
applicable.  In the course of the proof, it is shown that 
almost characters of $G^F$ are parametrized in terms of 
the twisted induction, which is compatible with the parametrization 
of $F$-stable irreducible characters of $G^{F^m}$ in terms of 
the Harish-Chandra induction for sufficiently 
divisible $m$.  Thus giving the relationship between two
parametrizations of $F$-stable irreducible characters of $G^{F^m}$
is equivalent to giving the relationship between two parametrizations of
almost characters of $G^F$.   
\par
We give a complete description of the relationship between those
two parametrizations of almost characters, by computing the
inner products of the characteristic functions of characters sheaves
with various modified generalized Gelfand-Graev characters.
Note that the inner product of characteristic functions with 
generalized Gelfand-Graev characters can be computed by using  
Lusztig's formula without difficulty.  However the computation 
in the case of modified generalized Gelfand-Graev characters 
is much more complicated since it involves  non-uni potent supports. 
Through this computation, we can describe the scalar constants appearing in 
Lusztig's conjecture.  Although the expression of 
such scalars are not so simple, they are explicitly computable.
(Here we need a result of Digne-Lehrer-Michel [DLM1] that the 
fourth root of unity occurring in the Lusztig's theory in [L7] 
is explicitly determined in the case of $SL_n$). 
\par
The computation of irreducible characters of $G^F$ is reduced to
the computation of characteristic functions of character sheaves.
In turn, the computation of those characteristic functions is 
reduced to the computation of generalized Green functions. 
Lusztig's algorithm of computing generalized Green functions  
contains certain unknown constants.  In the case of $SL_n$, 
we can determine such scalars, which will be discussed in [S3].
Thus our result makes it possible to compute the character table of 
$SL_n(\BF_q)$.  
\par\bigskip\noindent
{\bf Some notations.} \
For a finite group $\vG$, we denote by $\Irr \vG$ or $\vG\wg$ the set
of irreducible characters of $\vG$ over $\Ql$. 
If $F : \vG \to \vG$ is an
automorphism on $\vG$, we denote by $\vG\ssim_F$ the set of
$F$-twisted conjugacy classes in $\vG$, where $x, y \in \vG$ is 
$F$-twisted conjugate if there exists $z \in \vG$ such that
$y = z\iv xF(z)$.
In the case where $\vG$ is abelian, $\vG\ssim_F$ is naturally 
identified with the largest quotient of $\vG$ on which $F$ acts
trivially, which we denote by $\vG_F$. 
\par
For a reductive group $H$, we denote by $Z_H$ the center of $H$, and
denote by $Z_H^0$ the identity component of $Z_H$.
\par\vspace{1cm}
{\bf Contents}
\par\medskip
1. Parametrization of irreducible characters
\par
2. Generalized Gelfand-Graev characters
\par
3. Character formula for generalized Gelfand-Graev characters
\par
4. Shintani descent and almost characters
\par
5. Unipotently supported functions
\par
6. Character sheaves
\par
7. Lusztig's conjecture
\par
8. Parametrization of almost characters
\par
9. Proof of Theorem 8.6
\bigskip 
\section{Parametrization of irreducible characters}
\para{1.1.}
Let $k$ be an algebraic closure of a finite field $\Fq$ with 
ch $k = p$. In this and next section, we review the parametrization
of irreducible characters of $SL_n(\Fq)$ (or more generally, its 
Levi subgroups) following [S2].  we assume 
that $p$ is large enough so that the Dynkin-Kostant theory on 
Lie algebras can be 
applied.  For example, $p \ge 2n$ is enough for $G = SL_n$ (See the 
remark in 2.1).  
\par 
Let $\wt G = GL_{n_1} \times\cdots\times GL_{n_r}$.  We regard
$\wt G$ as a subgroup of $GL_n(k)$ with $n = \sum n_i$ in a natural 
way, and put $G = \wt G \cap SL_n$.  Thus $G$ is a Levi subgroup
of a parabolic subgroup of $SL_n$.
We consider a Frobenius map $F$ on $\wt G$ of the form 
$F = \f F_0$, where $F_0$ is the split Frobenius map on $\wt G$
with respect to the $\Fq$-structure, 
and $\f$ is a permutation of the factors in $\wt G$.
\par  
Put $\wt G_i = GL_{n_i}$.  We choose an $F_0$-split 
maximal torus $\wt T_i$ in $\wt G_i$ so that  
$\wt T = \wt T_1 \times \cdots \times \wt T_r$ is an 
$F$-stable maximal
torus in $\wt G$, which is maximally split with respect to $F$.
Let $W = N_{\td G}(\wt T)/\wt T$ be the Weyl group of $\wt G$.
Then $W \simeq W_1 \times \cdots \times W_r$ where
$W_i = N_{\wt G_i}(\wt T_i)/\wt T_i$ is the Weyl group of $\wt G_i$.
\par
Let $\wt G^* \simeq \wt G^*_1 \times \cdots \times \wt G_r^*$ be the
dual group of $\wt G$ over $\Fq$, and $G^*$ the dual group of $G$.  
We denote also by $F$ the corresponding Frobenius actions 
on $\wt G^*$ and $G^*$.  
The natural inclusion map $G \hra \wt G$ induces a map 
$\pi : \wt G^* \to G^*$, which is identified with the projection 
$\wt G^* \to \wt G^*/Z_1$, where $Z_1$ is the center of 
$GL_n$ under the identification of $\wt G^*$ with 
the subgroup of $GL_n^* = GL_n$.
Then $Z_{\wt G^*}$ and $Z_{G^*}$ are connected, and  
$Z_{G^*} = Z_{\wt G^*}/Z_1$.  We have natural inclusions and 
projections
\begin{align*} 
\wt G & \hla G \hla G\der = \wt G\der, \\
\wt G^* &\to G^* \to G^*/Z_{G^*} = \wt G^*/Z_{\wt G^*}, 
\end{align*}
where the dual group of $G\der$ (resp. $\wt G\der$) 
is identified with $G^*/Z_{G^*}$ (resp. $\wt G^*/Z_{\wt G^*}$). 
\par
We note that, for a connected reductive group $H$ defined over $\Fq$ 
with Frobenius map $F$, there exists an isomorphism of abelian groups
\begin{equation*}
\tag{1.1.1}
f: (Z^0_{H^*})^F \isom (H^F/H\der^F)\wg, 
\end{equation*}
where $H\der^F$ means $(H\der)^F$.
 Returning to the original setting, let 
$S = \wt T \cap \wt G\der$ be an $F$-stable  maximal torus of $\wt G\der$.
We have $\wt G = Z_{\wt G}\wt G\der$, $\wt T = Z_{\wt G}S$,  and 
$Z_{\wt G} \cap \wt G\der = Z_{\wt G} \cap S$ is finite.  
It follows that
$\wt G^F/\wt G\der^F \simeq \wt T^F/S^F$ and we have a natural 
inclusion map $Z_{\wt G}^FS^F/S^F \hra \wt T^F/S^F$, 
which induces a 
surjective map 
$(\wt G^F/\wt G\der^F)\wg \to (Z_{\wt G}^FS^F/S^F)\wg$.
Then we have the following lemma.
\begin{lem} 
Assume that $n_1 = \dots = n_r = d$, and put 
$(Z_{\wt G^*})_d = \{ z \in Z_{\wt G^*} \mid z^d = 1\}$. 
Then there exists an isomorphism 
$f_0: Z_{\wt G^*}^F/(Z_{\wt G^*})_d^F \to (\wt Z^{F}S^F/S^F)\wg$ 
which makes the following 
diagram commutative.
\begin{equation*}
\begin{CD}
Z_{\wt G^*}^F @>> f > (\wt G^F/\wt G\der^F)\wg \\
@VVV            @VVV     \\
Z_{\wt G^*}^F/(Z_{\wt G^*})_d^F  @>> f_0 > (Z_{\wt G}^FS^F/S^F)\wg,
\end{CD}
\tag{1.2.1}
\end{equation*}
where the vertical maps are natural surjections.
\end{lem}
\begin{proof}
By the isomorphism $f$, the subgroup $(Z_{\wt G^*})_d^F$ of 
$Z_{\wt G^*}^F$ 
is mapped onto the subgroup 
$A = \{ \th \in (\wt T^F/S^F)\wg \mid \th^d = 1\}$ of $(\wt T^F/S^F)\wg$.
We want to show that 
\begin{equation*}
\tag{1.2.2}
A = \{ \th \in (\wt T^F/S^F)\wg \mid \th|_{Z_{\wt G}^FS^F/S^F}= 1\}.
\end{equation*}
Let $B$ be the right hand side of (1.2.2).  Note that 
$\wt T = Z_{\wt G}S$, and $Z_{\wt G} \cap S$ consists 
of elements $z$ such that
$z^d = 1$ since $\wt G\der \simeq SL_d\times\cdots\times SL_d$.
It follows that for $zv \in (Z_{\wt G}S)^F$ with 
$z \in Z_{\wt G}, v \in S$,
we have $z^d \in Z_{\wt G}^F, v^d \in S^F$.  Hence for $\th \in B$, 
we have $\th^d(zv) = \th(z^dv^d) = 1$.  This implies that 
$\th \in A$, and we have $B \subseteq A$. 
Here $|A| = |(Z_{\wt G^*})^F_d|$.  On the other hand, 
$|B| = |\wt T^F|/|Z_{\wt G}^FS^F| = |Z_{\wt G}^F\cap S^F|$
since $|\wt T^F| = |Z_{\wt G}^F||S^F|$ (cf. [C, Prop. 3.3.7]). 
Under the identification $Z_{\wt G^*}^F \simeq Z_{\wt G}^F$, we have 
$(Z_{\wt G^*})^F_d \simeq Z_{\wt G}^F \cap S^F$.  It follows that 
$A = B$, and (1.2.2) follows. 
\par
Now we have natural isomorphisms
\begin{equation*}
Z_{\wt G^*}^F/(Z_{\wt G^*})^F_d 
    \simeq (\wt T^F/S^F)\wg/A \simeq (\wt T^F/S^F)\wg/B
             \simeq (Z_{\wt G}^FS^F/S^F)\wg, 
\end{equation*}
which makes the diagram (1.2.1) commutative.  This proves the lemma.
\end{proof} 
\para{1.3.}
We fix a dual torus 
$\wt T^*$ of $\wt T$ over $\Fq$ in $\wt G^*$.  Then the Weyl group
$N_{\td G^*}(\wt T^*)/\wt T^*$ may be identified with $W$.
For any semisimple element $\ds \in \wt T^*$ such that the conjugacy
class $\{ \ds\}$ of $\ds$ in $\wt G^*$ is $F$-stable, put
\begin{align*}
W_{\ds} &= \{ w \in W \mid w(\ds) = \ds\} \\
Z_{\ds} &= \{ w \in W \mid Fw(\ds) = \ds\}.
\end{align*}
We fix an $F$-stable  Borel subgroup $\wt B$ of $\wt G$ containing
$\wt T$.  Let $\vS$ (resp. $\vS^+$) be a root system (resp. positive
root system) with respect to the pair $(\wt B, \wt T)$.  
Then $Z_{\ds}$ may be written as $Z_{\ds} = w_1W_{\ds}$ for some $w_1 \in W$.
The element $w_1$ is determined uniquely by the condition that 
$w_1$ maps $\vS_{\ds}^+$ into $\vS^+$, where $\vS_{\ds}$ is the subroot system
of $\vS$ corresponding to $\ds$.
\par
Put $T = \wt T \cap G$ and $B = \wt B \cap G$.  
Then $T$ is the maximally split maximal torus of $G$, and 
$B$ is an $F$-stable Borel subgroup of $G$ containing $T$.  
We identify $W$ with the Weyl group $N_G(T)/T$ of $G$.
\para{1.4.}
The element $w_1$ induces an automorphism $\g: W_{\ds} \to W_{\ds}$ by 
$\g(w) = F(w_1ww_1\iv)$.
Let $W_{\ds}\lp\g\rp$ be the semidirect product of $W_{\ds}$ with the cyclic
group $\lp\g\rp$ generated by $\g$.
We denote by $(W_{\ds}\wg)^{\g}$ the set of $\g$-stable irreducible 
characters of $W_{\ds}$.
Each $E \in (W_{\ds}\wg)^{\g}$ is extendable to an irreducible character 
of $W_{\ds}\lp\g\rp$. We fix the preferred extension $\wt E$ of $E$ 
(see [L3, 17.2]).
\par
Let $\wt T_w$ be an $F$-stable maximal torus of $\wt G$ obtained from 
$\wt T$ by twisting by $w \in W$.
Then to any $w \in Z_{\ds}$, one can attach an irreducible character
$\th_w$ of $\wt T^{Fw} \simeq \wt T^F_w$ (see e.g., [S2, 2]).
Let $R^{\wt G}_{\wt T_w}(\th_w)$ be the Deligne-Lusztig character
of $\wt G^F$ associated to $\th_w \in (\wt T_w^F)\wg$.
For each $\wt E \in (W_{\ds}\wg)^{\g}$, put
\begin{equation*}
\wt\r_{\ds, E} = (-1)^{l(w_1)}|W_{\ds}|\iv
         \sum_{w \in W_{\ds}}\Tr(\g w, \wt E)
                  R^{\wt G}_{\wt T_{w_1w}}(\th_{w_1w}).
\end{equation*}
Then $\wt\r_{\ds, E}$ gives rise to an irreducible character of 
$\wt G^F$.  The set $\Irr \wt G^F$ is decomposed as 
\begin{equation*}
\Irr \wt G^F = \coprod_{\{\ds\}}\CE(\wt G^F, \{\ds\}),
\end{equation*}
where $\{\ds\}$ runs over all the $F$-stable semisimple 
classes in $\wt G^*$. 
The Lusztig series $\CE(\wt G^F, \{\ds\})$ associated to 
the $F$-stable class $\{\ds\}$ is given as 
\begin{equation*}
\CE(\wt G^F, \{\ds\}) = \{ \wt\r_{\ds,E} \mid E \in (W_{\ds}\wg)^{\g}\}.
\end{equation*}
\para{1.5.}
We now describe the irreducible characters of $G^F$ following [S2].
Let $\pi : \wt G^* \to G^*$ be as in 1.1. 
Let $T^* = \pi(\wt T^*)$ be the 
maximal torus of $G^*$.  Then $W = N_{\wt G^*}(\wt T^*)/\wt T^*$
is naturally identified with $N_{G^*}(T^*)/T^*$.
As in the case of $\wt G^F$, the set $\Irr G^F$ is partitioned as 
\begin{equation*}
\Irr G^F  = \coprod_{\{s\}}\CE(G^F, \{ s\}), 
\end{equation*}
where $\{ s\}$ runs over all the $F$-stable semisimple classes in $G^*$.
We fix $ s \in T^*$ for a given $F$-stable class 
$\{ s\} \subset G^*$.  There exists $\ds \in \wt T^*$ such that
$\pi(\ds) = s$ and that the class $\{ \ds\}$ is $F$-stable.
One can find $w_1 \in Z_{\ds}$ and an isomorphism 
$\g = Fw_1: W_{\ds} \to W_{\ds}$ as in 1.4. 
\par
Put $W_{s} = \{ w \in W \mid w(s) = s\}$.  Then $W_{\ds}$ is 
naturally regarded as a subgroup of $W_s$, and we have
$W_{s} \simeq W_{\ds}\rtimes \Om_s$, where $\Om_s$ is 
a cyclic group isomorphic to $Z_{G^*}(s)/Z^0_{G^*}(s)$.
$W_{\ds}$ is characterized as the largest reflection subgroup 
of $W_s$, and sometimes we denote it by $W_s^0$.
Let $F' = F\dw_1$, where $\dw_1$ is the representative of $w_1$ in 
$N_{\wt G^*}(\wt T^*)$.  Then $\pi$ is $F'$-equivariant, and 
$s \in (T^*)^{F'}$.  So $F'$ acts naturally on $W_s$, leaving
$W_{\ds}$ and $\Om_s$ invariant.
We consider the set $\Om_s\ssim_{F'}$ of $F'$-twisted 
classes in $\Om_s$.  Since $\Om_s$ is abelian, $\Om_s\ssim_{F'}$ 
is identified with $(\Om_s)_{F'}$, the largest 
quotient on which $F'$ acts trivially.   
For each $x \in (\Om_s)_{F'}$, one can choose $\ds_x \in (\wt T^*)^{xF'}$
such that $\pi(\ds_x) = s$, and obtain an isomorphism
$\g_x = xF': W_{\ds} \to W_{\ds}$. 
To each $\g_x$-stable irreducible character $E$ of $W_{\ds}$, one can
attach the irreducible character $\wt\r_{\ds_x, E}$ of $\wt G^F$ as before.
We denote by $\CT_{\ds_x, E}$ the set of irreducible characters of $G^F$
occurring in the restriction of $\wt\r_{\ds_x, E}$ on $G^F$. 
Then by [L5], we can decompose $\CE(G^F, \{ s\})$ as 
\begin{equation*}
\tag{1.5.1}
\CE(G^F, \{ s\}) = \coprod_{(x, E)}\CT_{\ds_x, E},
\end{equation*}
where the pair $(x, E)$ runs over all $x \in (\Om_s)_{F'}$ and
$E \in (W_{\ds}\wg)^{\g_x}/\Om_s^{F'}$ (the set of $\Om_s^{F'}$-orbits
in $(W_{\ds}\wg)^{\g_x}$).
\para{1.6.}
Following [S2], we shall modify the partition in (1.5.1). 
For $E \in W_{\ds}\wg$, let $\Om_{s,E}$ be the stabilizer of $E$ in
$\Om_s$. (In [S2], the notation $\Om_s(E)$ is used instead of 
$\Om_{s,E}$).  
If the $\Om_s$-orbit of $E$ in $W_{\ds}\wg$ is $F'$-stable, then
$\Om_{s,E}$ is $F'$-stable, and one can consider the largest quotient
$(\Om_{s,E})_{F'}$ as before.
If we put $\wt\Om_{s,E} = \{ x \in \Om_s \mid {}^{xF'}E = E\}$, 
then $\wt\Om_{s,E} \ne \emptyset$, and one can write it as
$\wt\Om_{s,E} = \Om_{s,E}a_E$ for some $a_E \in \Om_s$.
It follows that $\wt\Om_{s,E}\ssim_{F'}$ can be identified with
the set $(\Om_{s,E}\ssim_{F'})a_E$ and with $(\Om_{s,E})_{F'}a_E$.
We denote this set by $(\wt\Om_{s,E})_{F'}$.
By (4.4.2) in [S2], we have the following natural bijection
\begin{equation*}
\tag{1.6.1}
\coprod_{E \in (W_{\ds}\wg/\Om_s)^{F'}}(\wt\Om_{s,E})_{F'} \simeq
          \coprod_{x \in (\Om_s)_{F'}}(W_{\ds}\wg)^{\g_x}/\Om_s^{F'},
\end{equation*}   
where $(W_{\ds}\wg/\Om_s)^{F'}$ denotes the set of $F'$-stable 
$\Om_s$-orbits in $W_{\ds}\wg$.
(In [S2], $\Om_{s,E}$ is used instead of $\wt\Om_{s,E}$.
This is justified since we have a bijection 
$\wt\Om_{s,E} \simeq \Om_{s,E}$.  However this bijection depends on 
the choice of $a_E$, and so the form as in (1.6.1) is more convenient 
for our later purpose.)
\par
Let $E \in (W_{\ds}\wg/\Om_s)^{F'}$, i.e., the $\Om_s$-orbit
of $E$ is $F'$-stable.  Then for each 
$y \in (\wt\Om_{s,E})_{F'}$, 
one can associate the pair $(x, E')$, where $x \in (\Om_s)_{F'}$ and 
$E' \in (W_{\ds}\wg)^{\g_x}/\Om_s^{F'}$, by (1.6.1).
We denote by $\CT_{s, E}$ the union of various 
$\CT_{\ds_x, E'}$ where $(x, E')$ runs over all the pairs in the image
of $(\Om_{s,E})_{F'}$ under the bijection in (1.6.1).
Thus we can rewrite (1.5.1) as 
\begin{equation*}
\tag{1.6.2}
\CE(G^F, \{ s\}) = \coprod_{E \in (W_{\ds}\wg/\Om_s)^{F'}}\CT_{ s, E}.
\end{equation*}
For a pair $(s, E)$ with $E \in (W_{\ds}\wg/\Om_s)^{F'}$, put 
\begin{equation*}
\tag{1.6.3}
\ol\CM_{s, E} = (\Om_{s,E}^{F'})\wg\times(\wt\Om_{s,E})_{F'},
\end{equation*}
 where 
$\Om_{s,E}^{F'}$ is the $F'$-fixed subgroup of $\Om_{s,E}$. 
\par
It is known by [L5] that $\CT_{\ds_x, E}$ is in bijection with the set
$(\Om_{s,E}^{xF'})\wg =  (\Om_{s,E}^{F'})\wg$.  Hence by (1.5.1), 
$\CE(G^F, \{ s\})$ 
is parametrized by various $(\Om_{s,E}^{F'})\wg$.
However, this parametrization of $\CT_{\ds_x,E}$ is not canonical.  
It depends on the
choice of an irreducible character $\r_0 \in\Irr G^F$ occurring in the
decomposition of $\wt\r_{\ds_x, E}$.   In [S2], a 
bijective correspondence $\CT_{ s, E} \lra \ol\CM_{s, E}$ 
was constructed by making use of generalized
Gelfand-Graev characters. This bijection is determined uniquely
once we fix a representative $u \in C^F$ for
each $F$-stable unipotent class $C$ in $G$.   
Thus we have a parametrization of $\CE(G^F,\{ s\})$ as
\begin{equation*}
\tag{1.6.4}
\CE(G^F,\{ s\}) = \coprod_{E \in (W_{\ds}\wg/\Om_s)^{F'}}
                         \ol\CM_{ s,E}.
\end{equation*}
In the next section, we shall explain this parametrization in 
details.
\par\medskip
\section{Generalized Gelfand-Graev characters}
\para{2.1.}
In order to explain the bijection 
$\CT_{s, E} \lra \ol\CM_{s, E}$, we shall review results
on generalized Gelfand-Graev characters following [S2].
(Although no restriction on $p$ was assumed in [S2], 
this must be changed.  In fact Kawanaka's construction of 
generalized Gelfand-Graev characters of $GL_n$ or $SL_n$ 
requires no assumption on $p$.  However, our construction 
([S2, 2.3]) relies on the Dynkin-Kostant theory, which requires
that $p$ is not too small).
We also prove a character formula for modified generalized 
Gelfand-Graev characters, which will play an essential role
in later sections.
\par
Let $\Fg$ be the Lie algebra of $G$ with Frobenius map $F = \f F_0$.
Let $G\uni$ (resp. $\Fg\nil$) be the set of unipotent elements in $G$
(resp. nilpotent elements in $\Fg$).  We have a bijection 
$\log: G\uni \to \Fg\nil, v \mapsto v-1$.
Let $\CO$ be an $F$-stable nilpotent orbit in $\Fg$, and choose a
representative $N \in \Fg^F$. Correspondingly, we consider an $F$-stable
unipotent class $C$ containing $u = \log\iv N \in C^F$.
By Dynkin-Kostant theory, there exists a natural grading 
$\Fg = \oplus_{i \in \ZZ}\, \Fg_i$ associated to $N$.  Let 
$\Fu_i = \oplus_{j \ge i}\Fg_j$ for $j \ge 1$.  Then $\Fu_i$ is a
nilpotent subalgebra of $\Fg$, and there exists a connected unipotent
subgroup $U_i$ of $G$ defined over $\Fq$ such that $\log(U_i) = \Fu_i$.
Also one can find a parabolic subgroup $P =P_N$ and its Levi subgroup
$L = L_N$ such that $P = LU_1$, 
where $L$ is an $F$-stable Levi subgroup of $P$ with $\Lie(L) = \Fg_0$
and $U_1$ is the unipotent radical of $P$.  Moreover, we have 
$N \in \Fg_2$.
\par
Let $N^* \in \Fg^F_{-2}$
be the  element such that $\{N, N^*, H\}$ gives a TDS triple for some 
semisimple element $H \in \Fg_0$. 
We define a linear map $\la: \Fu_1 \to k$ by  
$\la(x) = \lp N^*, x\rp$ , where $\lp\ ,\ \rp$ is a fixed $G$-invariant
non-degenerate bilinear form on $\Fg$.  
It is known that the map $(x,y) \mapsto \la([x,y])$ defines a
symplectic form on $\Fg_1$, and according to [S2, 2.3] one can find an
$F$-stable Lagrangian subspace $\Fs$ of $\Fg_1$ satisfying 
the following properties.
Put $\Fu = \Fs + \Fu_2$ ($\Fu = \Fu_{1.5}$ in the notation in [S2]).  
Then $\Fu$ is a subalgebra of 
$\Fu_1$, and we obtain an $F$-stable closed subgroup $U$ of 
$U_1$ such that $\log(U) = \Fu$.  $U$ is a normal
subgroup of $U_1$.  Moreover, 
$U$ is stable by the conjugation action of $L$. (Note that 
the last property does not hold for a general Lagrangian subspace
(see [S2, 2.6]).) 
Now the map $\la\circ\log: U \to k$ turns out to be an
$F$-stable homomorphism from $U$ to $k$.  Thus we obtain a 
linear character $\vL_N$ on $U^F$ by 
$\vL_N = \p\circ\la\circ\log$, where $\p: \Fq \to \Ql^*$ is a fixed 
non-trivial additive character of $\Fq$.  
The generalized Gelfand-Graev character $\vG_N$ on $G^F$ 
associated to $N$ is defined
as $\vG_N = \Ind_{U^F}^{G^F}\vL_N$.
\par
Following [K3], we construct modified generalized Gelfand-Graev 
characters.
Let $A_{\la} = Z_L(\la)/Z_L^0(\la)$ for $\la: \Fu \to k$.
Then by [S2, 2.7] we have 
\begin{equation*}
\tag{2.1.1}
A_{\la} \simeq A_G(N) = Z_G(N)/Z_G^0(N).
\end{equation*}
  Since the latter group is abelian, 
$A_{\la}$ is an abelian group.  Moreover, we have a surjective
map $Z_G \to A_{\la}$.  Put 
\begin{equation*}
\tag{2.1.2}
\ol\CM = (A_{\la})_F\times (A_{\la}^F)\wg.
\end{equation*}  
For each $(c,\xi) \in \ol\CM$,   one can construct a character 
$\vG_{c,\xi}$ on $G^F$ as follows;
for $c \in A_{\la}$, we choose a representative 
$\dc \in Z_L(\la)$.  Then we find $\a_c \in L$ such that 
$\a_c\iv F(\a_c) = \dc$. Let us define a linear map 
$\la_c: \Fu \to k$ by $\la_c = \la\circ\Ad\a_c\iv$, where $\Ad$ is the
adjoint action of $L$ on $\Fu$.   
We define a linear character $\vL_c$ on $U^F$ by 
$\vL_c = \p\circ\la_c\circ\log$.
If we notice that $Z_L(\la_c)^F$ coincides with $Z_{L^F}(\vL_c)$, the
linear character $\vL_c$ can be extended to the character 
of $Z_L(\la_c)^FU^F$ so that it is trivial on $Z_L(\la_c)^F$, which we
denote by the same symbol $\vL_c$.
On the other hand, by the isomorphism
\begin{equation*}
Z_L(\la_c)^F/Z_L^0(\la_c)^F \simeq 
        Z_L(\la)^{\dc F}/Z_L^0(\la)^{\dc F} \simeq  
             A_{\la}^{\dc F} = A_{\la}^F,
\end{equation*}
the linear character $\xi \in (A_{\la}^F)\wg$ determines 
a linear character
$\xi\nat$ on $Z_L(\la_c)^F$ which is trivial on $Z_L^0(\la_c)^F$.
\par
Let $\wt P$ and $\wt L$ be the parabolic subgroup of 
$\wt G$ and its Levi subgroup associated to 
$N \in \wt\Fg = \Lie \wt G$.  Then we have $P = \wt P \cap G$ and 
$L = \wt L \cap G$, and so $Z_L(\la) \subset Z_{\wt L}(\la)$.
Let us take a linear character $\th$ of $Z_L(\la)^F$ of the following type.
\par\medskip\noindent
(2.1.3) \ $\th$ is the restriction to $Z_L(\la)^F$ of 
a linear character of $Z_{\wt L}(\la)^F$ which is trivial on 
$(Z_{\wt L}(\la)\der)^F$.
\par\medskip
Since $\dc \in Z_L(\la)$, we have 
$(Z_{\wt L}(\la)/Z_{\wt L}(\la)\der)^F 
       \simeq (Z_{\wt L}(\la)/Z_{\wt L}(\la)\der)^{\dc F}$.
It follows that $\th$ is regarded as a linear character of 
$Z_L(\la)^{\dc F}$, and it determines a linear character of 
$Z_L(\la_c)^F$ via the isomorphism  
$\ad \a\iv : Z_L(\la_c)^F \simeq Z_L(\la)^{\dc F}$, which we 
denote also by $\th$.
Then $\th\xi\nat$ gives rise to a character on $Z_L(\la_c)^FU^F$
under the surjective homomorphism 
$Z_L(\la_c)^FU^F \to Z_L(\la_c)^F$, which we denote also by
$\th\xi\nat$.    Under this setting, we define a modified 
generalized Gelfand-Graev character $\vG_{c,\xi,\th}$ by
\begin{equation*}
\tag{2.1.4}
\vG_{c,\xi,\th} = \Ind_{Z_L(\la_c)^FU^F}^{G^F}
            \bigl(\th\xi\nat\otimes\vL_c\bigr).
\end{equation*}
In the case where $\th = 1$, we simply write $\vG_{c,\xi,\th}$ as
$\vG_{c,\xi}$.
\par
For later use, we also define a generalized Gelfand-Graev character
$\vG_c$ associated to $c \in A_{\la}$ by 
$\vG_c = \Ind_{U^F}^{G^F}\vL_c$.  Under the isomorphism 
$A_{\la} \simeq A_G(N)$ (2.1.1), one can construct a nilpotent 
element $N_c \in \Fg_2^F$ twisted by $c$.  $\vG_c$ is nothing but
the generalized Gelfand-Graev character $\vG_{N_c}$ associated to $N_c$.
We remark that $\vG_{c,\xi}$ occurs as a direct summand of $\vG_c$.
\para{2.2.}
We choose $m$ large enough so that $F^m$ acts trivially on 
$A_{\la}$.  Replacing $F$ by $F^m$, we have a modified 
generalized Gelfand-Graev character $\vG_{c,\xi,\th}^{(m)}$ on 
$G^{F^m}$.    Now the parameter set $\ol\CM$ is replaced by 
$A_{\la} \times A_{\la}\wg$.
We denote by $\CM$ the subset of 
$A_{\la} \times A_{\la}\wg$ defined by 
\begin{equation*}
\CM = A_{\la}^F \times (A_{\la}\wg)^F,
\end{equation*}
where $(A_{\la}\wg)^F$ is the set of $F$-stable irreducible 
characters of $A_{\la}$.
We now construct, for a certain linear character $\th$ of 
$Z_L(\la_c)^{F^m}$, and 
for each $(c, \xi) \in \CM$, 
an $F$-stable modified generalized Gelfand-Graev character 
$\vG^{(m)}_{c,\xi,\th}$, and its extension to 
$G^{F^m}\lp\s\rp$, where $\s = F|_{G^{F^m}}$. The case where
$\th = 1$ is discussed in [S2, 1.8].    
For $c \in A_{\la}^F$, we choose $\dc \in L^F$. 
We construct the linear character 
$\vL^{(m)}_c$ of $U^{F^m}$ as in 2.1, i.e.,   
we choose $\b_c \in L$
such that $\b_c\iv F^m(\b_c) = \dc$, 
and define $\la_c$
by $\la_c = \la\circ\Ad\b_c\iv$, and put 
$\vL_c^{(m)} = \p_m\circ\la_c\circ\log$, where 
$\p_m = \p\circ\Tr_{{\mathbf F}_{q^m}/\Fq}$.  Put 
$\hat c = \b_cF(\b_c\iv) \in L^{F^m}$. Then  
$\vL_c^{(m)}$ turns out to be $\hat cF$-stable.
(Note that it is possible to choose $\dc \in T^F$.
Then we can choose $\b_c \in T$.)
\par
On the other hand, it can be checked that $\hat cF$ acts
on $Z_L(\la_c)$ commuting with $F^m$, and that
under the isomorphism 
\begin{equation*}
\ad \b_c\iv : Z_L(\la_c)^{F^m}/Z_L^0(\la_c)^{F^m} 
  \simeq Z_L(\la)^{\dc F^m}/Z_L^0(\la)^{\dc F^m}
  \simeq A_{\la},
\end{equation*}
the action of $\hat cF$ on $Z_L(\la_c)^{F^m}$
is transferred to the action of $F$ on $A_{\la}$. 
Hence if we take $\xi \in (A_{\la}\wg)^F$, it produces
a $\hat cF$-stable linear character $\xi\nat$ 
on $Z_L(\la_c)^{F^m}$.
\par
Furthermore, we take a linear character $\th$ of $Z_L(\la)^{F^m}$
of the following type.
\par\medskip\noindent
(2.2.1) \ $\th$ is the restriction to $Z_L(\la)^{F^m}$ 
of an $F$-stable  linear character of $Z_{\wt L}(\la)^{F^m}$ 
as in (2.1.3) by replacing $F$ by $F^m$.
\par\medskip
Then $\th$ is regarded as an  $F$-stable 
linear character of $Z_L(\la)^{\dc F^m}$.  
It follows, under the isomorphsim 
$\ad \b_c\iv: Z_L(\la_c)^{F^m} \simeq Z_L(\la)^{\dc F^m}$, that  
$\th$ determines a $\hat c F$-stable linear character of 
$Z_L(\la_c)^{F^m}$, which we denote also by $\th$.
Thus  
 $\th\xi\nat\otimes \vL^{(m)}_c$ is $\hat cF$-stable 
for $(c, \xi) \in \CM$, and we conclude that 
$\vG^{(m)}_{c,\xi,\th}$ is $F$-stable.
\par
Put $\hat c_0 = (\hat c\s)^m \in L^{F^m}$.
Then we have 
\begin{equation*} 
\hat c_0 = \b_cF^m(\b_c\iv) = \dc\iv
\end{equation*}
since $\b_c$ and $F^m(\b_c)$ commutes. 
We note that 
$\hat c_0 \in Z_L(\la_c)^{F^m} = Z_L(\vL_c^{(m)})^{F^m}$.
In fact, since $\vL_c^{(m)}$ is $\hat cF$-stable, 
it is stable by $(\hat c\s)^{m} = \hat c_0$.
Put $M_c = Z_L(\la_c)^{F^m}$ and $M_c^0 = Z_L^0(\la_c)^{F^m}$.
We consider a subgroup $M_cU^{F^m}\lp\hat c\s\rp$ of 
$G^{F^m}\lp\s\rp$ generated by $M_cU^{F^m}$ and $\hat c\s$. 
Since $\th\xi\nat \in M_c\wg$ is $\hat c\s$-stable, and 
$(\hat c\s)^m = \hat c_0 \in M_c$, $\th\xi\nat$ may be extended 
to a linear character $\wt{\th\xi}\nat$ of $M_c\lp\hat c\s\rp$
in $m$ distinct way.  
The extension $\wt{\th\xi}\nat$ is determined by the value 
$\wt{\th\xi}\nat(\hat c\s) = \mu_{c,\th\xi}$, where $\mu_{c,\th\xi}$ is
any $m$-th root of $\th\xi\nat(\dc\iv) = \th(\dc\iv)\xi(c\iv)$.
\par
We fix an extension $\wt{\th\xi}\nat$ of $\th\xi\nat$ to 
$M_c\lp\hat c\s\rp$.
Since $M_cU^{F^m}\lp\hat c\s\rp$ is the semidirect product
of $M_c\lp\hat c\s\rp$ with $U^{F^m}$, $\wt{\th\xi}\nat$ may be
regarded as a character of $M_cU^{F^m}\lp\hat c\s\rp$.
On the other hand, since $\vL_c^{(m)}$ is $\hat c\s$-stable,
it can be extended to a linear character $\wt\vL_c^{(m)}$ on 
$M_cU^{F^m}\lp\hat c\s\rp$ by $\wt\vL_c^{(m)}(\hat c\s) = 1$.
Thus we have a character $\wt{\th\xi}\nat\otimes\wt\vL_c^{(m)}$ of
$M_cU^{F^m}\lp\hat c\s\rp$ which is an extension of 
$\th\xi\nat\otimes\vL_c^{(m)}$ on $M_cU^{F^m}$.
We put
\begin{equation*}
\wt\vG_{c,\xi,\th}^{(m)} = 
   \Ind_{M_cU^{F^m}\lp\hat c\s\rp}^{G^{F^m}\lp\s\rp}
                            (\wt{\th\xi}\nat\otimes\wt\vL_c^{(m)}).
\end{equation*}
Then $\wt\vG_{c,\xi,\th}^{(m)}$ gives rise to 
an extension of $\vG^{(m)}_{c,\xi,\th}$
to $G^{F^m}\lp\s\rp$.  Note that 
$\mu_{c,\th\xi}\iv \wt\vG^{(m)}_{c,\xi,\th}|_{G^{F^m}\s}$ depends only on
the choice of $(c, \xi)$ and $\th$.
\para{2.3.}
In [L1], Lusztig defined, for a connected reductive group $H$ with 
connected center, a map from the set of irreducible characters of
$H^F$ to the set of $F$-stable unipotent classes in $H$. 
It is shown in [L7] that this map coincides with the map defined 
by Kawanaka [K1, K2, K3] in terms of generalized
Gelfand-Graev characters.  We denote by $C_{\r}$ 
(resp. $\CO_{\r}$ ) the unipotent
class in $H$ (resp. the nilpotent orbit in $\Lie H$) 
corresponding to $\r \in \Irr H$ under this map. 
We call $C_{\r}$ the unipotent class associated to $\r$ 
(the wave front set associated to $\r$ in the sense of 
Kawanaka).   
\par
In the case of $\wt G$, this map is given as follows.
Let $\wt\r = \wt\r_{\ds, E} \in \CE(\wt G^F, \{\ds\})$.  Put 
$E' = E \otimes\ve$ for the sign character $\ve$ of $W_{\ds}$.
Then $\Ind_{W_{\ds}}^WE'$ contains a unique irreducible character 
$\wh E$ of $W$ such that $b_{E'} = b_{\wh E}$. $C_{\wt\r}$
is defined as the unipotent class in $\wt G$ corresponding to $\wh E$ 
under the Springer correspondence. 
More precisely, we have the following. Assume that $\wt G = GL_n$.  
Then $W_{\ds}$ is a product
of various symmetric groups. Accordingly, $E \in W_{\ds}\wg$ 
is parametrized by a multipartition 
$\b = (\b_1, \dots, \b_k)$ of $n$.  By mixing and 
rearranging the parts in $\b_1, \dots, \b_k$, 
we regard $\b$ as a partition of $n$ which we denote by $\bar\b$. 
Let $\bar\b^*$ be the
dual partition of $\bar\b$.  Then $C_{\wt\r}$ is the unipotent class 
in $\wt G$ corresponding to $\bar\b^*$ through Jordan's normal form. 
For general $\wt G$, the description of $C_{\wt\r}$ is reduced to the 
case of $GL_n$ through the decomposition 
$\wt G = GL_{n_1}\times\cdots\times GL_{n_r}$. 
\para{2.4.}
Let $\wt\Fg = \Lie \wt G$. 
For a nilpotent element $N \in \wt\Fg^F$, we denote by $\CO_N$ the 
nilpotent orbit in  $\wt\Fg$ containing $N$.  Let $P = LU_1$ be as in 2.1, 
and let $\wt L \subset \wt P$ be as before. 
For each irreducible character $\th'$ of
$Z_{\wt L}(\la)^F$, the modified generalized Gelfand-Grave
character $\wt\vG_{N,\th'}$ is defined as 
$\wt\vG_{N,\th'} = 
  \Ind_{Z_{\wt L}(\la)^FU^F}^{\wt G^F}(\th'\otimes\vL_N)$. 
 
\par
Let $\wt\r = \wt\r_{\ds, E}$ be an irreducible character of $\wt G^F$,
such that $\CO_{\wt\r} = \CO_N$. 
Then it is known by 
[S2, Prop. 2.14] that there exists a unique linear
character $\vf$ of $Z_{\wt L}(\la)^F$ such that  
\begin{equation*}
\lp\wt\vG_{N,\th'}, \wt\r\rp_{\wt G^F} = \begin{cases}
                    1  &\quad\text{ if } \th' = \vf, \\
                    0  &\quad\text{ if }  \th' \ne \vf.
                                   \end{cases}
\end{equation*}
We denote by $\vD(\wt\r)$ the linear character $\vf$ determined as above.
We note that 
\par\medskip\noindent
(2.4.1)\ $\vD(\wt\r)$ is a linear character of $Z_{\wt L}(\la)^F$
which is trivial on $(Z_{\wt L}(\la)\der)^F$.
\par\medskip\noindent
In fact, by [S2, 2.13], $\vf = \vD(\wt\r)$ is determined in the following way.
There exists an $F$-stable Levi subgroup $\wt M$ of a parabolic subgroup
of $\wt G$ such that $Z_{\wt L}(\la) \subset \wt M$ and that 
$\ds \in Z_{\wt M^*}$,
where $\wt M^* \subset \wt G^*$ is the dual group of $\wt M$.  We choose an 
integer $m >0$ such that $\ds \in Z_{\wt M^*}^{F^m}$, and let 
$\wh\vf$ be a linear character of $Z_{\wt L}(\la)^{F^m}$ obtained by
restricting the linear character of $\wt M^{F^m}$ corresponding to $\ds$.
Then $\wh\vf$ is $F$-stable, and the Shintani descent 
$Sh_{F^m/F}(\wh\vf)$ coincides with $\vD(\wt\r)$.  Since the linear
character of $\wt M^{F^m}$ corresponding to $\ds$ has a trivial
restriction on $\wt M\der^{F^m}$, we see that $\wh\vf$ has a trivial
restriction on $(Z_{\wt L}(\la)\der)^{F^m}$, and so
$\vf$ is trivial on $(Z_{\wt L}(\la)\der)^F$.   This shows (2.4.1).
\par
In view of (2.4.1), the restriction $\th$ of $\vD(\wt\r)$ to $Z_L(\la)^F$
satisfies the property in (2.1.3).  Hnece we can consider 
$\vG_{c,\xi,\th}$ as in 2.1, which tursn out to be a direct 
summand of $\wt\vG_{N,\th'}|_{G^F}$.  
\para{2.5.}
Let $(s, E)$ be as in 1.6. Then 
$\ds \in \wt G^* = GL_{n_1}\times \cdots\times GL_{n_r}$ is written as 
$\ds = (\ds_1, \dots, \ds_r)$ with $\ds_i \in GL_{n_i}$, and we have 
$W_{\ds} = W_{1, \ds_1}\times\cdots\times W_{r, \ds_r}$. 
We now consider the following special setting for the pair $(s, E)$.
\par\medskip\noindent
(2.5.1) \ Let $t$ be a common divisor of $n_1, \dots, n_r$ which is 
prime to $p$.  We have $\Om_s \simeq \lp w_0\rp$, where
$w_0 \in W_{s}$ is an element of order $t$ permuting the factors of 
$W_{i,\ds_i}$ transitively, and  
$W_{i,\ds_i}$ is isomorphic to 
$\FS_{b_i} \times\cdots\times \FS_{b_i}$ ($t$ times) with 
$b_i = n_i/t$. Moreover, $E \in (W_{\ds}\wg)^{F'}$ is of the form
\begin{equation*}
E = E_1\boxtimes\cdots\boxtimes E_r \quad\text{where}\quad
E_i \simeq E_i' \boxtimes\cdots\boxtimes E_i' \in W_{i,\ds_i}\wg
\text{ with } E_i' \in \FS_{b_i}\wg.
\end{equation*}
\par
Assume that the pair $(s, E)$ satisfies the condition (2.5.1). 
Then $E$ is $\Om_s$-stable, and in particular, $\g_x$-stable for
$x \in \Om_s$. 
Since $E$ is $\Om_s$-stable, we have $E \in (W_{\ds}\wg/\Om_s)^{F'}$
and $(\wt\Om_{s,E})_{F'} = (\Om_s)_{F'}$ with $a_E = 1$.  Hence the
bijection in (1.6.1) leaves the pair $(x,E)$ invariant for 
$x \in (\Om_s)_{F'}$.  
It follows that the set $\CT_{s, E}$ coincides with the disjoint union of 
$\CT_{\ds_x, E}$ for $x \in (\Om_s)_{F'}$.
\par
Let $\wt\r = \wt\r_{\ds,E}$.  It is known that $\wt\r|_{G^F}$ is multiplicity 
free.  Let $\CT_{\wt\r} = \CT_{\ds,E}$ be as in 1.5.  Then
$\CT_{\wt\r}$ consists of $t'$ elements, where $t'$ is the order of
$\Om_{s,E}^{F'} = \Om_s^{F'}$.
It follows from [S2] that $(A_{\la})_F$ acts transitively on the set
$\CT_{\wt\r}$.  Thus there exists a quotient $(A_{\la})'_F$ of 
$(A_{\la})_F$ such that $(A_{\la})'_F$ is in bijection with 
$\CT_{\wt\r}$. $(A_{\la})_F'$ can be written also as $(\bar A_{\la})_F$ 
with  some quotient $\bar A_{\la}$ of $A_{\la}$, where 
$\bar A_{\la}$ is a cyclic group of order $t$ (see [S2, 2.19]). 
It can be checked from the proof in [S2, 2.21] that the map 
$A_{\la}^F \to \bar A_{\la}^F$ is surjective. 
Let us define a set $\ol\CM_{s, N}$ by 
\begin{equation*}
\tag{2.5.2}
\ol\CM_{s, N} = (\bar A_{\la})_F \times (\bar A_{\la}^F)\wg.
\end{equation*}
The set $(\bar A_{\la}^F)\wg$ is regarded as a subset of 
$(A_{\la}^F)\wg$ through the map $A_{\la}^F \to \bar A_{\la}^F$.
Also we have a surjective map $(A_{\la})_F \to (\bar A_{\la})_F$.
We define a subset $\ol\CM_0$ of $\ol\CM$ by 
$\ol\CM_0 = (A_{\la})_F \times (\bar A_{\la}^F)\wg$.
Thus we have a natural surjective map 
$f: \ol\CM_0 \to \ol\CM_{s, N}$.
The following result, which gives a parametrization of $\CT_{s, E}$
in terms of generalized Gelfand-Graev characters, is a generalization
of the results 2.16 and 2.21 in [S2].  The proof is done in a similar 
way as in [S2]. 
\begin{thm} 
Assume that the pair $(s, E)$ satisfies (2.5.1).  Let 
$\wt\r = \wt\r_{\ds, E} \in \Irr \wt G^F$. Let $\CO_N$ be 
the nilpotent orbit
in $\wt\Fg$ containing $N$.  Let $\th$ be a linear 
character of $Z_L(\la)^F$ as in (2.1.3), and $\th_0$ the restriction
of $\th$ to $Z_L^0(\la)^{F}$. 
Then for each pair 
$(c,\xi)\in \ol\CM $, the following holds.
\begin{enumerate}
\item
$\lp\vG_{c,\xi,\th}, \wt\r|_{G^F}\rp_{G^F} = 0$ unless 
                      $\CO_N \subseteq \ol\CO_{\wt\r}$.
\item
Assume that $\CO_{\wt\r} = \CO_N$, and let $\vD(\wt\r)$ be as in 2.4.
\begin{enumerate}
\item
If $\vD(\wt\r)|_{Z_L^0(\la)^F} \ne \th_0$, then  
$\lp\vG_{c,\xi,\th}, \wt\r|_{G^F}\rp_{G^F} = 0$.
\item
If $\vD(\wt\r)|_{Z_L^0(\la)^F} = \th_0$, then 
there exists a bijection
 $\CT_{s, E} \lra \ol\CM_{s, N}$ satisfying the following;
Let $\r_{c,\xi} \in \CT_{s, E}$ be the character corresponding
to $(c,\xi) \in \ol\CM_{s, E}$.
 For each pair $(c',\xi') \in \ol\CM_0$ we have 
\begin{equation*}
\lp\vG_{c',\xi',\th}, \r_{c,\xi}\rp_{G^F} = \begin{cases}
             1  &\quad\text{ if } f((c',\xi')) = (c,\xi), \\
             0  &\quad\text{ if } f((c',\xi')) \ne (c,\xi).
                                      \end{cases} 
\end{equation*}
We have $\lp\vG_{c',\xi',\th}, \r_1\rp_{G^F} = 0$ for any 
$\r_1 \in \CT_{s, E}$, if the pair $(c',\xi') \in \ol\CM$ is 
not contained in $\ol\CM_0$.
\par
Furthermore, $\vD(\wt\r)|_{Z_L(\la)^F}$ 
is expressed as $\th\xi\nat_1$ for a character 
$\xi_1$ of $A_{\la}^F$. Then 
$\xi_1$ is contained in $(\bar A_{\la}^F)\wg$, 
and we have 
\begin{equation*}
\wt\r|_{G^F} = \sum_{c \in (\bar A_{\la})_F}\r_{c, \xi_1}.
\end{equation*}
\end{enumerate}
\end{enumerate}
\end{thm}
\para{2.7.}
The above parametrization of $\CT_{s, E}$ is also interpreted
in terms of (not modified) generalized Gelfand-Graev characters of 
$G^F$ as follows.
\par\medskip\noindent
(2.7.1) \ For each $x \in (\Om_s)_{F'}$ and $c \in (A_{\la})_F$, 
we have 
$\lp\vG_c, \wt\r_{\ds_x, E}|_{G^F}\rp_{G^F} = 1$, i.e., there exists a
unique irreducible character of $G^F$ which occurs both in the
decomposition of $\wt\r_{\ds_x, E}|_{G^F}$ and of $\vG_c$.  Under the 
parametrization in Theorem 2.6, this character
is given by $\r_{c, \xi_x}$ for some 
$\xi_x \in (\bar A_{\la}^F)\wg$. 
In particular, we have
\begin{equation*}
\wt\r_{\ds_x, E}|_{G^F} = \sum_{c \in (\bar A_{\la})_F}\r_{c,\xi_x}.
\end{equation*}
\par\medskip
By using (2.7.1) we can identify 
$\ol\CM_{s, N}$ with 
$\ol\CM_{s, E}$ in (1.6.3).  Note that in this case, 
$\ol\CM_{s, E}$ is nothing but the set 
$(\Om_s^{F'})\wg \times (\Om_s)_{F'}$.
Also note that the map $x \mapsto \xi_x$ gives a bijection 
$h: (\Om_s)_{F'} \to (\bar A_{\la}^F)\wg$, where 
$\xi_x$ is given by $\vD(\wt\r_{\ds_x,E})|_{Z_L(\la)^F} = \th\xi_x$.
By the discussion in 2.4, we can choose $\th$ such that  
$\vD(\wt\r_{\ds,E})|_{Z_L(\la)^F} = \th$.
Then $\xi_1$ (the case where $x = 1$) is the trivial character of 
$\bar A_{\la}^F$.
Let us write $\ds_x = \ds z_x$ with $z_x \in Z_{\wt G^*}$.
Since $\ds_x$ is $xF'$-stable, we have
\begin{equation*}
\tag{2.7.2}
\ds\iv \dx\ds\dx\iv = z_xF(z_x)\iv, 
\end{equation*} 
where $\dx \in N_{G^*}(T^*)$ is a representative of 
$x \in (\Om_s)_{F'}$.  We may assume that 
$z_x \in Z_{\wt G^*}^{F^m}$ for a large $m$. 
Let $\wh\p'_x$ be the linear character of $\wt G^{*F^m}$
corresponding to $z_x$.  Since $\ds_x, \ds \in \wt T^*$, and 
$\ds_x$ is $xF'$-stable, $\ds$ is $F'$-stable, we see that
$\wh\p'_x$ is also $F'$-stable. 
As explained in [S2, 2.13], there exists an $F$-stable Levi 
subgroup $\wt M$ of $\wt G$ containing $\wt T$ such that 
$\wt M$ contains $Z_{\wt L}(\la)$ and that $\ds$ is contained
in the center of the dual group of $\wt M$.  This implies that 
the restriction $\wh\p_x$ of $\wh\p'_x$ on $Z_{\wt L}(\la)^{F^m}$
is $F$-stable (cf. [S2, Prop. 2.14]).  
We define a linear character $\p_x$ of 
$Z_{\wt L}(\la)^F$ by $\p_x = Sh_{F^m/F}(\wh\p_x)$.
Since $\xi_1 = 1$, we see that $\xi_x$ is obtained from the 
restriction of $\p_x$ to $Z_L(\la)^F$.
\par
Next, we shall describe the bijection between $(\bar A_{\la})_F$ and 
$(\Om_s^{F'})\wg$.  There exists a surjective
homomorphism  $f_1: \wt G^F/G^F \to (\bar A_{\la})_F$ defined as follows
(cf. [S2, 2.19]).  For $g \in \wt G^F$, we can write $g = g_1z$, with 
$g_1 \in G, z \in Z_{\wt G}$.  Then $g_1\iv F(g_1) \in Z_G$, and it
determines an element in $A_{\la} = Z_L(\la)/Z_L^0(\la)$, and so 
an element in $\bar A_{\la}$, which is unique up to
$F$-conjugacy.  
On the other hand, we construct $f_2: \wt G^F/G^F \to (\Om_s^{F'})\wg$
as follows. 
From (2.7.2), we have $\ds\iv\dx\ds\dx\iv \in Z_{\wt G^*}^F$ (we may 
choose $\dx \in N_{G^*}(T^*)^F$), and this defines a well-defined
injective homomorphism 
$f_2^*: \Om_s^{F'} \to Z_{\wt G^*}^F, x \mapsto \ds\iv\dx\ds\dx\iv$. 
Since $Z_{\wt G^*}^F \simeq (\wt G^F/G^F)\wg$, we have a surjective 
map $f_2$ as the transpose of $f_2^*$.
Then $\Ker f_1 = \Ker f_2$, and these maps induce the bijection 
$f: (\Om_s^{F'})\wg \to (\bar A_{\la})_{F}$. 
\par
Now the parametrization is given as follows.  There exists a unique
$\r_0 \in \Irr G^F$ such that $\r_0$ occurs in $\wt\r_{\ds_x, E}|_{G^F}$ 
and in $\vG_N$.  In our parametrization, then 
$\r_0 = \r_{1, \xi_x} = \r_{1, x}$ ($(1, \xi_x) \in \ol\CM_{s, N}$, 
$(1, x) \in \ol\CM_{s,E}$).  Then any $\r$ contained in 
$\wt\r_{\ds_x, E}|_{G^F}$ is obtained as ${}^g\r_0$ with 
$g \in \wt G^F/G^F$.  We then have $\r = \r_{c,\xi_x} = \r_{\e, x}$ with 
$c = f_1(g)$ and $\e = f_2(g)$.  
\par
By summing up the above argument, we obtain a bijection 
\begin{equation*}
(\Om_s^{F'})\wg \times (\Om_s)_{F'} \to 
     (\bar A_{\la})_F \times (\bar A_{\la}^F)\wg 
        \quad (\e, x) \mapsto (f(\e), \xi_x).
\end{equation*}
This gives the required bijection 
$\ol\CM_{s, E} \simeq \ol\CM_{s, N}$.
\para{2.8.} 
Slightly modifying the arguments in [S2, 4.5], (see the 
remark below),  
we establish a parametrization of 
$\CE(G^F, \{ s\})$ as in (1.6.4).  
We give a bijection 
$\CT_{s, E} \lra \ol\CM_{s, E}$ for each pair
$(s, E)$ such that $E \in (W_{\ds}\wg/\Om_s)^{F'}$.  
\par\medskip
(a) \ First we consider the case where the pair $(s, E)$
satisfies the property (2.5.1).  
If we put $\th = \vD(\wt\r_{\ds,E})|_{Z_L(\la)^F}$, then 
$\th$ satisfies the property (2.1.3).
Hence we have a natural bijection 
$\CT_{s, E} \lra \ol\CM_{s, N} \lra \ol\CM_{s, E}$
by Theorem 2.6 (ii), (b) together with the argument in 2.7. 
\par\medskip
(b) \ Next we consider the case where 
$W_s$ satisfies the same assumption as in (2.5.1), but 
$E$ is not of the form there.  So we assume that 
$\Om_s(E) \ne \Om_s$, and put $t' = |\Om_s(E)|$.  
Replacing $(s,E)$ by a certain $N_W(W_{\ds})$-conjugate, we may 
assume that $E$ can be written as 
$E \simeq E_1\boxtimes\cdots\boxtimes E_r$, ($E_i \in W_{i, \ds}\wg$)
with   
\begin{equation*}
E_i = (E_{i1}\boxtimes\cdots\boxtimes E_{i1})\boxtimes\cdots\boxtimes 
      (E_{ik}\boxtimes\cdots\boxtimes E_{ik}), 
\end{equation*}
where $E_{i1}, \dots, E_{ik}$ are distinct irreducible characters
of $\FS_{b_i}$ with $k = t/t'$, and $E_{ij}$ appears $t'$ times in 
the components of $E_i$.  Moreover, $\Om_s(E)$ acts transitively
on the factors $E_{ij}$.  
Since $E \in (W_{\ds}\wg/\Om_s)^{F'}$, there exists $a_E \in \Om_s$ such 
that $E \in (W_{\ds}\wg)^{F''}$ with $F'' = a_EF'$.
Let $\wt L = \wt L_1 \times\cdots\times \wt L_r$ 
be an $F$-stable Levi subgroup of $\wt G$ according to the 
decomposition of $E$, where 
$\wt L_i = \wt L_{i1} \times \cdots\times \wt L_{ik}$ with 
$\wt L_{ij} \simeq GL_{b_it'}$.  
Then $W_{\ds}$ coincides with $W_{\wt L^*, \ds}$, the stabilizer of $s$ in 
$W_{\wt L^*}$, and $F''$ can be written as $F'' = Fw_2$ with 
$w_2 \in W_{\wt L^*}$.  
Moreover, we have $\Om_s(E) = \Om_{s, L}$, a similar group as 
$\Om_s$ for $L = \wt L \cap G$, and the pair $(s,E)$ satisfies 
the condition in (2.5.1) with respect to $L$.  Hence by (a), 
the set $\CT_{s, E}^L$ is parametrized by 
$\ol\CM_{s, E}^L$ (the super script $L$ denotes the 
corresponding object in $L$).
Let $\wt P$ be the standard parabolic subgroup of $\wt G$ containing
$\wt L$ and put $P = \wt P \cap G$.  
Then by Lemma 4.2 in [S2], the map 
$\r_0 \mapsto \Ind_{P^F}^{G^F}\r_0$ gives a bijection 
between $\CT_{s, E}^L$ and $\CT_{s, E}$.
Since
\begin{equation*} 
\ol\CM_{s, E}^L = (\Om_{s,L}^{F''})\wg \times (\Om_{s,L})_{F''}
  = (\Om_{s,E}^{F'})\wg \times (\Om_{s,E})_{F'} = \ol\CM_{s, E},
\end{equation*} 
this gives a bijection $\CT_{s, E} \lra \ol\CM_{s,E}$.
\par\medskip
(c)\ We consider the general case.  Let 
$W_{\ds} = W_{1,\ds_1}\times \cdots\times W_{r,\ds_r}$, and 
$W_{s} = W_{\ds}\Om_s$.  Here 
we assume that there exists $i$ such that $\Om_s$ acts 
non-transitively on $W_{i,\ds_i}$.  
In this case, there exists a proper Levi subgroup $L^*$ of $G^*$
such that $W_{s}$ is contained in $W_{L^*}$ and that
$L^*$ is both $F$-stable and $F'$-stable.  
Then $Z_{G^*}(s)$ is contained in $L^*$.   
Under this condition, it is known that the twisted induction 
$R^G_L(\dw_1)$ (see, e.g., [S2, 3.1]) induces a bijection between 
$\CE(L^{F'}, \{s\})$ and $\CE(G^F,\{ s\})$.
By induction hypothesis, we may assume that there exists a 
bijection $\CT^L_{s, E} \lra \ol\CM_{s, E}^L$.
Since $\Om_s(E) = \Om_{s,L}$, $\ol\CM_{s, E}^L$ is 
identified with $\ol\CM_{s, E}$.
Hence we have a bijection $\CT_{s, E} \lra \ol\CM_{s,E}$
as asserted.
\par\medskip\noindent
\remark{2.9.}\ 
In [S2], 4.5, the parametrization is done through three steps as above.
However, in the step (a), only the pair $(s, E)$ such that 
$\vD(\wt\r_{\ds,E}) = 1$ is treated, and it is stated that 
other cases are
reduced to this the case by considering the linear character 
$\th$ of $G^F$ corresponding to the central element 
$\dz \in Z_{G^*}^F$.  But this is not true in general.  
In fact, if the $F$-stable class $\{ \ds'\}$ in $G^*$ satisfies 
the property in (2.5.1), then $\ds'$ can be written as 
$\ds' = \dz\ds$ for an $F$-stable class $\{\ds\}$ such that 
$\vD(\wt\r_{\ds,E}) = 1$ with $\dz \in Z_{G^*}$. 
However, it occurs that $\dz \notin Z_{G^*}^F$ even if 
the classes $\{\ds\}$ and $\{ \dz\ds\}$ are $F$-stable.
In that case one cannot find a linear character $\th$ of $G^F$
corresponding to $\dz$.  Hence the step (a) in [S2] 
does not cover all the cases, and one needs to consider the cases
where $\vD(\wt\r_{\ds,E}) \ne 1$ discussed as in 2.8. 
\par\medskip
\section{Character formula for generalized Gelfand-Graev characters}
\para{3.1.}
For later use, we shall prove a character formula for $\vG_{c,\xi}$,
which is a variant of the formula given in [K3].
Note that $\vG_{c,\xi}$ is constructed by using a specific Lagrangian 
subspace $\Fs$ of $\Fu_1$. Following [S2, 2.3], we recall the
construction of $\Fs$.  Assume, for simplicity, that $G = SL_n$.
The weighted Dynkin diagram of $N$ is given as follows.  Let 
$\Pi \subset \vS^+$ be the set of simple roots and 
the set of positive roots of $G$, which is written in the form 
$\vS^+ = \{ \ve_i - \ve_j \mid 1 \le i < j \le n\}$ for certain 
basis vectors $\ve_1, \dots, \ve_n$ of $\BR^n$, 
and $\Pi = \{ \a_1, \dots, \a_{n-1}\}$ with 
$\a_i = \ve_i - \ve_{i+1}$.  
Assume that $N$ corresponds to a partition 
$\mu = (\mu_1 \ge \mu_2 \ge\cdots\ge \mu_r>0)$ of $n$ via Jordan's
normal form.  
For each $\mu_i$, put 
\begin{equation*}
Y_i = \{ \mu_i-1, \mu_i-3, \dots, -\mu_i+1 \}
\end{equation*}
consisting of $\mu_i$ integers.  Then $Y = \coprod_i Y_i$ is 
a set of $n$ integers (with multiplicities), and we arrange 
its elements in a decreasing order, 
\begin{equation*}
\tag{3.1.1}
Y = \{ \n_1 \ge \n_2 \ge \cdots \ge \n_n\}.
\end{equation*} 
The weighted Dynkin diagram $h : \Pi \to \BZ$ is given by 
$h(\a_i) = \n_i - \n_{i+1}$ for $1 \le i \le n-1$.
Let $\Pi_1$ (resp. $\vS_1$) be the set of $\a \in \Pi$ 
(resp. $\a \in \vS^+$) such that $h(\a) = 1$.  Clearly we have
$\Fg_1 = \bigoplus_{\a \in \vS_1}\Fg_{\a}$. 
The set $\vS_1$ is described as follows. 
For a given $\a_i \in \Pi_1$, let $j$ be the smallest integer such
that $j > i$ and that $h(\a_j) >0$, and let $k$ be the largest 
integer such that $k < i$ and that $h(\a_k) > 0$.
We define a subset $\Psi_i$ of $\vS^+$ by 
\begin{equation*}
\Psi_i = \{ \ve_p - \ve_q \mid k+1 \le p \le i, i+1 \le q \le j\}.
\end{equation*}
(If $j$ or $k$ does not exist, we put $k = 0$ or $j = n$.)
Then it is easy to see that $\Psi_i$ are mutually disjoint and that
\begin{equation*}
\vS_1 = \coprod_{\a_i \in \Pi_1}\Psi_i.
\end{equation*}
For $\a_i, \a_j \in \Pi_1$ such that $i < j$, we say that 
$\Psi_i$ and $\Psi_j$ are adjacent if $\a_k \notin \Pi_1$ 
for $i < k < j$. 
There exists a subset $\Psi$ of $\vS_1$ satisfying the following
properties;
$\Psi$ is a union of the $\Psi_i$ which are not adjacent 
each other, and $\vS_1 = \Psi \coprod \s(\Psi)$, where 
$\s$ is the permutation of $\vS^+$ induced from the graph 
automorphism of $\Pi$. Note that $\Psi$ is uniquely determined up
to the action of $\s$. 
Put 
$\Fs = \bigoplus_{\a \in \Psi}\Fg_{\a}$.
Then it was shown in [S2, 2.3] that $\Fs$ is a Lagrangian subspace 
in $\Fg_1$, stable by the action of $L$.
\par
In the discussion below, we follow the notation in 2.1.
Let $V$ be an $n$-dimensional vector space over $k$ on which $G$ acts 
naturally.  We can find a basis 
$\{ N^jv_i \mid 1 \le i \le r,  0 \le j < \mu_i\}$ of $V$ such that
$N^{\mu_i}v_i = 0$ and that 
$H$ acts on $N^jv_i$ by a scalar multiplication 
$-\mu_i +1 + 2j$.  
\par
Put $M = Z_L(\la)$. 
Take $t \in M$ such that $t$ stabilizes each basis vector
$N^jv_i$ up to scalar.  It follows that $t \in T_1$, where $T_1$
is a maximal torus in $L$ related to the weighted Dynkin 
diagram of $N$. 
Since $G$ is simply connected, $Z_G(t)$ is connected.
Put $\Fz_{t} = \Lie Z_G(t)$.
Since $M = Z_G(N) \cap Z_G(N^*)$, 
$\Fz_t$ contains $N, N^*$, and so it contains $H$.
If we put $(\Fz_t)_j = \Fg_j \cap \Fz_t$,
$\Fz_t = \bigoplus_j (\Fz_t)_j$ gives the grading of $\Fz_t$
associated to $N \in \Fz_t$.  Put 
$(\Fu_t)_j = \bigoplus_{j' \ge j}(\Fz_t)_{j'}$.  Then we have 
$(\Fu_t)_j = \Fu_j \cap \Fz_t$.  In particular, 
$P_t = P \cap Z_G(t)$ (resp. $L_t = L \cap Z_G(t)$)
is the parabolic subgroup of $Z_G(t)$ (resp. its Levi subgroup)
 associated to $N$.
Moreover, the restriction of $\la : \Fu_1 \to k$ to $(\Fu_t)_1$
coincides with the corresponding linear map $\la_t$ with respect
to $N$. 
We have the following lemma. 
\begin{lem}  
  $\Fs_t = \Fs \cap \Fz_t$ is 
  a Lagrangian subspace of $(\Fz_t)_1$, which 
  is stable by $L_t$.
\end{lem}
\begin{proof}
  From the above discussion, the symplectic form on $(\Fz_t)_1$ is 
  obtained as the restriction to $(\Fz_t)_1$ of the symplectic form 
  on $\Fg_1$.  Hence 
  $\Fs \cap \Fz_t$ is an anisotropic subspace of $(\Fz_t)_1$.
  Also it is clear that $\Fs \cap \Fz_t$ is stable by $L_t$.
  Note that $T_1$ is a maximal torus on which the root system $\vS$ 
  is defined.  Then the 
  subroot system $\vS_t$ of $\vS$ associated to the group 
  $Z_G(t)$ consists of roots $\ve_p - \ve_q \in \vS $ such that the
  corresponding basis vectors $N^kv_l$ and $N^{k'}v_{l'}$ in $V$ have
  the eigenvalue 0 for $t$.  (We identify the basis $\{ \ve_i\}$ and 
  $\{ N^jv_i\}$ via the total order $\nu_1 , \dots, \nu_n$ in
  (3.1.1)).
  Since $t \in M$, $N^jv_i$ 
  have the same eigenvalue for all $1 \le j \le \mu_i-1$ ($i$ is fixed).
  It follows that the subroot system $\vS_t$ is invariant under 
  the action of $\s$.  Hence $\vS_1 \cap \vS_t$ is also 
  $\s$-invariant.  If we put $\Psi_t = \Psi \cap \vS_t$, then 
  we have $\vS_1 \cap \vS_t = \Psi_t \coprod \s(\Psi_t)$. 
  In particular, $|\Psi_t| = |\vS_1 \cap \vS_t|/2$.
  Since $\Fs \cap \Fz_t = \bigoplus_{\a \in \Psi_t}\Fg_{\a}$, 
  we see that $\Fs \cap \Fz_t$ is a Lagrangian subspace of 
  $(\Fz_t)_1$.  The lemma is proved.
\end{proof}  
\para{3.3.}
For $c \in A_{\la}$, we choose $\dc \in Z_L \subset T_1$, and take 
$\a_c \in T_1$ such that $\a_c\iv F(\a_c) = \dc$.  We consider 
the $c$-twisted version of the previous results.  In the following, 
we denote by $X_c$ the object obtained from $X$ related to $G$ or 
$\Fg$ by the conjugation or adjoint action by $\a_c$. 
Then $M_c = \a_c M \a_c\iv$ coincides with  $Z_L(\la_c)$, and we 
have $M_c = Z_G(N_c) \cap Z_G(N_c^*)$, where $N_c, N_c^*, H_c$ are
$F$-stable TDS-triple.  
Since $\a_c \in L$, we have $\Fs_c = \Fs$ and $\Fu_{2,c} = \Fu_2$.
It follows that $\Fu_c = \Fu$ for $\Fu = \Fs + \Fu_2$,
and so $U_c = U$.  
Let $t$ be a semisimple element in $M_c^F$. 
Put $\Fz_t = \Lie Z_G(t)$, etc., as before.
If we put $t' = \a_c\iv t\a_c \in M$, $t'$ is conjugate to 
an element in $T_1$ under $M^0 \subset L$ (since $M_c = M_c^0Z_G$). 
It follows that 
$\Fs \cap \Fz_{t'}$ is a Lagrangian subspace of $(\Fz_{t'})_1$, and so
$\Fs \cap \Fz_t = {}^{\a_c}(\Fs \cap (\Fz_{t'})_1)$ is 
a Lagrangian subspace of $(\Fz_t)_1$, which is stable by $F$.  
We put $\Fs_t = \Fs \cap \Fz_t$.
Let $\Fu_t = \Fs_t + (\Fu_t)_2$.  Then 
we have an $F$-stable subgroup $U_t$ of $Z_G(t)$ such that 
$\Lie U_t = \Fu_t$, which is stable by 
$L_t$.  Moreover, $\Fu_t = \Fu \cap \Fz_t$.  
It follows from this that 
\begin{equation*}
\tag{3.3.1}
U \cap Z_G(t) = U_t.
\end{equation*}
\par
Now $\la_c: \Fu_1 \to k$ is the linear map defined as $\la$, 
by using $N_c$ instead of $N$.  Then $\la_{c,t} = \la_c|_{(\Fu_t)_1}$ is 
the linear map on $(\Fu_t)_1$ defined by $N_c \in \Fz_t$.
It follows that the restriction of $\vL_c : U^F \to \Ql^* $ 
on $U_t^F$ coincides with the linear character of $U_t^F$ defined in terms
of $N_c$, which we denote by $\vL_{c,t}$.  
\para{3.4.}
Take a semisimple element $s \in G^F$, and assume that there exists
$g \in G^F$ such that $g\iv sg \in M_c^F$.  
(Do not confuse $s$ with an element in the dual group $G^*$).
 By fixing $s$, we put 
$P_g = gPg\iv \cap Z_G(s)$ and $L_g = gLg\iv \cap Z_G(s)$.
We apply the previous argument for $t = g\iv sg$.
Then $N_g = {}^gN_c$ is an $F$-stable nilpotent element in $\Lie Z_G(s)$, and 
$U_g = gU_tg\iv$ is the unipotent subgroup of $Z_G(s)$ associated to
$N_g$.  
We have $P_g = gP_tg\iv$ and $L_g = gL_tg\iv$.
Moreover $\la_g = \Ad g\circ \la_{c,t}$ is the linear map 
of ${}^g(\Fu_t)_1 = (\Fu_s)_1$, and 
$\vL_g = \ad g\circ \vL_{c,t}$ coincides with the linear 
character of $U_g^F$ associated to $N_g$.
We define a modified generalized 
Gelfand-Grave character  $\vG^{Z_G(s)}_{N_g,1}$ of $Z_G(s)$ by
\begin{equation*}
\tag{3.4.1}
\vG^{Z_G(s)}_{N_g,1} = 
\Ind_{Z_{L_g}(\la_g)^FU_g^F}^{Z_G(s)^F}(1\otimes\vL_g)
\end{equation*}
associated to $N_g \in \Lie Z_G(s)$.
\par
The following result gives a description of modified generalized
Gelfand-Graev characters in terms of various modified 
generalized Gelfand-Graev characters of smaller groups, 
which is an extended version of the formula 
stated in [K3, Lemma 2.3.5]. 
\begin{prop} 
 Assume that $s, v \in G^F$ such that $sv = vs$, where $s$ is
 semisimple and $v$ is unipotent.
 Then we have
  \begin{equation*}
    \vG_{c,\xi,\th}(sv) = \frac{1}{|Z_G(s)^F|}
               \sum_{\substack{g \in G^F \\ g\iv sg \in Z_L(\la_c)^F}}
   \frac{|Z_{L_g}(\la_g)^F|}{|Z_L(\la_c)^F|}
    \th\xi\nat(g\iv sg)\vG^{Z_G(s)}_{N_g, 1}(v).
  \end{equation*}
\end{prop}
\begin{proof}  
By definition,  
  \begin{align*}
    \tag{3.5.1}
    \vG_{c,\xi,\th}(sv) &= 
    \bigl(\Ind_{M_c^FU^F}^{G^F}\th\xi\nat\otimes\vL_c\bigr)(sv) \\
    &= |M_c^FU^F|\iv
    \sum_{\substack{ g \in G^F \\ g\iv svg \in M_c^FU^F}}
    (\th\xi\nat\otimes\vL_c)(g\iv svg).
  \end{align*}
  Here the condition $g\iv svg \in M_c^FU^F$ in the sum 
  is equivalent to the condition that $g\iv vg \in M_c^FU^F$ and 
  $g\iv sg \in M_c^FU^F$.    
We note that 
\par\medskip\noindent
(3.5.2) \ Any semisimple element in $M_c^FU^F$ is contained in 
$\bigcup_{x \in U^F}xM_c^Fx\iv$.
\par\medskip
In fact, let $T_1$ be an $F$-stable maximal torus in $L$ as in 3.1.
Then $T_2 = T_1 \cap M$ is a maximal torus in $M$.  
Since we can choose $\dc \in T_1$, $T_2$ is also 
contained in $M_c$. Thus $T_2$ is a maximal torus in $M_cU$.  We have
$N_{M_cU}(T_2) \simeq N_{M_c}(T_2)Z_U(T_2)$.   
Since $U$ is a product of one parameter subgroups $U_{\a}$ 
associated to roots $\a$ with respect to $T_1$, $Z_U(T_2)$ is a
product of $U_{\a}$ such that $\a|_{T_2} = \Id$.  It follows that 
$Z_U(T_2)$ is connected, and $N_{M_cU}(T_2)^0 = T_2Z_U(T_2)$.
We see that $N_{M_cU}(T_2)/N_{M_cU}(T_2)^0 \simeq N_{M_c}(T_2)/T_2$.
This implies that any $F$-stable maximal torus in $M_cU$ is taken from
$M_c$ up to $U^F$-conjugate. 
Since any semisimple element in $M_c^FU^F$ is contained in an $F$-stable
maximal torus, we obtain (3.5.2).
\par\medskip
It follows from (3.5.2) that 
$g\iv sg \in x M_c^Fx\iv$ for some $x \in U^F$, i.e., 
  $(gx)\iv s(gx) \in M_c^F$.
It is easy to see that the set 
  $\{ x_1 \in U^F \mid (gx_1)\iv s(gx_1) \in M_c^F\}$ is given 
  by $xZ_U((gx)\iv s(gx))^F$ for some $x \in U^F$ such that 
  $(gx)\iv s(gx) \in M_c^F$.  Hence the last formula in (3.5.1) 
implies that 
  \begin{equation*}
\vG_{c,\xi,\th}(sv) = |M_c^FU^F|\iv\sum_{\substack{ g \in G^F, x \in U^F \\ 
        g\iv vg \in M_c^FU^F \\
        g\iv sg \in x M_c^Fx\iv}}|Z_U((gx)\iv s(gx))^F|\iv
    (\th\xi\nat\otimes\vL_c)(g\iv svg).
  \end{equation*} 
By replacing $gx$ by $g$, we have 
  \begin{align*}
\vG_{c,\xi,\th}(sv) &= |M_c^F|\iv\sum_{\substack{ g \in G^F \\
        g\iv vg \in M_c^FU^F \\
        g\iv sg \in M_c^F}}
    |Z^F_U(g\iv sg)|\iv\th\xi\nat(g\iv sg)      
    (\th\xi\nat\otimes\vL_c)(g\iv vg)  \\
    &= |M_c^F|\iv\sum_{t }
          \sum_{\substack{ y\in Z_G(t)^F \\
                    y\iv v_1y \in M_c^FU^F \\
             }}
     |Z_U^F(t)|\iv\th\xi\nat(t)(\th\xi\nat\otimes\vL_c)(y\iv v_1y), 
  \end{align*}
where in the first sum in the last formula, $t$ runs over all the
semisimple element in $M_c^F$ such that $t = g\iv sg$ for some $g \in
G^F$.  We fix such $g$ for each $t$, and put $v_1 = g\iv vg$.
\par
Hence we have $v_1 \in Z_G(t)^F$.  Since $t$ normalizes $M_c^F$ and
$U^F$, we have 
\begin{align*}
Z_G(t)^F \cap M_c^FU^F &= 
   (Z_G(t)^F \cap Z_L(\la_c)^F)(Z_G(t)^F \cap U^F) \\
    &= Z_{M_c}(t)^FU_t^F 
\end{align*}
by (3.3.1).  Also we have $Z_U(t) = U_t$ by (3.3.1).  
It follows that we have 
\begin{equation*}
\tag{3.5.3}
\begin{split}
\vG_{c,\xi,\th}(sv) = 
    &\sum_{\substack{t \in M_c^F \\ g\iv sg = t}}
          \frac{|Z_{M_c}(t)^F|}{|M_c^F|}\th\xi\nat(t)\times \\ 
      &\times \biggl\{|Z_{M_c}(t)^F|\iv|U_t^F|\iv
       \sum_{\substack{y \in Z_G(t)^F \\ y\iv v_1y \in Z_{M_c}(t)^FU_t^F}}
              (\th\xi\nat\otimes\vL_c)(y\iv v_1y)\biggr\}.
\end{split}
\end{equation*} 
Here we note that $y\iv v_1y$ is unipotent.  Hence the component
of $y\iv v_1y$ in $Z_{M_c}(t)^F$ is unipotent.  Since $\xi\nat$ is 
a character of $M_c^F$ which is trivial on $M_c^{0F}$, it is trivial 
on the set of unipotent elements in $M_c^F$.  Also by (2.1.3) 
$\th$ is trivial on the set of unipotent elements in $M_c^F$.  
It follows that 
\begin{equation*}
(\th\xi\nat\otimes\vL_c)(y\iv v_1y) = (1\otimes\vL_c)(y\iv v_1y).
\end{equation*}
Then the expression in the parenthesis in (3.5.3) coincides with
\begin{equation*}
\Ind_{Z_{M_c}(t)^FU_t^F}^{Z_G(t)^F}(1\otimes\vL_c)(v_1)
  = \vG^{Z_G(s)}_{N_g, 1}(v)
\end{equation*}
under the conjugation by $g \in G^F$.  Substituting this 
into (3.5.3) we obtain the proposition.
\end{proof}
\section{Shintani descent and almost characters}
\para{4.1.}
We consider the group $G^{F^m}$ for a positive integer $m$.  
We denote 
by $G^{F^m}\ssim_F$ the set of $F$-twisted conjugacy classes in 
$G^{F^m}$. (In the case where $m = 1$, the set of $F$-twisted classes 
coincides with the set of conjugacy classes, which we denote simply by
$G^F\ssim$.)  A norm map 
\begin{equation*}
N_{F^m/F} : G^{F^m}\ssim_F \ \to G^F\ssim
\end{equation*}
is defined by attaching $x = F^m(\a)\a\iv$ to $\hat x = \a\iv F(\a)$
where $x \in G^F, \hat x \in G^{F^m}$ and $\a \in G$. 
Let $C(G^{F^m}\ssim_F)$ (resp. $C(G^F\ssim)$) be the space of
$F$-twisted class functions on $G^{F^m}$ (resp. class functions on
$G^F$).
A Shintani descent map 
\begin{equation*}
Sh_{F^m/F} : C(G^{F^m}\!\ssim_F) \to C(G^F\ssim)
\end{equation*}
is given by $Sh_{F^m/F} = (N_{F^m/F}^*)\iv$, 
which is a linear isomorphism of vector spaces.
\par
Let $\s = F|_{G^{F^m}}$.  We consider the semidirect product 
$G^{F^m}\lp\s\rp$ of $G^{F^m}$ with the cyclic group $\lp\s\rp$
of order $m$ generated by $\s$.  Then the coset $G^{F^m}\s$ is 
invariant under the conjugation action of $G^{F^m}$, and the set 
$G^{F^m}\s\ssim$ is identified with the set $G^{F^m}\ssim_F$
via the map $x\s \lra x$.
Now each $F$-stable irreducible character $\r$ of $G^{F^m}$ can be
extended to an irreducible character $\wt\r$ of $G^{F^m}\lp\s\rp$
(in $m$-distinct way), and the restriction $\wt\r|_{G^{F^m}\s}$ to
the coset $G^{F^m}\s$ determines an element in $C(G^{F^m}\ssim_F)$
under the above bijection.  The function $\wt\r|_{G^{F^m}\s}$ does 
not depend on the choice of the extension up to a scalar multiple, 
and the collection of those $\wt\r|_{G^{F^m}\s}$ for 
$\r \in (\Irr G^{F^m})^F$ gives a basis of $C(G^{F^m}\ssim_F)$. 
In what follows, we often regard a character $f$ of $G^{F^m}\lp\s\rp$ 
as an element in $C(G^{F^m}\ssim_F)$ by considering its restriction 
to $G^{F^m}\s$, if there is no fear of confusion.
\para{4.2.}
We shall describe the Shintani descent of the modified generalized 
Gelfand-Graev characters.  We follow the setting in 2.2.  Recall the 
set $\ol\CM$ in (2.1.2) and $\CM$ in 2.2.  Let 
$\th$ be as in (2.2.1).  Hence it is the restriction to 
$Z_L(\la)^{F^m}$ of an $F$-stable 
linear character $\th'$ of $Z_{\wt L}(\la)^{F^m}$.  We denote by 
$\th_0$ the linear character
of $Z_L(\la)^F$ obtained by restricting the linear character 
$Sh_{F^m/F}(\th')$ of $Z_{\wt L}(\la)^F$.
Hence $\th_0$ satisfies the condition in (2.1.3).  
We consider the modified generalized Gelfand-Graev characters 
$\vG^{(m)}_{c,\xi, \th}$ and $\vG_{c_1,\xi_1,\th_0}$ for 
$(c,\xi) \in \CM$ and $(c_1,\xi_1) \in \ol\CM$.
\par
Let us consider an extension $\wt\vG^{(m)}_{c,\xi,\th}$ as 
in 2.2, which is determined by the choice of 
an extension $\wt{\th\xi}\nat$ of $\th\xi\nat$ to 
$M_c\lp\hat c\s\rp$.  Since $c \in A_{\la}^F$, we may choose 
$\dc \in Z_L(\la)^F$.
Note that, under the isomorphism 
$\ad\b_c\iv: Z_L(\la_c)^{F^m} \simeq Z_L(\la)^{\dc F^m}$, 
the linear character $\th\xi\nat$ corresponds 
to a linear character $\th\xi$ of $Z_L(\la)^{\dc F^m}$, 
and $\wt{\th\xi}\nat$ corresponds to its extension 
$\wt{\th\xi}$ to $Z_L(\la)^{\dc F^m}\lp\s\rp$.  
Take $c_1 \in (A_{\la})_F$. As $A_{\la} = A_{\la}^{F^m}$, 
we may choose an element $\dc_1 \in Z_L(\la)^{F^m}$ whose image 
on $A_{\la}$ gives a representative of $c_1 \in (A_{\la})_F$.
Now the following proposition 
describes the Shintani descent of $\wt\vG^{(m)}_{c,\xi,\th}$ 
in terms of $\vG_{c_1,\xi_1, \th_0}$.  The proof is done in a similar
way as in [S2].  In fact, Theorem 1.10 in [S2] can be extended 
to our setting, and the proposition is the direct consequence of the 
theorem (cf. [S2, 4.11]).
\begin{prop}
Let the notations be as above.  Assume that $m$ is sufficiently 
divisible.  Then we have 
\begin{equation*}
Sh_{F^m/F}(\mu_{c,\th\xi}\iv \wt\vG^{(m)}_{c,\xi,\th})
  = |A_{\la}^F|\iv\wt\xi(\s)\sum_{(c_1,\xi_1) \in \ol\CM}
               \xi(c_1)\xi_1(c)\vG_{c_1,\xi_1, \th_0}.
\end{equation*}
\end{prop}
\para{4.4.}
We shall describe the set of $F$-stable irreducible characters
of $G^{F^m}$ in the case where $m$ is sufficiently divisible.
Let $\{s\}$ be an $F$-stable class in $G^*$, and we assume that
$s \in T^*$.  As in the case of $G^F$, one can find 
$\ds \in \wt T^*$ such that $\pi(\ds) = s$ and that the class 
$\{ s\}$ is $F$-stable.  Hence $F'(\ds) = \ds$ for $F' = Fw_1$.    
We choose $m$ large enough so that $\ds \in \wt T^{*F^m}$ and that
$F^m$ acts trivially on $\Om_s$.  For each $E \in W_{\ds}\wg$, we denote 
by $\ol\CM_{s, E}^{(m)}$ and $\CT^{(m)}_{s, E}$ the set
$\ol\CM_{s, E}$ and $\CT_{s, E}$ as given in 1.4, but 
replacing $F'$ by $F^m$.  Then  $\CE(G^{F^m}, \{ s\})$ is a 
disjoint union of various $\CT^{(m)}_{s, E}$, and 
the latter set is in bijection with $\ol\CM_{s, E}^{(m)}$.
By our assumption on $m$, we have 
$\ol\CM_{s, E}^{(m)} = \Om_{s,E}\wg \times \Om_{s,E}$.
Let us define a subset $\CM_{s, E}$ of 
$\ol\CM^{(m)}_{s, E}$ by 
$\CM_{s, E} = (\Om_{s,E}\wg)^{F'}\times \Om_{s,E}^{F'}$, where
$(\Om_{s,E}\wg)^{F'}$ means the set of $F'$-stable irreducible characters
of $\Om_{s,E}$.  
Then by [S2, (4.6.1)], the set $(\CT_{s,E}^{(m)})^F$ of $F$-stable
irreducible characters in $\CT_{s, E}^{(m)}$ is parametrized by 
$\CM_{s, E}$, and so $\CE(G^{F^m}, \{ s\})^F$ 
can be described as 
\begin{equation*}
\tag{4.4.1}
\CE(G^{F^m}, \{s\})^F = 
    \coprod_{E \in (W_{\ds}\wg/\Om_s)^{F'}}\CM_{s, E}.
\end{equation*}
\par
In the case where $(s, E)$ is of the form 2.8 (a), the set
$\CT_{s, E}$ is also parametrized in terms of $\ol\CM_{s, N}$.
Since $m$ is large enough, $\CT_{s,E}^{(m)}$ is parametrized 
by $\ol\CM_{s, N}^{(m)} = \bar A_{\la} \times \bar A_{\la}\wg$.
Then under this parametrization, $(\CT_{s,E}^{(m)})^F$ is parametrized
by $\bar A_{\la}^F \times (\bar A_{\la}\wg)^F$.
\para{4.5.}
We define a pairing 
$\{ \ , \ \} : \CM_{s, E} \times \ol\CM_{s, E} \to \Ql^*$ 
as follows.  For $x = (\e, z) \in \CM_{s, E}$ and 
$y = (\e',z') \in \ol\CM_{s, E}$, 
\begin{equation*}
\tag{4.5.1}
\{ x, y\} = |\Om_{s,E}^{F'}|\iv\e(z')\e'(z).
\end{equation*}    
(Note that $\e \in (\Om_{s,E}\wg)^{F'}$ can be viewed as a character of
the group $(\Om_{s,E})_{F'}$.)
\par
We define a function $R_x \in C(G^F\ssim)$ for each 
$x \in \CM_{s, E}$ by 
\begin{equation*}
\tag{4.5.2}
R_x = \sum_{y \in \ol\CM_{s, E}}\{ x, y\}\r_y.
\end{equation*}
\par
In the case where $(s,E)$ satisfies the property in 2.8 (a), 
the set $\CT_{s,E}$ is also parametrized by 
$\ol\CM_{s, N} = (\bar A_{\la})_F\times (\bar A_{\la}^F)\wg$, 
and we have a bijection between $\ol\CM_{s, E}$ and 
$\ol\CM_{s, N}$ by 2.7.  Put 
$\CM_{s,N} = \bar A_{\la}^F\times (\bar A_{\la}\wg)^F$.  
Then the set $(\CT^{(m)}_{s,E})^F$ is parametrized by 
$\CM_{s,N}$.  By modifying the argument in 2.7 appropriately
to the situation in $G^{F^m}$, we have a bijection between 
$\CM_{s,E}$ and $\CM_{s,N}$. 
Let $\th_0$ be the linear character of $Z_L(\la)^F$ obtained 
by restricting $\vD(\wt\r_{\ds,E})$ to $Z_L(\la)^F$. 
The linear character $\th$ of $Z_L(\la)^{F^m}$ is also defined
by using the Shintani descent of $Z_{\wt L}(\la)$ (cf. 4.2).
We say that $\th$ (resp. $\th_0$) is the linear character 
associated to $\CM_{s,N}$ (resp. $\ol\CM_{s,N}$).
\par
We define a pairing 
$\{\ ,\ \}: \CM_{s,N} \times \ol\CM_{s,N} \to \Ql^*$ for
$x = (c,\xi)\in \CM_{s,N}$ and 
$y = (c',\xi') \in \ol\CM_{s,N}$,
\begin{equation*}
\tag{4.5.3}
\{ x,y\} = |\bar A_{\la}^F|\iv\xi(c')\xi'(c).
\end{equation*}   
Then the bijections $\ol\CM_{s,E} \simeq \ol\CM_{s,N}$, etc., 
are compatible with those pairings.
This property was  used in [S2, 4.11] to connect almost 
characters defined in terms of $\ol\CM_{s,N}$ to that of 
$\ol\CM_{s,E}$ (in the case where $\th = 1$, but the proof was omitted
there).  
We give a proof of this property. 
\begin{lem} 
Assume that $(s,E)$ is as in 2.8 (a).  Then
under the bijections $\ol\CM_{s, E} \simeq \ol\CM_{s,N}, 
(\e,z) \lra (c,\xi)$ and 
$\CM_{s,E} \simeq \CM_{s,N}$, $(\e',z') \lra (c',\xi')$, 
we have
\begin{equation*}
 |\Om_{s,E}^{F'}|\iv\e(z')\e'(z) 
     = |\bar A_{\la}^F|\iv\xi(c')\xi'(c).
\end{equation*}
\end{lem}
\begin{proof}
By our assumption, we have $\Om_{s,E}^{F'} = \Om_s^{F'}$.
It follows from the parametrization of $\Irr G^F$ in 1.6 and 2.7,
we see that $|\Om_s^{F'}| = |\bar A_{\la}^F|$, which coincides with
the number of irreducible components in $\wt\r_{\ds,E}|_{G^F}$.   
Thus, in order to prove the lemma, it is enough to show that
\begin{equation*}
\tag{4.6.1}
\e(z') = \xi'(c), \quad \e'(z) = \xi(c').
\end{equation*} 
We recall the bijection $\ol\CM_{s,E} \simeq \ol\CM_{s,N}$ given 
by $h: (\Om_s)_{F'} \to (\bar A_{\la}^F)\wg$ and 
$f: (\Om_s^{F'})\wg \to (\bar A_{\la})_{F}$ in 2.7. 
A similar construction gives  bijections 
\begin{align*}
h'&: \Om_s \to \bar A_{\la}\wg, \qquad 
h'': (\Om_s)^{F'} \to (\bar A_{\la}\wg)^{F}, \\ 
f'&: \Om_s\wg \to \bar A_{\la}, \qquad
f'': (\Om_s\wg)^{F'} \to \bar A_{\la}^F, 
\end{align*}
and $h''\times f'': (\Om_s\wg)^{F'} \times \Om_s^{F'}
       \to \bar A_{\la}^F \times (\bar A_{\la}\wg)^F$
gives the bijection $\CM_{s,E} \to \CM_{s,N}$.
We have inclusions 
$(\Om_s)^{F'} \hra \Om_s, (\Om_s\wg)^{F'} \hra \Om_s\wg$ 
and natural surjections 
$\Om_s \to (\Om_s)_{F'}, \Om_s\wg \to (\Om_s^{F'})\wg$.
Also we have inclusions 
$(\bar A_{\la}\wg)^F \hra \bar A_{\la}\wg, 
 \bar A_{\la}^F \hra \bar A_{\la}$, and natural surjections
$(\bar A_{\la})\wg \to (\bar A_{\la}^F)\wg, 
 \bar A_{\la} \to (\bar A_{\la})_F$.
\par
We want to show that the maps $h, h', h''$ and $f,f',f''$ are
compatible with various inclusions and surjections given above.
First we note that the map $h: (\Om_s)_{F'} \to (\bar A_{\la}^F)\wg$
is compatible with the extension of the filed, i.e., 
the following diagram commutes.
\begin{equation*}
\tag{4.6.2}
\begin{CD}
(\Om_s)_{F'} @>h>> (\bar A_{\la}^F)\wg \\
@AAA                 @AAA             \\ 
(\Om_s)_{F^k} @>h^{0}>> (\bar A_{\la}^{F^k})\wg,
\end{CD}
\end{equation*}
where $F^k$ is the map such that $F^k(\ds) = \ds$ and 
that $F^k$ acts trivially on $\Om_s$, and on $\bar A_{\la}$. 
$h^0$ is a similar map as $h$ defined by replacing $F'$ by $F^k$.
\par
We show (4.6.2).  We choose $m$ large enough so that $m$ is divisible 
by $k$.  Let $\p_x, \wh\p_x, \wh\p_x'$, be the maps given in 2.7 
with respect to
$F'$.   Let $y \in (\Om_s)_{F^k}$ such that its canonical image in 
$(\Om_s)_{F'}$ coincides with $x$.  We may assume that 
$\dx = \dy \in N_{G^*}(T)$.  
We denote by $\p_y, \wh\p_y, \wh\p'_y$ similar maps
constructed by using $F^k$ instead of $F'$.
In particular, $\wh\p'_y$ is the linear 
character of $\wt G^{F^m}$ corresponding to 
$z_y \in Z_{\wt G^*}^{F^m}$ such that 
$\ds\iv\dx\ds\dx\iv = z_yF^k({z_y}\iv)$.
Since $z_yF^k({z_y}\iv) = z_xF(z_x\iv)$, 
one can choose $z_x$ and $z_y$ so that they satisfy the relation
\begin{equation*}
z_x = z_yF(z_y)\cdots F^{k-1}(z_y).
\end{equation*}  
It follows that 
\begin{equation*}
\wh\p_x = \wh\p_yF\iv(\wh\p_y)\cdots F^{-k+1}(\wh\p_y).
\end{equation*} 
Put $\p_x' = Sh_{F^m/F^k}(\wh\p_x)$.  Since 
$\p_y = Sh_{F^m/F^k}(\wh\p_y)$, we have
\begin{equation*}
\p_x' = \p_yF\iv(\p_y)\cdots F^{-k+1}(\p_y), 
\end{equation*}
and $\p_x = Sh_{F^k/F}(\p'_x)$. 
We shall compute the value $\p_x(t)$ for $t \in T^F$.
Take $\a \in T$ such that $t = F^k(\a)\a\iv$, and put
$\hat t = \a\iv F(\a)$.  Then $\hat t \in T^{F^k}$, and 
we have
\begin{equation*}
\p_x(t) = \p'_x(\hat t) = \p_y(\hat tF(\hat t)\cdots F^{k-1}(\hat t)) 
        = \p_y(\a\iv F^k(\a)) = \p_y(t)
\end{equation*}
since $\a\iv F^k(\a) = F^k(\a)\a\iv = t$.  It follows that 
\begin{equation*}
\tag{4.6.3}
\p_x|_{T^F} = \p_y|_{T^F}.
\end{equation*}
Now for $c \in \bar A_{\la}^F$, one can choose a representative
$\dc \in Z_L(\la)^F$ of $c$ so that $\dc \in T^F$.
Then by (4.6.3), we have
\begin{equation*} 
h(x)(c) = \p_x(\dc) = \p_y(\dc) = h^0(y)(c). 
\end{equation*}
This proves the commutativity of (4.6.2).
\par
Next we show that the map $f:(\Om_s)_{F'} \to (\bar A_{\la})_F$ 
is compatible with the extension of the field, i.e., the following 
diagram commutes.
\begin{equation*}
\tag{4.6.4}
\begin{CD}
(\Om_s^{F'})\wg @>f>> (\bar A_{\la})_F \\
@AAA                    @AAA \\
(\Om_s^{F^k})\wg @>f^0>> (\bar A_{\la})_{F^k}, 
\end{CD}
\end{equation*}
where $f^0$ is a similar map as $f$ defined by replacing $F'$ by $F^k$.
\par
In fact, we consider the following diagram
\begin{equation*}
\tag{4.6.5}
\begin{CD}
\wt G^F/G^F @>\pi_1>> Z_F @>\pi_2>> (\bar A_{\la})_F \\
@A N_{F^k/F}AA  @AAA         @AAA  \\
\wt G^{F^k}/G^{F^k} @>\pi_1^0>> Z_{F^k} @>\pi_2^0>> (\bar A_{\la})_{F^k}, 
\end{CD}
\end{equation*}
where $\pi_2\circ\pi_1 = f_2$, $\pi_2^0\circ\pi_1^0 = f_2^0$, and 
the second and the third vertical maps are natural
surjections.  This diagram turns out to be commutative. 
In order to show this, it is enough to see the commutativity of the
left square.  Take $g \in \wt G^F$ and write it as $g = g_1z$ with 
$g_1 \in G, z \in Z_{\wt G}$.
Then one can find $\b \in G, \g \in Z_{\wt G}$ such that 
$g_1 = F^k(\b)\b\iv$, $z = F^k(\g)\g\iv$.  Put 
$\hat g_1 = \b\iv F(\b), \hat z = \g\iv F(\g)$.  Since 
$\g \in Z_{\wt G}$, we see that $\hat g = \hat g_1\hat z$ 
satisfies the condition that $N_{F^k/F}(\hat g) = g$. 
Hence we have $\pi_1(g) = g_1\iv F(g_1) = z F(z\iv)$ and 
$\pi_1^0(\hat g) = \hat z F^k(\hat z\iv)$.  But by using 
$\hat z = \g\iv F(\g)$, we see easily that 
$zF(z\iv) = \hat z F^k(\hat z\iv)$.  This shows that 
the diagram (4.6.5) is commutative. 
\par
On the other hand, since the map 
$f_2^*: \Om_s^{F'} \to Z_{\wt G^*}^F$ is compatible with the 
inclusions $\Om_s^{F'} \hra \Om_s^{F^k}, 
               Z_{\wt G^*}^F \hra Z_{\wt G}^{F^k}$,
we have the commutative diagram
\begin{equation*}
\tag{4.6.6}
\begin{CD}
\wt G^F/G^F @>>> (\Om_s^{F'})\wg \\
@A N_{F^k/F}AA                 @AAA  \\
\wt G^{F^k}/G^{F^k} @>>> (\Om_s^{F^k})\wg.
\end{CD}
\end{equation*}
The commutativity of the diagram (4.6.4) follows from (4.6.5) and
(4.6.6).
\par
Finally, it is easy to check that the maps $h', h''$ are compatible
with $(\Om_s)^{F'} \hra \Om_s$ and 
$(\bar A_{\la}\wg)^F \hra \bar A_{\la}\wg$, and the corresponding
results hold also for $f',f''$. 

\par
Now by using the commutativity of $h,h',h''$ and $f,f',f''$, one can
check that (4.6.1) is reduced to showing that 
\par\medskip\noindent
(4.6.7) \ Assume that $F(s) = s$, and that $F$ acts trivially 
on $\Om_s$ and on $A_{\la}$.  Then, 
for $(\e,z) \in \Om_s\wg \times \Om_s$ and 
$(c,\xi) \in \bar A_{\la} \times \bar A_{\la}\wg$, we have
$\e(z) = \xi(c)$.
\par\medskip
But (4.6.7) is proved in a similar way as the proof of Lemma 3.16
in [ShS], where a similar problem for $(\CT^{(m)}_{s,E})^{F^2}$
is discussed.  Thus Lemma 4.6 is proved.
\end{proof} 
\par
The following result gives a description of the Shintani descent 
of $G^{F^m}$ in the case where $m$ is sufficiently divisible.
In the following, we denote by $\r_x^{(m)}$ the $F$-stable 
irreducible character
of $G^{F^m}$ belonging to the set ($\CT_{s, E}^{(m)})^F$ 
corresponding to $x \in \CM_{s, E}$. 
\begin{thm}[{[S2, Theorem 4.7]}]   
Assume that $m$ is sufficiently divisible.  
For each $\r_x^{(m)} \in (\Irr G^{F^m})^F$ corresponding to 
$x  = (\e,z) \in \CM_{s, E}$, we fix an extension $\td\r_x^{(m)}$ of 
$\r_x^{(m)}$ to $G^{F^m}\lp\s\rp$.  Then we have
\begin{equation*}
Sh_{F^m/F}(\td\r_x^{(m)}|_{G^{F^m}\s}) = \mu_xR_x, 
\end{equation*}
where $\mu_x$ is a certain rot of unity.  In the case where
$(s,E)$ is in (a) of 2.8, $\mu_x$ is given by an $m$-th root of unity
of $\th(\dc\iv)\xi(c\iv)$ under the correspondence 
$(\e, z) \lra (c,\xi)$ (see 2.2 and 4.5 for the notation). 
\end{thm}
\remark{4.8.}
$\mu_x$ is not given explicitly in [S2].  But the determination of 
$\m_x$ is reduced to the case where $(s, E)$ is in (a) of 2.8.
In this case the extension $\td\r^{(m)}_x$ of $\r_x^{(m)}$ is 
determined by the extension $\wt\vG^{(m)}_{c,\xi,\th}$ of 
$\vG^{(m)}_{c,\xi,\th}$ for 
$x = (c,\xi) \in \CM_{s, N}$, which is determined by the choice of
$\mu_{c,\th\xi}$ as in 2.2.  Then the argument in 4.11 in [S2]
gives the description of $\mu_x = \mu_{c,\th\xi}$.
\para{4.9.}
Let $L$ be a Levi subgroup of a standard parabolic subgroup $P$
of $G$ containing $T$.    
Let $\d = \d^{(m)}$ be an irreducible cuspidal 
character of $L^{F^m}$.  Let  $\CW = N_G(L)/L$, and put 
\begin{align*}
\CW_{\d} &= \{ w \in \CW \mid {}^w\d = \d \}, \\
\CZ_{\d} &= \{ w \in \CW \mid {}^{Fw}\d = \d \}.
\end{align*}   
$\CW_{\d}$ is naturally regarded as a subgroup of $W$, and 
according to Howlett and Lehrer [HL], $\CW_{\d}$ can be decomposed as 
$\CW_{\d} = \CW_{\d}^0\Om_{\d}$, where $\CW_{\d}^0$ is a normal subgroup 
of $\CW_{\d}$ which is a reflection group with a set of simple 
reflections associated to some root system $\vG \subset \vS$, and 
$\Om_{\d}$ is given by 
\begin{equation*}
\Om_{\d} = \{ w \in \CW_{\d} \mid w(\vG^+) \subset \vG^+\},
\end{equation*}
where $\vG^+ = \vG \cap \vS^+$ is the set of positive roots of 
$\vG^+$.  Assume that $\CZ_{\d} \ne \emptyset$. Then 
$\CZ_{\d}$ can be written as 
$\CZ_{\d} = w_{\d} \CW_{\d}$ for some $w_{\d} \in \CW$.
We choose $w_{\d}$ so that $Fw_{\d}(\vG^+) \subset \vG^+ $, 
and let $\dw_{\d} \in N_G(L)$ be a representative of $w_{\d}$.
Note that this condition determines $w_{\d}$ only up to 
the coset of $\Om_{\d}$.
Let $\g_{\d}: \CW_{\d} \to \CW_{\d}$ be the automorphism induced 
by the map $Fw_{\d}$.  Then $\g_{\d}$ stabilizes $\CW_{\d}^0$.
Let $\wt \CW_{\d} = \CW_{\d}\lp\g_{\d}\rp$ be the semidirect product
of $\CW_{\d}$ with the cyclic group generated by $\g_{\d}$. 
We denote by $(\CW_{\d}\wg)^{\g_{\d}}$ the set of  
$\g_{\d}$-stable irreducible characters of $\CW_{\d}$.
\par
Let $\CP_{\d} = \Ind_{P^{F^m}}^{G^{F^m}}\d$ be the Harish-Chandra
induction of $\d$.  We review the results from [S2, 3.5, 3.6]. 
The irreducible characters of $G^{F^m}$ appearing in the
decomposition of $\CP_{\d}$ are parametrized by $\CW_{\d}\wg$.  We denote
by $\r_E = \r^{(m)}_E$ the irreducible character of $G^{F^m}$ corresponding to 
$E \in \CW_{\d}\wg$. 
Let $M$ be the subgroup of $N_G(L)$ generated by $L$ and $w \in
\CW_{\d}$.
Then it is known ([G], [Le]) that $\d$ can be extended 
to a representation $\wt\d$ of $M^{F^m}$.  $F\dw_{\d}$ stabilizes
$M^{F^m}$, and the restriction of $F\dw_{\d}$ on $M^{F^m}$ is written
as $\s\dw_{\d}$.  $F\dw_{\d}$ stabilizes $\wt\d$, and one can extend
$\wt\d$ to a representation of $M^{F^m}\lp\s\dw_{\d}\rp$.  We fix such 
an extension of $\wt\d$, and denote it also by $\wt\d$.
\par
Now we have an action of $F$ on $\CP_{\d}$.  
$\r_E$ is $F$-stable if and only if 
$E \in (\CW_{\d}\wg)^{\g_{\d}}$.
The choice of an extension $\wt E$ of $E$ to $\wt \CW_{\d}$-module 
(and of $\wt\d$)  
determines an extension of $\r_E$ to $G^{F^m}\lp\s\rp$, 
which we denote by $\wt\r_{\wt E}$.
We consider the Shintani descent of $\wt\r_{\wt E}$.  Then by 
Theorem 3.4, one can write 
$Sh_{F^m/F}(\wt\r_{\wt E}|_{G^{F^m}\s}) = \mu_{\wt E}R_E$, where 
$R_E$ is a certain almost character of $G^F$, and $\mu_{\wt E}$ is 
a root of unity depending on the choice of $\wt E$.  
Similarly, for each $y \in W_{\d}$, $\wt\d$ is a character of 
$L^{F^m}\lp\s\dw_{\d}\dy\rp$.  
Hence the Shintani descent of $\wt\d$ can be written as  
\begin{equation*}
Sh_{F^m/F\dw_{\d}\dy}(\wt\d|_{L^{F^m}\s\dw_{\d}\dy}) 
   = \mu_{\wt\d, y}R_{\d, y}
\end{equation*} 
where $R_{\d, y}$ is the almost character of $L^{F\dw_{\d}\dy}$, 
and $\mu_{\wt\d, y}$ is a root of unity depending on the choice 
of $\wt\d$ and on $y$.
\par
Now the twisted induction 
$R_{L(\dw)}^G: C(L^{F\dw}\ssim) \to C(G^F\ssim)$ is defined as in 
[S2, 3.1].
By using the specialization argument of the Shintani descent 
identity (see [S2, Remark 4.13]), we obtain the following.
\begin{prop}  
For each $w = w_{\d}y \in \CZ_{\d}$, we have 
\begin{equation*}
R_{L(\dw)}^G(\mu_{\wt\d, y}R_{\d, y}) 
    = \sum_{E \in (\CW_{\d}\wg)^{\g_{\d}}}
                \Tr(\g_{\d}y, \wt E)\mu_{\wt E}R_E. 
\end{equation*}
\end{prop}
\remark{4.11.}
The formula in [S2, Remark 4.13] contains a linear character 
$\ve : \CW_{\d} \to \{ \pm 1\}, y \mapsto \ve_y$ 
which is trivial on $\Om_{\d}$.  However, we have $\ve = 1$
in our case.  In fact, since $\ve$ is a character of 
$\CW_{\d}^0$, $\ve$ is determined by the corresponding 
formula for $R_{\wt L(w)}^{\wt G}$ with $y \in \CW_{\d}^0$.
In that case, the formula is nothing but the decomposition of 
the Deligne-Lusztig character $R_{\wt T_x}^{\wt G}(\th)$ 
into irreducible characters for some $\th \in (\wt T_x^F)\wg$, and 
the assertion is verified by using the explicit description in 
1.4.
\section{Unipotently supported functions}
\para{5.1.}
Let $G\uni$ be the unipotent variety of $G$.
Let $\CI_G$ be the set of all pairs $(C, \CE)$ where 
$C$ is a unipotent class in $G$ and $\CE$ is an
irreducible $G$-equivariant local system on $C$.
If we fix $u \in C$, the set of $G$-equivariant local systems 
on $C$ is in bijection with $A_G(u)\wg$.  Thus the pair
$(C,\CE)$ is represented by the pair $(u, \t)$ for $\t \in A_G(u)\wg$.
Let $\CM_G$ be the set of triples $(L, C_0, \CE_0)$,
 up to $G$-conjugacy, where $L$ is a Levi subgroup of some 
parabolic subgroup $P$ of $G$, and $\CE_0$ is a cuspidal local system 
on a unipotent class $C_0$ in $L$.
It is known by Lusztig [L2, 6.5] that there exists a natural bijection
\begin{equation*}
\tag{5.1.1}
\CI_G \simeq \coprod_{(L, C_0, \CE_0) \in \CM_G}(N_G(L)/L)\wg, 
\end{equation*}
which is called the generalized Springer correspondence between
unipotent classes in $G$ and irreducible characters of various
Coxeter groups.  (Note that $N_G(L)/L$ is a Coxeter group for
any $(L, C_0, \CE_0) \in \CM_G$.) The set $\CM_G$ gives 
a partition of $\CI_G$.  A subset of $\CI_G$ corresponding to
some triple $(L,C_0,\CE_0) \in \CM_G$ is called a block.
The correspondence in (5.1.1) is given more precisely 
as follows.  For each triple
$(L, C_0, \CE_0)$, one can associate a semisimple perverse
sheaf $K$ on $G$ such that $\End K \simeq \Ql[\CW]$ with 
$\CW = N_G(L)/L$.  Let
$K_E$ be the simple component of $K$ corresponding to 
$E \in \CW$.  Then
\begin{equation*}
\tag{5.1.2}
K_E|_{G\uni} = \IC(\ol C, \CE)[\dim C + \dim Z_L^0]
\end{equation*}
for some pair $(C, \CE) \in \CI_G$.
The correspondence $(C, \CE) \lra E$ gives the required 
bijection.
\par
Now $F$ acts naturally on $\CI_G$ and $\CM_G$ by 
$(C, \CE) \mapsto (F\iv(C), F^*\CE), 
(L, C_0, \CE_0) \mapsto (F\iv(L), F\iv(C_0), F^*\CE_0)$.
Let $(L, C_0, \CE_0) \in \CM_G^F$ and $\CI_0$ the block
corresponding to it.  
Then one can choose
$L$ an $F$-stable Levi subgroup of an $F$-stable parabolic 
subgroup $P$ of $G$. In that case, $C_0$ is an $F$-stable unipotent
class and $\CE_0$ is an $F$-stable local system. 
Then $F$ acts on $\CW$, and we consider
the semidirect product $\wt \CW = \CW\lp c\rp$ , where $c$ is 
the automorphism on $\CW$ induced by $F$.
For each $\io = (C, \CE) \in \CI_0$, 
we put $K_{\io} = K_E$ if $\io = (C, \CE)$ corresponds to 
$E \in \CW\wg$ under (5.1.1).  
Then  $K_{\io}$ is $F$-stable if and 
only if $E$ is $F$-stable.  
We choose an isomorphism $\f_0: F^*\CE_0 \isom \CE_0$ so that 
it induces a map of finite order at the stalk of each point
in $C_0^F$. Then it induces an isomorphism $F^*K \isom K$, and by choosing
a preferred extension of $E$ to $\wt\CW$, induces an isomorphism 
$\f_E: F^*K_E \isom K_E$.
Since $\CH^{a_0}(K_E)|_C = \CE$ with $a_0 = -\dim C - \dim Z_L^0$,
$\f_E$ induces an isomorphism $F^*\CE \isom \CE$.
We define $\p_{\io}: F^*\CE \isom \CE$ by the condition that  
$q^{(a_0+r)/2}\p_{\io}$ coincides with the map 
$\f_E: F^*\CH^{a_0}(K_E) \isom \CH^{a_0}(K_E)$, where
$r = \dim\supp K_E$.  Note that we have 
\begin{equation*}
a_0 + r = (\dim G - \dim C) - (\dim L - \dim C_0). 
\end{equation*}
Then by [L3, 24.2], $\p_{\io}$ induces a map of finite order 
at the stalk of each point in $C^F$.
\para{5.2.}
Let $\CV_G = C(G^F\ssim)$ be the space of $G^F$-invariant 
functions on $G^F$, and $\CV\uni$ the subspace of $\CV_G$ 
consisting of functions whose supports lie in $G^F\uni$. 
For each pair $\io = (C, \CE) \in \CI_G$, we define 
$\CY_{\io} \in \CV^F\uni$ by 
\begin{equation*}
\CY_{\io}(v) = \begin{cases}
                   \Tr(\p_{\io}, \CE_v) &\quad\text{ if } v \in C^F \\
                   0                    &\quad\text{ otherwise},
               \end{cases}
\end{equation*}  
where $\CE_v$ is the stalk of $\CE$ at $v$.  
Then $\{ \CY_{\io} \mid \io \in \CI_G \}$ gives rise to a basis 
of $\CV\uni$. 
We have a natural decomposition 
\begin{equation*}
\tag{5.2.1}
\CV\uni = \bigoplus_{\CI_0}\CV_{\CI_0},
\end{equation*}
where $\CI_0$ runs over all the $F$-sable blocks, and  
$\CV_{\CI_0}$ is the subspace of $\CV\uni$ 
spanned by $\CY_{\io}$ for $\io \in \CI_0$.
\par
Let $\CI_0$ be an $F$-stable block associated to the triple 
$(L,C_0,\CE_0)$.
We assume that $L$ is an $F$-stable Levi subgroup of an $F$-stable
parabolic subgroup of $G$.  We denote by $L_w$ an $F$-stable
Levi subgroup twisted by $w \in \CW = N_G(L)/L$. 
  For a pair $\io = (C, \CE) \in \CI_G$, we put
$\supp(\io) = C$.  For each $\io,\io' \in \CI_0^F$, put
\begin{equation*}
\tag{5.2.2}
\begin{split}
\w_{\io,\io'} = &|\CW|\iv q^{-(\codim C + \codim C')/2 + \dim Z_L^0} \\ 
                &\times |G^F|\sum_{w \in \CW}|Z_{L_w}^{0F}|\iv
                     \Tr(w, E_{\io})\Tr(w, E_{\io'}),
\end{split}
\end{equation*}
where $C = \supp(\io), C' = \supp(\io')$, and 
$E_{\io}, E_{\io'} \in \CW\wg$ are the ones 
corresponding to 
$\io, \io'$ via the generalized Springer correspondence.  
If $\io, \io' \in \CI_G^F$ are not in the same block, we put
$\w_{\io,\io'} = 0$.
\par
For $K_{\io} = K_E$, put $\f_{\io} = \f_E$.  
We define $\CX_{\io} \in \CV\uni$ by 
\begin{equation*}
\CX_{\io}(g) = \sum_a(-1)^{a+a_0}\Tr(\f_{\io}, \CH^a_g(K_{\io}))
                   q^{-(a_0+r)/2} 
\qquad (g \in G^F\uni).
\end{equation*}
\par
We define an
equivalence relation $\sim$ in $\CI_G$ by $\io \sim \io'$ if 
$\supp\io = \supp\io'$.
Also we define a partial order on $\CI_G$ by $\io \le \io'$ if 
$\ol{\supp \io} \subseteq \ol{\supp\io'}$.
Assume that $\io \in \CI_0^F$.  Then it is known
that $\CX_{\io}$ can be written as 
\begin{equation*}
\tag{5.2.3}
\CX_{\io} = \sum_{\io' \in \CI_0}P_{\io',\io}\CY_{\io'},
\end{equation*}
where $P_{\io', \io} \in \BZ$.  Actually there exists a polynomial
$\BP_{\io',\io}(t) \in \BZ[t]$ such that 
$P_{\io',\io} = \BP_{\io',\io}(q)$.  Moreover, 
$P_{\io',\io} = 0$ if $\io' \not\le \io$ or if 
$\io' \sim \io, \io' \ne \io$.  $P_{\io,\io} = 1$.
In particular, $\{ \CX_{\io} \mid \io \in \CI_0^F \}$ gives rise
to a basis of $\CV_{\CI_0}$.
Moreover, we have
\begin{equation*}
\tag{5.2.4}
\lp \CX_{\io}, \CX_{\io'}\rp_{G^F} = |G^F|\iv \w_{\io,\io'}.
\end{equation*}
\para{5.3.}
Take $(L,C_0,\CE_0) \in \CM^F$.  In our case (i.e., $G$ is given
as in 1.1),  $C_0$ is the 
regular unipotent class in $L_0$. 
We choose $u_0 \in C_0^F$ as Jordan's normal form, 
and define $\f_0: F^*\CE_0 \isom \CE_0$ by the condition 
that it induces the identity map on the stalk at $u_0$.
Let  
$\CI_0$ be the block corresponding to $(L, C_0, \CE_0)$.
For $\io = (C, \CE) \in \CI_0^F$, we fix $u_1 \in C^F$ in 
Jordan's normal form.  $A_G(u_1)$ is abelian, on which 
$F$ acts naturally. The set of $G^F$-conjugacy classes 
in $C^F$ is in bijective correspondence with the group
$A_G(u_1)_F$.  We denote by $u_a$ a representative of 
the $G^F$-class in $C^F$ corresponding to $a \in A_G(u_1)_F$.   
Assume that $\CE$ corresponds to an $F$-stable
irreducible character $\t \in A_G(u_1)\wg$.  We define 
a function $\x_{u_1,\t} \in \CV\uni$ by   
\begin{equation*}
\x_{u_1, \t}(g) = \begin{cases}
                 \t(u_a) &\quad\text{ if } g \sim_{G^F} u_a, \\
                  0      &\quad\text{ if } g \notin C^F.
             \end{cases}
\end{equation*}
\par
The following result determines the function $\CY_{\io}$ explicitly. 
\begin{prop}[{[S3]}]  
Assume that $\io = (C, \CE)$ is represented by 
$(u_1,\t)$ as above.  Then we have
$\CY_{\io} = \x_{u_1, \t}$.
\end{prop}
\para{5.5.}
By making use of the map $\log: G\uni \to \Fg\nil$ (see 2.1), we
identify $\CV\uni$ the space of $G^F$-invariant functions 
$\Fg\nil^F$.  Then the function $\CY_{\io}$ can be regarded as 
a function on $\CO^F$, where $\CO$ is the nilpotent orbit 
corresponding to $C$ such that $C = \supp (\io)$.
For each $F$-stable nilpotent orbit $\CO$, we choose a 
representative $N \in \CO^F$ via Jordan's normal form corresponding 
to $u_1 \in C^F$.  Let $\{ N, N^*, H\}$ be the TDS-triple.
Then the associated parabolic subgroup $P_N$ and its Levi subgroup
$L_N$ are defined, and we have the group 
$A_{\la} = Z_{L_N}(\la)/Z^0_{L_N}(\la)$ as in 2.1. 
For $c \in (A_{\la})_F$, we consider the twisted element
$N_c$.  Then the generalized Gelfand-Graev character
$\vG_c$ associated to $N_c$ is defined as in 2.1, which 
gives an element of $\CV\uni$. 
Now Lusztig gave 
a formula
expressing $\vG_c$ in terms of the linear combination of $\CX_{\io}$
as follows.
\begin{thm}[{[L7, Theorem 7.3]}]  
Let $\CI_0$ be an $F$-stable block corresponding to 
$(L, C_0, \CE_0)$.  Let $(\vG_c)_{\CI_0}$ be the 
projection of $\vG_c$ onto the subspace $\CV_{\CI_0}$
in (5.2.1). Then 
\begin{equation*}
\tag{5.6.1}
\begin{split}
(\vG_c)_{\CI_0} = \sum_{\io, \io',\io_1 \in \CI_0}q^{f(\io, \io_1)}
        \z_{\CI_0}\iv&
         |\CW|\iv\sum_{w \in \CW}
             \Tr(w, E_{\io})\Tr(w, E_{\io_1}\otimes\ve)  \\
         &\times |Z^{0F}_{L_w}|\BP_{\io',\io}(q\iv)
            \ol{\CY_{\io'}(-N_c^*)}\CX_{\io_1},
\end{split}
\end{equation*}
where
\begin{equation*}
\begin{split}
f(\io,\io_1) = &-\dim\supp(\io_1)/2 + \dim\supp(\io)/2  \\
               &-\dim \CO_N/2 + \dim (G/Z^0_L)/2, 
\end{split}
\end{equation*}
and $\z_{\CI_0}$ is a fourth root of unity attached to the block $\CI_0$.
$\ve$ is the sign representation of $\CW$ (cf. [L7, 5.5]).
\end{thm} 
\remark{5.7.} \ 
The restriction of the Fourier transform of $\CX_{\io}$ 
($\io \in \CI_0$) on $\Fg\nil$ coincides with $\CX_{\io}$ 
up to scalar.  The fourth root of unity $\z_{\CI_0}$ occurs in the
description of this scalar ([L7, Proposition 7.2]).
$\z_{\CI_0}$ depends only on the
pair $(C_0, \CE_0)$ and does not depend on $G$. 
In our case, $\CI_0$ is always a regular block, i.e., $C_0$ is 
the regular unipotent class in $L$.  In such a case, 
Digne, Lehrer and Michel [DLM1, Proposition 2.8] determined the value
$\z_{\CI_0}$ explicitly.
\para{5.8.}
In order to apply the formula (5.6.1), we need to describe 
$-N_c^*$ for a nilpotent element $N \in \Fg^F$.
Since $-N^* \in \Fg^F$ is $G$-conjugate to $N$, one can write 
$-N^* = N_{c_0}$ for some 
$c_0 \in (A_{\la})_F$, i.e., $N_{c_0} = \Ad(\a_{c_0})N$ with 
$\a_{c_0}\iv F(\a_{c_0}) = \dc_0$.  We consider $N_c$ for 
$c \in (A_{\la})_F$. Then $N_c = \Ad(\a_c)N$ with 
$\a_c\iv F(\a_c) = \dc$.  Since $Z_G \to A_{\la}$ is surjective, 
we may choose $\dc \in Z_G$.  Now 
$-N^*_c$ is obtained as 
$-N^*_c = \Ad(\a_c)(-N^*)$.  Hence $-N^*_c$ is $G^F$-conjugate
to $\Ad(\a_c\a_{c_0})N$.  But since $\dc \in Z_G$, we have
\begin{equation*}
(\a_c\a_{c_0})\iv F(\a_c\a_{c_0}) 
  = \a_{c_0}\iv \dc F(\a_{c_0}) 
 = \dc\dc_0,
\end{equation*}
It follows that 
$-N_c^*$ is $G^F$-conjugate to $N_{cc_0}$.
\para{5.9.}
Following [L2, LS], we describe the generalized Springer
correspondence for $G$ explicitly.
Let $n'$ be the largest common divisor of $n_1, \dots, n_r$ which is 
prime to $p$.  Then $Z_G$ is the cyclic group of order $n'$.  
Let $u = u_{\mu}$ 
be a unipotent element in $G$ corresponding to 
$\mu = (\mu_1, \dots, \mu_r)$, where 
$\mu_i = (\mu_{i1} \ge \mu_{i_2} \ge\cdots)$ is a
partition of $n_i$. 
Put $n_{\mu}'$ be the greatest common divisor of $n'$ and 
$\{\mu_{ij}\}$.  Then $A_G(u)$ is a cyclic group of order 
$n_{\mu}'$.
For each $\t \in A_G(u)\wg$, $Z_G/Z_G^0$ acts on the representation 
space $V_{\t}$ of $\t$ via the homomorphism $Z_G/Z_G^0 \to A_G(u)$. 
For each $\e \in (Z_G/Z_G^0)\wg$, we denote by $A_G(u)\wg_{\e}$ the set of
irreducible characters $\t \in A_G(u)\wg$ such that $Z_G/Z_G^0$ acts 
on $V_{\t}$ via the character $\e$.  
We have 
\begin{equation*}
\tag{5.9.1}
|A_G(u)\wg_{\e}| = \begin{cases}
                       1 &\quad\text{ if } d|n_{\mu}', \\
                       0 &\quad\text{ otherwise,}
                    \end{cases}
\end{equation*}
where $d$ is the order of $\e$.
Now the generalized Springer correspondence in (5.1.1) is 
described as follows:  We have a partition
\begin{equation*}  
\CI_G = \coprod_{\e \in (Z_G/Z_G^0)\wg}(\CI_G)_{\e},
\end{equation*}
where $(\CI_G)_{\e}$ is the set of pairs $(u,\t)$ with 
$\t \in A_G(u)\wg_{\e}$.  Note that $\t$ is uniquely determined 
by $u$ if $(u,\t) \in (\CI_G)_{\e}$ by (5.9.1), which we denote by
$\t(u)$. 
For each $\e \in (Z_G/Z_G^0)\wg$ of order $d$, there exists
a unique Levi subgroup up to conjugacy such that the type of $L$
is $A_{d-1} + \cdots + A_{d-1}$, and a unique cuspidal pair
$(L, C_0, \CE_0)$.  Here  $C_0$ is regular unipotent in $L$ and
for $u_0 \in C_0$, $A_L(u_0) \simeq Z_L/Z_L^0$.  $\CE_0$ is the
unique local system on $C_0$ corresponding to 
$\e_0 \in (Z_L/Z_L^0)\wg$ such that $\e_0\circ f = \e$ for 
a natural homomorphism
$f : Z_G/Z_G^0 \to Z_L/Z_L^0$.  Then 
$N_G(L)/L \simeq \FS_{n_1/d} \times\cdots\times \FS_{n_r/d}$, 
and the map $E_{\mu} \mapsto (u_{d\mu}, \t(u_{d\mu}))$ 
($d\mu = (d\mu_{ij})$ for 
$\mu = (\mu_{ij})$) gives the generalized Springer 
correspondence 
\begin{equation*}
\tag{5.9.2}
(N_G(L)/L)\wg \simeq (\CI_G)_{\e}.
\end{equation*}
\para{5.10.}
Assume that $\wt G = \wt G_1 \times \cdots\times \wt G_r$ with 
$n_1 = \cdots = n_r = t$.
In this case, $n'$ is the largest divisor of $t$ which is prime to   
$p$. From the description of the generalized Springer correspondence,
the partition of $\CI_G$ into blocks is nothing but the partition of 
$\CI_G$ into  $(\CI_G)_{\e}$.  Assume that $\CI_0 = (\CI_G)_{\e}$ 
with $\e \in (Z_G/Z_G^0)\wg$ of order $d$.  We shall make the formula
(5.6.1) more explicit in the case where $N$ is regular nilpotent,
i.e., where $\vG_{c}$ is the modified Gelfand-Graev characters.
\begin{lem}  
Let $G$ be as before.  Assume that $N$ is regular nilpotent. 
Then for any $\io \in \CI_0 = (\CI_G)_{\e}$, 
\begin{equation*}
\lp\, (\vG_c)_{\CI_0}, \CX_{\io}\rp_{G^F} = 
                 \begin{cases}
        q^{(\dim Z^0_L - \codim\supp(\io))/2}\z_{\CI_0}\iv\e(cc_0)\iv
            &\text{ if } E_{\io} = \ve, \\
        0   &\text{ otherwise.}
                 \end{cases}                 
\end{equation*} 
\end{lem}
\begin{proof}
We apply the formula (5.6.1).  Since $-N^*_c$ is regular nilpotent,  
the non-zero contribution of $\CY_{\io'}(-N_c^*)$ occurs only when  
$\supp(\io') = \CO_N$, the regular nilpotent orbit.
By 5.8, $-N_c^*$ is $G^F$-conjugate to $N_{cc_0}$.  Since 
$A_G(N) \simeq Z_G/Z_G^0$, we see that $\io' = (N, \e)$ under the 
identification $\e \in A_G(N)\wg$.  This implies that 
$\CY_{\io'}(-N_c^*) = \e(cc_0)$ by Proposition 5.4.
On the other hand, since $\supp(\io') = \CO_N$, 
$P_{\io',\io} = 0$ unless $\io = \io'$, and we have 
$P_{\io',\io'} = 1$ by the property in 5.2.
Moreover, under the generalized Springer correspondence, 
$E_{\io'}$ is the identity character of $\CW$.
Thus (5.6.1) can be written as 
\begin{equation*}
(\vG_c)_{\CI_0} = \e(cc_0)\iv\sum_{\io_1 \in \CI_0}q^{f(\io',\io_1)}
            \z_{\CI_0}\iv|\CW|\iv\sum_{w \in \CW}
                 |Z_{L_w}^{0F}|\Tr(w, E_{\io_1}\otimes\ve)\CX_{\io_1}.
\end{equation*}
Now by (5.2.2) and (5.2.4), we have
\begin{align*}
\lp\, (\vG_c)_{\CI_0}, \CX_{\io}\rp_{G^F}  
    &= \z_{\CI_0}\iv\e(cc_0)\iv|\CW|^{-2}\sum_{\io_1 \in \CI_0}
          q^{(\dim Z^0_L - \codim\supp(\io))/2} \\
    &\times\sum_{w, w' \in \CW}|Z_{L_w}^{0F}||Z_{L_{w'}}^{0F}|\iv
     \ve(w)\Tr(w, E_{\io_1})\Tr(w', E_{\io_1})\Tr(w', E_{\io}) \\
    &= q^{(\dim Z^0_L - \codim\supp(\io))/2}\z_{\CI_0}\iv
       \e(cc_0)\iv|\CW|\iv\sum_{w \in \CW}\ve(w)\Tr(w, E_{\io}).
\end{align*}
The lemma follows from this.
\end{proof}
For later applications, we also consider the general case where 
$N$ is arbitrary. 
\begin{lem} 
Assume that $N$ is an arbitrary nilpotent element.  
\begin{enumerate}
\item
For any 
$\io \in \CI_0^F$, we have
\begin{equation*}
\lp\vG_c, \CX_{\io}\rp_{G^F} 
  = q^{g(\io',\io)}\z_{\CI_0}\iv\BP_{\io'',\io'}(q\iv)\ol{\CY_{\io''}(-N_c^*)},
\end{equation*}
where 
\begin{equation*}
g(\io',\io) = (-\codim_G\supp(\io') + \dim Z_L)/2
                   +(\dim\supp(\io) - \dim\CO_N)/2
\end{equation*}
and $\io'\in \CI_0$ is such that $E_{\io} = E_{\io'}\otimes\ve$, 
$\io''$ is the unique element in $\CI_0$ 
such that $\supp(\io'') = \CO_N$.
\item
$t^{(\dim\supp(\io') - \dim\supp(\io''))/2}\BP_{\io'',\io'}(t\iv)$ 
is a polynomial in $t$.  It is divisible by $t$ if $\io'' \ne \io'$.
\item
There exists the ring of integers $\CA$ of some fixed cyclotomic field 
independent of the field $\Fq$ such that 
\begin{equation*}
\lp\vG_c, \CX_{\io}\rp_{G^F} \in q^{(\dim Z_L - \codim_G\supp(\io))/2}\CA.
\end{equation*}
\end{enumerate}
\end{lem}
\begin{proof}
First we show (i).  By applying Theorem 5.6, we have
\begin{align*}
\lp\vG_c, \CX_{\io}\rp_{G^F} 
  &= \lp\ (\vG_c)_{\CI_0}, \CX_{\io}\rp_{G^F}  
= \sum_{\io',\io'', \io_1}q^{f(\io',\io_1)}\z_{\CI_0}\iv|\CW|\iv \times \\
   &\times \sum_{w \in \CW}\Tr(w, E_{\io'})\Tr(w, E_{\io_1})\ve(w)
|Z_{L^w}^{0F}|\BP_{\io'',\io'}(q\iv)
   \ol{\CY_{\io''}(-N_c^*)}\lp\CX_{\io_1}, \CX_{\io}\rp_{G^F}.
\end{align*}
Substituting (5.2.4) into $\lp\CX_{\io_1}, \CX_{\io}\rp_{G^F}$, 
we have
\begin{align*}
\lp\vG_c, \CX_{\io}\rp_{G^F} = 
\sum_{\io',\io''}&q^{g(\io',\io)}\z_{\CI_0}\iv|\CW|^{-2}  \\
&\times \sum_{w, w' \in \CW}
|Z_{L^w}^{0F}||Z_{L^{w'}}^{0F}|\iv
\biggl(\sum_{\io_1}\Tr(w, E_{\io_1})\Tr(w',E_{\io_1})\biggr)   \\
 &\times \ve(w)\Tr(w, E_{\io'})\Tr(w', E_{\io})
       \BP_{\io'',\io'}(q\iv)\ol{\CY_{\io''}(-N_c^*)} \\
= \sum_{\io',\io''}&q^{g(\io',\io)}\z_{\CI_0}\iv 
     \biggl(|\CW|\iv\sum_{w \in \CW} 
        \Tr(w, E_{\io'}\otimes\ve)\Tr(w,E_{\io})\biggr) \\
          &\times\BP_{\io'',\io'}(q\iv)\ol{\CY_{\io''}(-N_c^*)}
\end{align*}
with $g(\io',\io)$ as in (i). 
Hence in the sum, the only $\io'$ such that 
$E_{\io} = E_{\io'}\otimes \ve$ gives the contribution.
On the other hand, the condition $\CY_{\io''}(-N_c^*) \ne 0$
implies that $\supp(\io'') = \CO_N$.
It follows that 
\begin{equation*}
\lp\vG_c, \CX_{\io}\rp_{G^F} = 
q^{g(\io',\io)}\z_{\CI_0}\iv
   \BP_{\io'',\io'}(q\iv)\ol{\CY_{\io''}(-N_c^*)},
\end{equation*}
where $\io''$ is the unique elements in $(\CI_0)^F$ such that
$\supp(\io'') = \CO_N$.
This proves the first statement.
\par
Next we show (ii). 
We may assume that $G = SL_n$. By the generalized 
Springer correspondence $(\CI_G)_{\e} \lra \FS_{n/d}$, we have 
$\supp(\io') = \CO_{\la}$ and $\supp(\io'') = \CO_{\mu}$, where 
$\la, \mu$ are partitions of $n$ such that each part is divisible by 
$d$. Put $\io' = \io_{\la}, \io''= \io_{\mu}$.  
We denote by $\la/d, \mu/d$ the partitions of $n/d$ obtained
from $\la, \mu$ by dividing each part by $d$.  We consider the group
$GL_{n/d}$ and nilpotent orbits $\CO_{\la/d}, \CO_{\mu/d}$ of 
$\Fg\Fl_{n/d}$.  We have a polynomial $\BP_{\io_{\mu/d}, \io_{\la/d}}$
associated to $GL_{n/d}$, defined in a similar way as 
$\BP_{\io_{\mu}, \io_{\la}}$, where $\io_{\mu/d}, \io_{\la/d}$ are elements in 
$\CI_{GL_{n/d}}$ corresponding to $\CO_{\la/d}, \CO_{\mu/d}$ under the
Springer correspondence.  Then it is known by [DLM2, Theorem 8.1] that
\begin{equation*}
t^{(\dim \CO_{\mu} - \dim \CO_{\la})/2}\BP_{\io_{\mu}, \io_{\la}}(t)
  = t^{(\dim\CO_{\mu/d} - \dim \CO_{\la/d})/2}
          \BP_{\io_{\mu/d}, \io_{\la/d}}(t).
\end{equation*}
It follows that 
\begin{equation*}
t^{(\dim\CO_{\la} - \dim\CO_{\mu})/2}\BP_{\io_{\mu}, \io_{\la}}(t\iv)
  = t^{(\dim\CO_{\la/d} - \dim\CO_{\mu/d})/2}
         \BP_{\io_{\mu/d},\io_{\la/d}}(t\iv).
\end{equation*}
But the right hand side of this formula coincides with the 
Kostka polynomial $K_{\la/d, \mu/d}(t)$ associated to partitions 
$\la/d, \mu/d$ of $n/d$ (cf. [M, III]), hence it is a polynomial in $t$.  
Then the second assertion of the lemma follows from the well-known properties
of Kostka polynomials.
\par
Finally we show (iii).  
We may assume that $\z_{\CI_0} \in \CA$ and 
$\ol{\CY_{\io''}(-N_c^*)} \in \CA$.  Thus (iii) follows from 
(i) and (ii). 
The lemma is proved.
\end{proof}
\par\bigskip
\section{Character sheaves}
\para{6.1.}  Following [L3, IV], 
we review the classification of  
character sheaves in the case of type $A$.
So let $G$ be as before, and we denote by $\wh G$ the set
of character sheaves on $G$.  Let $\CS(T)$ be the set
of local systems of rank 1 on $T$ such that 
$\CL^{\otimes m} \simeq \Ql$ for some $m$ prime to $p$.
Then for each $\CL \in \CS(T)$, the subset 
$\wh G_{\CL}$ of $\wh G$ is defined, and we have
\begin{equation*}
\wh G = \coprod_{\CL \in \CS(T)}\wh G_{\CL}. 
\end{equation*}
\par
For each $\CL \in \CS(T)$, we put 
$W_{\CL} = \{ w \in W \mid w^*\CL \simeq \CL\}$.
Then there exists a subroot system $\vS_{\CL}$ of $\vS$, 
and $W_{\CL}$ can be written as 
$W_{\CL} = \Om_{\CL}\ltimes W_{\CL}^0$, where 
$W_{\CL}^0$ is the reflection group associated to the root
system $\vS_{\CL}$, and 
$\Om_{\CL} = \{ w \in W_{\CL} \mid w(\vS_{\CL}^+) = \vS_{\CL}^+\}$.
If $\CL \in \CS(T)$ is $F^m$-stable for some integer $m > 0$, we 
fix an isomorphism $\f: (F^m)^*\CL \isom \CL$ so that it induces an
identity map on $\CL_1$ (the stalk at $1 \in T^{F^m}$).  
Then the characteristic function 
$\x_{\CL, \f}$ gives rise to a character $\th \in (T^{F^m})\wg$, which
induces an isomorphism between $\CS(T)^{F^m}$ and 
$(T^{F^m})\wg$.  Thus we have $\CS(T)^{F^m} \simeq (T^*)^{F^m}$.
  Since this isomorphism is compatible with the extension of the 
filed $\BF_{q^m}$, we 
can identify $\CS(T)$ with $T^*$ in this way.  Now assume that 
$\CL \in \CS(T)$ corresponds to $s \in T^*$.  Then we have 
$W_{\CL} = W_s, W_{\CL}^0 = W_s^0$ and $\Om_{\CL} = \Om_s$ in the 
notation of Section 1.
\par
In [L3, 17], families in $(W^0_{\CL})\wg$ and in  $W\wg_{\CL}$ were
introduced.  In our case, $W^0_{\CL}$ is a product of symmetric
groups and a family $\CF$ in $(W^0_{\CL})\wg$ is of the form 
$\CF = \{ E \}$ with $E \in (W_{\CL}^0)\wg$.
Let $\Om_{\CL,E}$ be the stabilizer of $E$ in $\Om_{\CL}$. 
Then $E$ can be extended to a character $\wt E$ on 
$\Om_{\CL,E}W_{\CL}^0$. (We choose a canonical extension
so that for each $\s \in \Om_{\CL,E}$, $\wt E$ gives the preferred
extension of $E$ to $\lp\s\rp W^0_{\CL}$). 
For each $\th \in \Om_{\CL,E}\wg$, 
put $\wt E_{\th} = 
  \Ind_{\Om_{\CL,E}W_{\CL}^0}^{\Om_{\CL}W_{\CL}^0}
        (\th \otimes \wt E)$.  
Then $\wt E_{\th}$ is irreducible, and the family $\CF'$ in $W_{\CL}$
associated to $\CF$ consists of 
$\{ \wt E_{\th} \mid \th \in \Om_{\CL,E}\wg\}$.
Put 
\begin{equation*}
\CM_{\CL, E} = \Om_{\CL,E} \times \Om_{\CL,E}\wg.
\end{equation*}
We have an embedding $\CF' \hra \CM_{\CL, E}$  by
$\wt E_{\th} \mapsto (1, \th)$.   
\para{6.2.}
For each $\CL \in \CS(T)$ and $w \in W_{\CL}$, we choose a 
representative $\dw \in N_G(T)$.  By [L3, 2.4], a complex 
$K_{\dw}^{\CL} \in \DD G$ is defined. The set $\wh G_{\CL}$ is 
defined as the set of character sheaves $A$ such that $A$ is 
a consistuent of ${}^pH^i(K_{\dw}^{\CL})$ for some 
$w \in W_{\CL}$ and some $i \in \BZ$. Let 
$E \in (W_{\CL}^0)\wg$ and $\wt E_{\th} \in W_{\CL}\wg$ be as 
in 6.1.  We define 
\begin{equation*}
R_{\wt E_{\th}}^{\CL} = |W_{\CL}|\iv\sum_{y \in W_{\CL}}
       \Tr(y\iv, \wt E_{\th})\sum_i(-1)^{i+\dim G}
                  \, {}^p\!H^i(K_y^{\CL}),
\end{equation*}
which is an element of the subgroup of the Grothendieck group of 
the perverse sheaves on $G$ spanned by the character sheaves of $G$.
\par
By [L3, Corollary 16.7], a natural map from 
$\wh G_{\CL}$ to the set of two sided cells  of $W_{\CL}$
was constructed. We denote by $\wh G_{\CL, \CF'}$ the set of 
character sheaves $A \in \wh G_{\CL}$ such that the corresponding 
two sided cell coincides with $\CF'$.  In our case, the family 
$\CF'$ is determined by the choice of $E \in (W_{\CL}^0)\wg$. 
  Thus we write $\wh G_{\CL, \CF'}$ as $\wh G_{\CL, E}$.
We have the following partition of $\wh G_{\CL}$.
\begin{equation*}
\wh G_{\CL} = \coprod_{E \in (W_{\CL}^0)\wg/\Om_{\CL}}
                   \wh G_{\CL, E}.
\end{equation*} 
The following result gives a parametrization of $\wh G$.
\begin{prop}[{[L3, Proposition 18.5]}] 
There exists a bijection $\wh G_{\CL, E} \lra \CM_{\CL,E}, 
A \lra (x_A, \th_A)$ satisfying the following property:
for any $\th \in \Om_{\CL,E}\wg$, 
\begin{equation*}
(A: R_{\wt E_{\th}}^{\CL}) = |\Om_{\CL,E}|\iv \ve_A\th(x_A)\iv,
\end{equation*}
where $\ve_A = (-1)^{l(x_A)}$.
($l$ is the restriction of the length function of $W$ to $\Om_{\CL,E}$). 
\end{prop}
\para{6.4.}
Let $\wh G_0$ be the set of cuspidal character sheaves on $G$.
We shall describe the set $\wh G_0$ explicitly.  By Lemma 18.4 and 
by the proof of Proposition 18.5 in [L3], $\wh G_{\CL}$ contains 
a cuspidal character sheaf if and only if $W^0_{\CL} = \{ 1 \}$ and 
$\Om_{\CL} = W_{\CL}$ is a cyclic group generated by a 
Coxeter element in $W$, which is isomorphic to $Z_G/Z_G^0$.  They are
indexed by the pair $(x, z)$ where $x$ is a generator of $\Om_{\CL}$
and $z$ is a representative of $Z_G/Z_G^0$.    In particular, $\wt G$   
is of the form $\wt G = \wt G_1 \times \cdots \times \wt G_r$, where
$\wt G_i \simeq GL_t$ with $t = n/r$.  
Also we note that the character sheaf 
$A_{x,z}$ corresponding to the pair $(x,z)$ has its support
in $zZ_G^0\times G\uni$. 
Under the parametrization in Proposition 6.3, this is also 
given in the following form. 
\par\medskip\noindent
(6.4.1) \ Let $\wt G$ be as above.  Assume that $W^0_{\CL} = \{ 1\}$
        and that $\Om_{\CL}$ is a cyclic group generated by a Coxeter 
element, which is isomorphic to $Z_G/Z_G^0$.  Hence 
$\CM_{\CL, E} = \Om_{\CL} \times \Om_{\CL}\wg$.  Let 
$(\wh G_{\CL})_0$ be the set of cuspidal
character sheaves in $\wh G_{\CL}$.  Then we have
\begin{equation*}
(\wh G_{\CL})_0 = 
   \{ A_{x, \th} \mid x \in \Om_{\CL,0}, \th \in \Om_{\CL}\wg\},   
\end{equation*}
where $\Om_{\CL,0}$ is the set of $x \in \Om_{\CL}$ 
such that $x$ is a generator of $\Om_{\CL}$. 
\para{6.5.}
The set $(\wh G_{\CL})_0$ is also given in terms of intersection 
cohomology complexes as follows.  
Let $C$ be the  regular unipotent class in $G$, and we choose a 
representative $u_1 \in G^F$. For each $z \in Z_G/Z_G^0$, 
we choose a representative $\dz \in Z_G$.  
Then $\dz C$ is a conjugacy class 
of $G$ containing $\dz u_1$ and the component group 
$A_G(\dz u_1)$ coincides with $A_G(u_1) \simeq Z_G/Z_G^0$.  
We denote by $\CE_{\e}$ the irreducible local system on $\dz C$
corresponding to $\e \in A_G(\dz u_1)\wg$.  
Put $\vS = zZ_G^0C = Z_G^0\times \dz C$, and consider a local system 
$\Ql\boxtimes \CE_{\e}$ on $\vS$ associated 
to $(z, \e) \in Z_G/Z_G^0 \times A_G(\dz u_1)\wg$.
We define a perverse sheaf $A_{z,\e}$ on $\vS$ by 
\begin{equation*}
\tag{6.5.1}
A_{z,\e} = \IC(\ol\vS, \Ql\boxtimes \CE_{\e})[\dim \vS]  
       = \Ql\boxtimes \IC(\ol{\dz C}, \CE_{\e})[\dim Z_G^0 + \dim C].
\end{equation*}
Then $A_{z,\e} \in \wh G$.  Let $\CE$ be a 
local system of $G$ of rank 1 obtained as the inverse image 
under the map $G \to G/G\der$ of a local system $\CE' \in \CS(G/G\der)$. 
We have $\CE\otimes A_{z,\e} \in \wh G$. 
Let $A_G(\dz u_1)\wg_0$ be the subset of $A_G(\dz u_1)\wg$ 
consisting of faithful characters of $A_G(\dz u_1)$.
Then we have
\begin{equation*} 
\tag{6.5.2}
\wh G_0 = 
 \{ \CE\otimes A_{z, \e} \mid z \in Z_G/Z_G^0, 
             \e \in A_G(\dz u_1)\wg_0, \CE' \in \CS(G/G\der) \}.  
\end{equation*}
\par
We now consider the $\Fq$-structure of $G$, and let 
$\wh G^F$ (resp. $\wh G^F_0$) be the set of $F$-stable
character sheaves (resp. $F$-stable cuspidal character 
sheaves ) of $G$.   The regular unipotent class $C$ is 
$F$-stable, and we choose $u_1 \in C^F$ such that 
$u_1$ is given by Jordan's normal form.  Also we can 
choose a representative $\dz \in Z_G^F$.
Then $\wh G_0^F$ is given as
\begin{equation*}
\tag{6.5.3}
\wh G_0^F = 
  \{ \CE\otimes A_{z,\e} \mid z \in (Z_G/Z_G^0)^F, 
 \e \in (A_G(\dz u_1)\wg_0)^F, \CE|_{\dz Z_G^0} \text{ :$F$-stable} \}.
\end{equation*}
Put $y = (z, \e)$, and assume that $A_y = A_{z,\e}$ is $F$-stable.
We have $F^*\CE_{\e} \simeq \CE_{\e}$, and choose an 
isomorphism $\vf_0: F^*\CE_{\e} \isom \CE_{\e}$
by the requirement that $\vf_0$ induces an identity map 
on the stalk at $\dz u_1 \in (\dz C)^F$.   
$\vf_0$ induces an isomorphism $\wt\vf_0: F^*A_y \simeq A_y$.
Note that 
$\CH^{-d}(A_y)|_{zZ_G^0C} = \Ql\boxtimes\CE_{\e}$, where
$d = \dim Z_G^0 + \dim C$. 
Then we define $\f_y: F^*A_y \isom A_y$ by the condition that  
$\f_y = q^{(\dim G - d)/2}\wt\vf_0 = 
   q^{(\codim C - \dim Z_G^0)/2}\wt\vf_0$ 
on $\CH^{-d}_g(A_y)$
($g \in (zZ_G^0\ol C)^F$).
We denote by $\x_y$ the characteristic function 
$\x_{A_y,\f_y}$ associated to $A_y$ and $\f_y$. 
Now, for each $\dz \in Z_G^F$, the left multiplication 
$\dz: C \to \dz C$ induces, under the 
identification $A_G(\dz u_1) \simeq A_G(u_1)$, the isomorphism  
$\dz^*\IC(\ol{\dz C}, \CE_{\e}) \isom \IC(\ol C, \CE_{\e})$
compatible with the $\Fq$-structure.  
Recall that for the cuspidal pair 
$\io_0  = (C, \CE_{\e}) \in \CI_G$, one can attach 
the function $\CX_{\io_0} \in \CV\uni$ as in 5.2.  Then we see easily
that   
\begin{equation*}
\tag{6.5.4}
\x_y(g) = \begin{cases}
  q^{(\codim C - \dim Z_G^0)/2}\CX_{\io_0}(v) 
                       &\quad\text{ if  $g = \dz z_1v$ with 
                             $z_1 \in Z_G^{0F}, v\in \ol C^F$}, \\
          0            &\quad\text{ otherwise}.
                    \end{cases}     
\end{equation*}
\par
More generally, we consider $\CE\otimes A_y$ with 
$\CE' \in \CS(G/G\der)$ such that $\CE|_{\dz Z_G^0}$ is $F$-stable.
Put $\CE_1 = \CE|_{\dz Z_G^0}$.
The isomorphism $\vf_1: F^*\CE_1 \isom \CE_1$ is described 
as follows. Put $\CE_2 = \CE|_{T}$.  There exists an integer 
$m > 0$ such that $\CE_2$ is 
$F^m$-stable.  We choose $\vf_2: (F^m)^*\CE_2 \isom \CE_2$ so that 
it induces the identity map on the stalk at $1 \in T^{F^m}$.  Then 
the characteristic function $\x_{\CE_2, \vf_2}$ coincides with a linear
character $\p \in (T^{F^m})\wg$.  Let $V$ be a one dimensional 
$T^{F^m}$-module affording $\p$.  
We consider the quotient $T\times^{T^{F^m}}V$ of $T\times V$ by 
$T^{F^m}$ under the action $t_1: (t,v) \mapsto (tt_1\iv, t_1v)$.
Then $\CE_2$ is 
realized as the local system associated to the locally trivial
fibration $f: T\times^{T^{F^m}}V \to T$, $(t,v) \mapsto t^{q^m-1}$.
For $\dz \in T^F$, one can choose $\a \in T$ such that 
$\a^{q^m-1} = \dz$.  
Now $f\iv(\dz Z_G^0)$ can be identified with 
$\a Z_G^0 \times^{Z_G^{0F^m}}V$, and 
$\CE_1$ is the local system associated to 
$f_1: \a Z_G^0\times^{Z_G^{0F^m}}V \to \dz Z_G^0$ 
with $f_1 = f|_{f\iv(\dz Z_G^0)}$.
Since $\CE|_{\dz Z_G^0}$ is $F$-stable, the restriction of 
$\p$ on $\dz Z_G^{0F^m}$ is $F$-stable.  Since $\dz \in Z_G^F$, 
$\p|_{Z_G^{0F^m}}$ is also $F$-stable.  Hence we may assume that 
$V$ satisfies the property that $F(t)v = tv$ for
$t \in Z^{0F^m}, v\in V$.  Moreover, since $\dz \in T^F$, we have
$\hat z = \a\iv F(\a) \in T^{F^m}$.  It follows that we have 
a well-defined automorphism $F$ on $\a Z_G^0\times^{Z_G^{0F^m}}V$ given by 
\begin{equation*}
\tag{6.5.5}
F(\a t, v) = (F(\a t),v) = (\a F(t), \hat z v)
\end{equation*} 
for $(\a t, v) \in \a Z_G^0\times^{Z_G^{0F^m}}V$.
The map $f_1$ is compatible with the actions of $F$ on 
$\a Z_G^0\times^{Z_G^{0F^m}}V$ and on $\dz Z_G^0$ (the natural action).
Hence $F$ defines an isomorphism 
$\vf_1: F^*\CE_1 \isom \CE_1$.
\par
We define a character $\th_0 \in (Z_G^{0F})\wg$ by 
$\th_0 =  Sh_{F^m/F}(\p|_{Z_G^{0F^m}})$.  
Then the previous argument implies, in view of (6.5.5), that
\begin{equation*}
\tag{6.5.6}
\x_{\CE_1,\f_1}(\dz z_1) = \p(\hat z)\th_0(z_1) \quad (z_1 \in Z_G^{0F}).
\end{equation*} 
We now define an isomorphism
$\f_0 : F^*(\CE_1\otimes A_y) \isom \CE_1\otimes A_y$ by 
$\f_0 = \vf_1\otimes \f_y$, which is regarded as an isomorphism 
$F^*(\CE\otimes A_y) \isom \CE\otimes A_y$.
We denote by $\x_{\CE,y}$ the characteristic function 
$\x_{\CE\otimes A_y, \f_0}$.
It follows from (6.5.6) that we have 
\begin{equation*}
\tag{6.5.7}
\x_{\CE,y}(z_1g) = \p(\hat z)\th_0(z_1)\x_{y}(g)
  \qquad (z_1 \in Z_G^{0F}, g \in \dz\ol C^F).
\end{equation*}
\para{6.6.}
We consider the map $f: Z_G^0 \hra G \to G/G\der$.  Then 
$f(Z_G^0)$ can be identified with $Z_G^0/Z_G^0\cap G\der$, 
and $Z_G^0 \to Z_G^0/Z_G^0 \cap G\der$ is a finite \'etale
covering. It follows that we have an isomorphism 
$\CS(Z_G^0)\simeq \CS(Z_G^0/Z_G^0\cap G\der)$ and a surjective
map $\CS(G/G\der) \to \CS(Z_G^0/Z_G^0\cap G\der)$.  This implies 
that all the tame local systems (resp. $F$-stable tame local 
systems) on $Z_G^0$ are obtained as the 
pull back from tame local systems on $G/G\der$ (resp. tame local
systems on $G/G\der$ such that
its restriction on $Z_G^0$ is $F$-stable) by the map $f$.
\par
This fact can be explained in the following way from a view point of 
Shintani descent. We consider a map 
$f: Z_G^{0F^m} \hra G^{F^m} \to G^{F^m}/G^{F^m}\der$ for sufficiently 
divisible integer $m$.  Let $\th'$ be a linear character of 
$G^{F^m}/G^{F^m}\der$.  Then $\th_1 = \th'\circ f$ is a linear 
character of $Z_G^{0F^m}$ whose restriction to 
$Z_G^{0F^m}\cap G^{F^m}\der$ is trivial, and all such characters of 
$Z_G^{0F^m}$ are obtained by the pull back by $f$ from a linear
character of $G^{F^m}/G\der^{F^m}$.  We now consider the norm map
$N_{F^m/F}: Z_G^{0F^m} \to Z_G^{0F}$, which is a surjective homomorphism 
with kernel 
$K = \{ z \in Z^{0F^m} \mid zF(z)\cdots F^{m-1}(z) = 1\}$.
Now the Shintani descent  
gives a bijection between the set of irreducible characters of
$Z_G^{0F}$ and the set of irreducible characters of $Z_G^{0F^m}$ whose
restriction on $K$ is trivial.  Since $Z_G^0 \cap G\der$ is an
$F$-stable finite set, we see that $Z_G^0\cap G\der \subset K$ if 
$m$ is chosen large enough.  In particular, any irreducible character
$\th_0$ of $Z_G^{0F}$ can be obtained as $Sh_{F^m/F}(\th_1)$ with 
$\th_1 = \th'\circ f$, where $\th'$ is a linear character 
of $G^{F^m}/G^{F^m}\der$ such that $\th_1$ is $F$-stable. 
\para{6.7.}  
Let $\wt L$ be a Levi subgroup of the standard parabolic subgroup
$\wt P$ of $\wt G$ containing $T$, and put $L = G \cap \wt L$.  
Assume that $L$ contains 
a cuspidal character sheaf. Then by 6.4, $\wt L$ is of the form
$\wt L = \wt L_1 \times \cdots \times \wt L_r$, where 
$\wt L_i \simeq GL_d \times\cdots\times GL_d$ ($n_i/d$-times) 
for a fixed $d$.   
Under the notation of 6.5 applied to $L$, let 
$A_0 = \CE\otimes A_{z,\e}$
be the cuspidal character sheaf on $L$, where $z \in Z_L/Z_L^0$ 
and $\e \in A_L(\dz u_1)\wg$ with $u_1 \in C$ ($C$ is the regular 
unipotent class in $L$), and $\CE' \in \CS(L/L\der)$.  
Let $K = \ind_P^G A_0$ be the induced 
complex of $G$. Then $K$ is a semisimple perverse sheaf whose 
simple components are character sheaves on $G$.
All the characters sheaves on $G$ are obtained in this way by 
decomposing a suitable $K$.
Let $\CE_1$ be the tame local system on $zZ_L^0$ obtained 
by restricting $\CE$ on $zZ_L^0$.  
Let $W_L$ be the Weyl subgroup of $W$ corresponding to $L$, and 
put 
\begin{equation*}
\CW_{\CE_1} = \{ w \in N_W(W_L) \mid w^*\CE_1 \simeq \CE_1 \}/W_L.
\end{equation*}
Let $\End K$ be the endomorphism algebra of $K$ in $\CM G$, and 
$\Ql[\CW_{\CE_1}]$ be the group algebra of $\CW_{\CE_1}$ over
$\Ql$.  We note that 
\begin{equation*}
\tag{6.7.1}  \End K \simeq \Ql[\CW_{\CE_1}].
\end{equation*}
In fact it is known by [L3, 10.2] that $\End K$ is isomorphic to 
the twisted group algebra of $\CW_{\CE_1}$.  
By a general principle, we have only to show that $K$ contains a 
character sheaf with multiplicity one. 
Let $A_1 = \CE\otimes A_{1, \e}$ be the cuspidal character sheaf 
in the case where $z = 1$, and put $K_1 = \ind_P^G A_1$.  We may 
choose a representative $\dz$ of $Z_L/Z_L^0$ so that $\dz \in Z_G$. 
Then we have $\dz^*K \simeq K_1$, and it is enough to consider the
case where $z = 1$, i.e., $K = K_1$.  But in this case, it is known
by [L4, 2.4] that $\End K \simeq \Ql[\CW_{\CE_1}]$.  This shows 
(6.7.1).
\par
In view of the discussion in 6.6, (6.7.1) is interpreted also in the
following form. Take $m$ large enough so that $\CE'$ is $F^m$-stable.
Let $\th'$ be the linear character of $(L/L\der)^{F^m}$ corresponding
to $\CE'$, and $\th$ the linear character of $L^{F^m}$ obtained 
by the pull back of $\th'$.  Let $\th_1$ be the restriction of $\th$
on $Z_L^{0F^m}$ and put 
\begin{equation*}
\CW_{\th_1} = \{ w \in N_W(W_L) \mid {}^w\th_1 = \th_1\}/W_L.
\end{equation*}
Then we have
\begin{equation*}
\tag{6.7.2}
\CW_{\CE_1} \simeq \CW_{\th_1}.
\end{equation*}
\para{6.8.} 
We keep the notation in 6.7.
We denote by $A_E$ the character sheaf 
occurring as the simple component in $K$ corresponding to 
$E \in \End K$.  For each $A_E \in \wh G^F$, we shall determine 
$\f_{A_E} : F^*A_E \isom A_E$.  Assume that 
$A_{z,\e}$ is 
$F$-stable, and put 
\begin{equation*}
\CZ_{\CE_1} = \{ w \in N_W(W_L) \mid (Fw)^*\CE_1 \simeq \CE_1\}/W_L.
\end{equation*}
By [L3, II, 10.2], we see that 
\begin{equation*}
\CZ_{\CE_1} \simeq  
    \{ w \in N_W(W_L) \mid (Fw)^*A_0 \simeq A_0\}/W_L.
\end{equation*}
\par
Now any $w \in \CZ_{\CE_1}$ can be written as $w = w_{\CE_1}y$ with
$y \in \CW_{\CE_1}$.  $Fw_{\CE_1}$ induces an automorphism 
$\g_{\CE_1} : \CW_{\CE_1} \to \CW_{\CE_1}$.
For each $w \in \CZ_{\CE_1}$, let $L^{w} \subset P^{w}$ be the 
$F$-stable Levi, and parabolic 
 subgroup of $G$ obtained from $L \subset P$ by twisting $w$, 
so that $(L^w)^F \simeq L^{F\dw}$, etc.  (Note that this definition 
of $L^w$ etc. is not the same as $L^w$ given in [L3, II, 10.6].  
$L^w$ defined there coincides with our $L^{w'}$ with $w' = F\iv(w\iv)$.
The formulas below are derived from [L3, II] under a suitable
modification.)   
We denote by $A_0^w, K^{w}$ 
the corresponding cuspidal character sheaf on $L^{w}$ and the
induced complex on $G$.  Then we have 
$K^{w} = \ind_{P^w}^GA_0^w$.
By applying the arguments in 6.5 to the $Fw$-stable subgroups 
$L^w \subset P^w$, one can construct 
$\vf_0^w = \vf^w_1\otimes \wt\vf^w_0: F^*A^w_0 \isom A^w_0$,
where $\vf_1^w, \wt\vf_0^w$ are constructed from 
$\vf_1: (Fw)^*\CE_1 \isom \CE_1$, $\wt\vf_0: (Fw)F^*A_y \isom A_y$ 
as in 6.5 (replacing $F$ by $Fw$).  Put 
$\f^w_0 = q^{b_0}\vf^w_0$, where 
\begin{equation*}
b_0 = (\codim_LC_0 - \dim Z_L^0)/2
\end{equation*}
with the regular unpotent class $C_0$ in $L$. 
Now $\vf^w_0: F^*A^w_0 \isom A^w_0$ induces an isomorphism 
$\vf^w: F^*K^w \isom K^w$.  
\par
We look at the special case where $w = w_{\CE_1}$.  Then $A_E$ 
occurring in $K^w$ is $F$-stable if and only if $E \in \CW_{\CE_1}$ is 
$\g_{\CE_1}$-stable.  We fix a preferred extension $\wt E$ of $E$
for each $E \in (\CW_{\CE_1}\wg)^{\w_{\CE_1}}$.
Then $\vf^w: F^*K^w \isom K^w$ induces an isomorphism 
$\vf_{A_E} : F^*A_E \isom A_E$ satisfying the following formula
(cf. [L3, 10.4, 10.6]).  For each $w = w_{\CE_1}y \in \CZ_{\CE_1}$, we have 
\begin{equation*}
\tag{6.8.1}
\x_{K^w, \vf^w} = \sum_{E \in (\CW_{\CE_1}\wg)^{\g_{\CE_1}}}
\Tr(\g_{\CE_1} y, \wt E)\x_{A_E,\vf_{A_E}},
\end{equation*}
\par
We define a normalized isomorphism $\f_{A_E}: F^*A_E \isom A_E$ by 
$\f_{A_E} = q^{b_0}\vf_{A_E}$.  In general, for an $F$-stable character sheaf 
$A = A_E$, we put  $\f_{A} = \f_{A_E}$, and we denote by  
$\x_A$ the characteristic function 
$\x_{A, \f_{A}}$ for $A \in \wh G^F$.
   
\par\bigskip
\section{Lusztig's conjecture}
\para{7.1.}
We give a description of cuspidal irreducible characters of 
$G^F$.  Assume that $\ol\CM_{s, E}$ contains a cuspidal character.
Then there exists $\ds_x \in \wt T$ such that $\pi(\ds_x) = s$,
$\ds_x$ is a regular semisimple element such that 
$Z_{\ds_x}$ consists of a Coxeter element in $W$.  
$W_{\ds_x} = \{ 1\}$,
and $\wt\r_{\ds_x, E}$ gives a cuspidal character of $\wt G^F$, where
$E$ is the trivial character of $W_{\ds_x} = W_s^0$.  
We assume further that the pair $(s, E) = (s,1)$ is of the form
as in (2.5.1). Then $W_s = \Om_s$, and 
$\Om_s$ is a cyclic group generated by $w_0$. $\wt G$ is of the
form    
$\wt G = \wt G_1 \times \cdots\times \wt G_r$, where 
$\wt G_i = GL_t$ with $t = n/r$ for an integer $t$ prime to $p$.  
Moreover, $w_0 = (w_c, \dots, w_c)$ with a
cyclic permutation $w_c = (1, 2, \dots, t) \in \FS_t$.
Then $\{ s\}$ is an $F$-stable regular semisimple class in $G^*$ 
such that $\Om_s \simeq \BZ/t\BZ$, which is unique modulo $Z_{G^*}$.
We make the following specific choice of $s$.
\par\medskip\noindent
(7.1.1) \ Let $\ds = (\ds_1, \dots, \ds_r) \in \wt T^* = 
      \wt T^*_1\times\cdots\times \wt T_r^*$ 
such that $\ds_i = \Diag(1, \z, \dots, \z^{t-1})$,
where $\z$ is a primitive $t$-th root of unity in $\ol\BF_q$.
Put $s = \pi(\ds)$.
\par\medskip
Let $s'$ be a regular semisimple element such that 
$\{ s'\}$ is $F$-stable and that $\Om_{s'} \simeq \BZ/t\BZ$.  
Since they are unique modulo 
$Z_{G^*}$, $s'$ can be written as $s' = zs$ with $z \in Z_{G^*}$.
If $z \in Z_{G^*}^F$, then 
$z$ determines a linear character $\th$ of $G^F$ under   
the natural isomorphism $Z^F_{G^*} \simeq (G^F/G^F\der)\wg$. 
The parameter set $\ol\CM_{s, E}$ and $\ol\CM_{sz, E}$ are then 
naturally 
identified, and we have a bijection $\CT_{s,E} \simeq \CT_{sz, E}$
via $\r_y \lra \th\otimes\r_y$ with $y \in \ol\CM_{s,E}$. 
\par
Assume that $\ds$ is as in (7.1.1).  
We assume further that $F$ acts trivially on $Z_G/Z_G^0$.  
Since $Z_G/Z_G^0 \simeq \BZ/t\BZ$, 
$F$ stabilizes $\z$.  Hence $F(\ds) = \ds$ and so $F' = F$.
Moreover, $F'$ acts trivially on 
$\Om_s$, and we have 
\begin{equation*}
\ol\CM_{s, E} = \Om_s\wg \times \Om_s. 
   \end{equation*}
By 2.8 (a), 
$\ol\CM_{s,E}$ is in bijection with $\ol\CM_{s, N}$.  Since $\ds$ is
regular semisimple, $\CO_{\wt\r}$ for $\wt\r = \wt\r_{\ds,E}$ is the regular 
nilpotent orbit.  So we assume that $N$ is regular nilpotent. 
In this case, $\wt P_N = \wt B$ and $\wt L_N = \wt T$, and 
$A_{\la} \simeq Z_G/Z_G^0$.
Since
$A_{\la}^F \simeq \BZ/t\BZ$, we have $\ol A_{\la} = A_{\la}$, and
so $\ol\CM_{s, N} = A_{\la}\times A_{\la}\wg$.  
Let $(A_{\la})\wg_0$ be the set of faithful characters of $A_{\la}$. 
We show
\begin{lem} 
Let $\ds$ and $s = \pi(\ds)$ be as in (7.1.1).  
Assume that $F$ acts trivially on $Z_G/Z_G^0$.
\begin{enumerate}
\item
For each $x \in \Om_s$, we have $\D(\wt\r_{\ds_x, E})|_{Z_G^{0F}} = 1$.
\item
Assume further that $q \equiv 1 \pmod 4$ if $t = 2$.  Then for each 
$(c, \xi) \in \ol\CM_{s, N}$,  
$\r_{c,\xi} \in \CT_{s, E}$ is cuspidal if and only if 
$\xi \in (A_{\la})\wg_0$.
\end{enumerate}
\end{lem}
\begin{proof}
For each $x \in (\Om_s)_{F'} = \Om_s$, we consider 
$\wt\r_{\ds_x, E} \in \Irr \wt G^F$.  
$\D(\wt\r_{\ds_x, E})$ is determined from $\ds_x$ as in 2.4 
by using an $F$-stable Levi subgroup $\wt M$ of an $F$-stable 
parabolic subgroup of $\wt G$. 
Assume that $\ds_x \in \wt T^{*xF}$  for $x \in \Om_s$. 
In the case where $N$ is regular 
nilpotent, $\wt M$ coincides with $\wt T$.  We now choose $m \ge 1$
such that $\ds_x \in \wt T^{*F^m}$ and that $(xF)^m = F^m$.  
Then by the duality of the torus, 
there exists an $xF$-stable linear character 
$\vT$ of $\wt T^{F^m}$ corresponding 
to $\ds_x$. The Shintani descent $Sh_{F^m/xF}(\vT) = \vT_0$  gives
rise to a linear character of $\wt T^{xF}$.
The restriction of $\vT_0$ on $Z_G^{xF} = Z_G^F$ gives the character 
$\D(\wt\r_{\ds_x, E})$.    
\par
Following the arguments in [S2, Corollary 2.21], we give a more
precise description of $\vT$ and $\vT_0$.  
Since $\pi(\ds_x) = \pi(\ds) = s$, there exists 
$z_x \in (\ker\pi)^{F^m} \subset Z_{\wt G^*}^{F^m}$
such that $\ds_x = \ds z_x$.
  Hence there exists an $F$-stable
linear character $\wh\vT$ of $\wt T^{F^m}$ and a linear character
$\w_x$ of $\wt G^{F^m}$ such that $\vT = \wh\vT\w_x$. 
Now there exists a decomposition 
$\wt T = \wt T^+_1 \times \cdots\times  \wt T^+_{t}$ 
such that $\wt T^+_i$ are
all $F$-stable, and $x$ permutes the factors $\wt T^+_i$. 
According to this decomposition of $\wt T$, $\w_x$ can be written as
\begin{equation*}
\w_x|_{\wt T^{F^m}} = \w^1_x\boxtimes\cdots\boxtimes \w^1_x,
\end{equation*}
where $\w^1_x$ is a linear character of 
$\wt T^{+F^m}_i \simeq \wt T^{+F^m}_1$.
Here $\w^1_x$ is expressed as $\w^1_x(y_1) = \ol\w_x(\det y_1)$ 
for $y_1 \in\wt T_1^{+F^m}$, where $\ol\w_x$ is a homomorphism
$\BF_{q^m}^* \to \Ql^*$.
Under the decomposition of $\wt T$ into $\wt T^+_i$, $Z_G^0$ can 
be identified with the set of $y = (y_1, y_1, \dots, y_1)$ 
($t$-times) such that $\det y_1 = 1$.
It follows that $\w_x|_{Z_G^{0F^m}} = 1$. 
On the other hand, 
our choice of $\ds$ implies that $\wh\vT$ is written on $\wt T^{F^m}$ as
\begin{equation*}
\wh\vT = 1\boxtimes a\boxtimes\cdots\boxtimes a^{t-1}
\end{equation*} 
with $a^t = 1$.
We may assume that the generator $w_0 \in \Om_s$ permutes the 
factors $\wt T^+_i$ by 
$w_0(\wt T^+_i) = \wt T^+_{i+1}$ for $i \in \BZ/t\BZ$.
Put $x = w_0$.
Since $\wh\vT$ is $F$-stable and $\wh\vT\w_x$ is $xF$-stable,
we see that $a = (\w^1_x)\iv F(\w^1_x) = (\w^1_x)^{q-1}$.
This implies that $\wh\vT|_{Z_G^{0F^m}} = 1$. 
Hence  $\vT|_{Z_G^{0F^m}} = 1$, and we see that 
$\vT_0|_{Z_G^{0F}} = 1$.  This proves the first assertion.
\par 
By the previous argument, the restriction of $\vT_0$ on 
$Z_G^F$ gives rise to a character $\xi_x \in A_{\la}\wg$.
We know that $\wt\r_{\ds_x, E}$ is cuspidal if and only if 
$x$ is a generator of $\Om_s$, and that 
all the cuspidal irreducible characters in $\CT_{s, E}$ are 
obtained as the irreducible components $\r_{c, \xi_x}$ 
in $\wt\r_{\ds_x, E}|_{G^F}$. 
Hence in order to prove the second assertion, we have only to
show that $\xi_x \in (A_{\la})_0\wg$ if $x$ is a generator
of $\Om_s$. 
First we show the following. 
\par\medskip\noindent
(7.2.1) \ Assume that $x = w_0$.  Then $\xi_x$ is a faithful
character.
\par\medskip 
It is easy to see that $\w_x|_{Z_{\wt G}^{F^m}}$ is $F$-stable.
Put $\w_0 = Sh_{F^m/F}(\w_x|_{Z_{\wt G}^{F^m}})$.
Let us choose $y \in Z^F_G$ such that the image of $y$ in 
$Z_G/Z_G^0 = A_{\la}$ is a generator of $A_{\la}$.
Then $y = (y_1, \dots, y_1) \in Z^F_G$ with $\det y_1$ a
primitive $t$-th root of unity. By a similar argument 
as in [S2, p.208 - p.209], we see that $\w_0(y)$ is a primitive
$t$-th root of unity.
On the other hand, put 
$\vT_1 = Sh_{F^m/F}(\wh\vT|_{Z_{\wt G}^{F^m}})$. Any element 
$z \in Z_{\wt G}^{F^m}$ can be written as $z = (z_1, \dots, z_1)$
with $z_1 \in \wt T^{+F^m}_1$ and 
\begin{equation*}
\wh\vT(z) = a^{t(t-1)/2}(z_1) = (\w_x^1)^{t(t-1)(q-1)/2}(z_1).
\end{equation*}  
Since $a^t = 1$, we have $(\w_x^1)^{t(q-1)} = 1$.
It follows that $\wh\vT^2|_{Z_{\wt G}^{F^m}} = 1$ and so 
$\vT_1^2 = 1$.  Since $\vT_0|_{Z_{\wt G}^F} = \vT_1\w_0$, 
we see that $\xi_x$ is  a faithful character if $t > 2$.  
So assume that $t = 2$.  Let 
$y = (y_1, y_1) \in Z_{\wt G}^F$
be such that $\det y_1 = -1$.  We have only to show that 
$\vT_1(y) = 1$.
Take $z = (z_1, z_1) \in Z_{\wt G}^{F^m}$ such that
$N_{F^m/F}(z) = y$.  Then it is easy to see that 
$N_{F^m/F}(\det z_1) = \det y_1 = -1$.
Then we have 
\begin{equation*}
\vT_1(y) = \wh\vT(z) = (\w_x^1)^{q-1}(z_1) = \ol\w_x^{q-1}(\det z_1).
\end{equation*} 
Since $(\w_x^1)^{2(q-1)} = 1$, we may assume that $(\w_x^1)^{q-1}$
is a character of order 2.  Hence $(\ol\w_x)^{q-1}$ is the unique 
character of order 2 of $\BF_{q^m}^*$, and one can write as 
$(\ol\w_x)^{q-1} = \th\circ N_{F^m/F}$, where $\th$ is the unique 
character of order 2 of $\Fq^*$, i.e., $\th(x) = x^{(q-1)/2}$ for 
$x \in \Fq^*$.  It follows that 
\begin{equation*} 
\ol\w_x^{q-1}(\det z_1) = \th(-1) = (-1)^{(q-1)/2} = 1
\end{equation*}
since $q \equiv 1 \pmod 4$ by our assumption.
Hence we have $\vT_1(y) = 1$, and $\xi_x$ is faithful in this case
also.  Thus (7.2.1) is proved.
\par
If we replace $x$ by $x^j$ for some $j$, then we have
$Sh_{F^m/F}(\w_{x^j}|_{Z_{\wt G}^{F^m}}) = \w_0^j$ and the 
previous argument shows that 
$Sh_{F^m/x^jF}(\vT|_{Z_{\wt G}^{F^m}}) = \vT_1\w_0^j$. 
Since $\w_0^j(y^i)$ are 
all distinct for $i = 1,\dots, t$ if $j$ is prime 
to $t$, we see that $\xi_{x^j}$ is faithful if $j$ is 
a generator of $\Om_s$.  This proves the second assertion, and the lemma 
follows.
\end{proof}
\para{7.3.}
We preserve the setting in 7.1. 
Removing the assumption on $F$, we consider the sets 
$\ol\CM_{s, N} = (A_{\la})_F \times (A_{\la}^F)\wg$ and 
$\CM_{s,N} = A_{\la}^F\times (A_{\la}\wg)^F$, which 
are in bijection with $\ol\CM_{s, E}$ and 
$\CM_{s,E}$.  Since $A_{\la} = \bar A_{\la}$, 
we have $\ol\CM_0 = \ol\CM_{s, N}$.  
Applying Lemma 7.2 to the situation in $G^{F^m}$, we have the 
following corollary.
\begin{cor}  
Assume that $m$ is sufficiently divisible so that 
$F^m$ satisfies the assumption in Lemma 7.2 with respect to
$F$.  Let $s'$ be a regular
semisimple element in $G^*$ such that $\CM_{s', E}$
contains an $F$-stable cuspidal irreducible character of $G^{F^m}$.
Then 
\begin{enumerate}
\item
The pair $(s', E)$ is of the from as in 2.8 (a), and 
$\wt G \simeq \wt G_1 \times \cdots\times \wt G_r$, where 
$\wt G_i \simeq GL_t$ with $t = n/r$, and 
$\{ s'\}$ is the unique regular semisimple class modulo $Z_{G^*}$ 
such that 
$W_{s'} = \Om_{s'} \simeq \BZ/t\BZ$.  In particular, 
$(\CT_{s',E}^{(m)})^F$ is parametrized by 
$\CM_{s', N}$, where $N$ is a regular nilpotent element in $\Fg^F$.
\item
Let $(A_{\la}\wg)_0^F$ be the set of $F$-stable, faithful
characters of $A_{\la} = Z_G/Z_G^0$. 
Then under the parametrization $(\CT_{s',E}^{(m)})^F \lra \CM_{s', N}$, 
$\r_{c,\xi}^{(m)}$ is cuspidal if and only if 
$\xi \in (A_{\la}\wg)^F_0$.
\end{enumerate}
\end{cor}
\begin{proof}
Since $s'$ is regular semisimple, the pair $(s',E)$ is of 
type (a) or (c) in 2.8.  But it is easy to see that in case (c),
$\CM_{s, E}$ does not contain a cuspidal irreducible character.
(Note that by the Shintani descent theory, Lusztig induction 
$R_L^G(\dw_1)$ corresponds to the Harish-Chandra induction from
$L^{F^m}$ to $G^{F^m}$).  Hence $(s',E)$ is of type (a), and (i) follows
from 7.1.  Now by 7.1, $s'$ can be written as $s' = zs$, where 
$s$ is as in (7.1.1) and $z \in Z_{G^*}$.
We may assume that $z \in Z_{G^*}^{F^m}$ by choosing $m$ large enough.  
Thus we have a natural bijection
$\CT^{(m)}_{s', E} \simeq \CT^{(m)}_{s,E}$ by 7.1 under the 
identification $\CM_{s',N}^{(m)} = \CM_{s,N}^{(m)}$.   
Since $F^m$ acts trivially on $A_{\la}$, (ii) follows from
Lemma 7.2.    
\end{proof}
We note the following lemma. 
\begin{lem} 
Let $s'$ be as in 7.1 and 
$\th = \vD(\wt\r_{\ds',E})|_{Z_G^F}$ the linear character of 
$Z_G^F$.  
Let $R_{z,\e}$ be the almost character
of $G^F$ for $(z,\e) \in \CM_{s', N}$ under the bijection
$\CM_{s', N} \simeq \CM_{s', E}$.  
Let $\vG_{c,\xi,\t'}$ be the modified generalized Gelfand-Graev 
character associated to $(c,\xi) \in \ol\CM_0 = \ol\CM_{s',N}$ and 
to a linear character $\th'$ of $Z_G^F$.  Then we have
\begin{equation*}
\lp\vG_{c,\xi,\th'}, R_{z,\e}\rp_{G^F}
          = \begin{cases}
               \e(c)\xi(z)|(Z_G/Z_G^0)^F|\iv 
                      &\quad\text{ if } \th'|_{Z_G^{0F}} = 
                                        \th|_{Z_G^{0F}} ,   \\
                0     &\quad\text{ if }  
                          \th'|_{Z_G^{0F}} \ne \th|_{Z_G^{0F}}.
             \end{cases}
\end{equation*}
\end{lem}
\begin{proof}
In our case, $\ol A_{\la} = A_{\la} = Z_G/Z_G^0$.  Thus by 
applying (4.5.2) and (4.5.3), one can write as 
\begin{equation*}
R_{z,\e} = |(Z_G/Z_G^0)^F|\iv\sum_{(c_1,\xi_1)}\e(c_1)\xi_1(z)\r_{c_1,\xi_1}.
\end{equation*}
Then by Theorem 2.6, (ii)-(b), we have 
\begin{equation*}
\lp\vG_{c,\xi,\th}, R_{z,\e}\rp_{G^F} 
               = |(Z_G/Z_G^0)^F|\iv\e(c)\xi(z) 
\end{equation*}
if $\th'|_{Z_G^{0F}} = \th|_{Z_G^{0F}}$.
This proves the first formula. 
The second formula also follows from Theorem 2.6, (ii)-(a) 
together with (4.5.2).
\end{proof}
\para{7.6.}
We preserve the notations  $\ol\CM_{s, N}$, etc. as 
in 7.1 for regular nilpotent element $N$.
Recall that $\CE\otimes A_{z,\e}$ is an $F$-stable cuspidal character sheaf on
$G$ as given in (6.5.3) for $z \in (Z_G/Z_G^0)^F$, 
$\e \in (A_G(\dz u_1)\wg_0)^F$ and $\CE|_{\dz Z_G^0}$: $F$-stable.  
Under the identification 
$A_G(\dz u_1) = A_G(u_1) \simeq A_{\la}$, we regard 
$\e$ as an element in $(A_{\la}\wg)^F$.  Hence 
$(z,\e)$ is regarded as an element in $\CM_{s, N}$.  
Let $\x_{\CE,y} = \x_{\CE, z,\e}$ be the characteristic function of 
$\CE\otimes A_{z,\e}$
defined in 6.5.  
Also recall the linear character $\th_0$ of $Z_G^{0F}$ associated to
$\CE'$ in 6.5.  
Let $\vG_{c,\xi,\th'}$ be as in Lemma 7.5.
We shall compute the 
inner product of $\vG_{c,\xi,\th'}$ with $\x_{\CE,z,\e}$. 
\begin{lem}  
Under the setting in 7.6, we have 
\begin{equation*}
\begin{split}
\lp \vG_{c,\xi,\th'}, &\x_{\CE, z,\e}\rp_{G^F} \\
    &=
       \begin{cases}
                \z_{\CI_0}\iv\xi(z)\e(cc_0)\iv\p(\hat z\iv)\th'(\dz)
                       |(Z_G/Z_G^0)^F|\iv
                 &\quad\text{ if } \th'|_{Z_G^{0F}}= \th_0, \\
                 0  &\quad\text{ if } 
                         \th'|_{Z_G^{0F}} \ne \th_0.
         \end{cases} 
\end{split}
\end{equation*}                 
where $\z_{\CI_0}$ is the fourth root of unity associated to the 
block $\CI_0 = \{ \io_0 \}$ for the cuspidal pair 
$\io_0 = (C,\CE_{\e})$.  $\dz \in Z_G^F$ is a representative of $z$, 
and $\hat z$ is as in 6.5. 
\end{lem}
\begin{proof}
We follow the notation in 2.1.  Since $N$ is regular nilpotent, $U$ is
the maximal unipotent subgroup of $G$ and $Z_L(\la_c) = Z_G$. 
$\xi\nat$ is a linear character of $Z_G^F$ which is trivial on $Z_G^{0F}$,
hence it is naturally identified with $\xi \in (A_{\la}^F)\wg$.
$\vG_{c,\xi,\th'}$ is given as 
$\vG_{c,\xi,\th'} = \Ind_{Z_G^FU^F}^{G^F}(\th'\xi\nat\otimes \vL_c)$.
It is easy to see that $\vG_{c, \xi,\th'}$ has the support in  
$Z_G^F G^F\uni$, and that
\begin{equation*}
\tag{7.7.1}
\vG_{c,\xi,\th'}(z'v) = \th'\xi(z')|Z_G^F|\iv\vG_c(v) 
             \qquad (z' \in Z_G^F, v \in G^F\uni).
\end{equation*}
Then by (6.5.4) and (6.5.7), we have
\begin{equation*}
\tag{7.7.2}
\lp\vG_{c,\xi,\th'}, \x_{\CE,z,\e}\rp_{G^F} \\ 
 = \begin{cases}
      A\lp \vG_c, \CX_{\io_0}\rp_{G^F}
          &\quad\text{ if } \th'|_{Z_G^{0F}} = \th_0, \\
      0   &\quad\text{ if } \th'|_{Z_G^{0F}} \ne \th_0
    \end{cases}
\end{equation*}
with 
\begin{equation*}
A = q^{(\codim C - \dim Z^0_G)/2}\p(\hat z\iv)\th'(\dz)\xi(z)
      |(Z_G/Z_G^0)^F|\iv.
\end{equation*}
Hence we have only to compute $\lp\vG_c, \CX_{\io_0}\rp_{G^F}$.
We compute it by making use of
Lusztig's formula (Theorem 5.6). 
By (5.2.4), $\CX_{\io_0}$ is orthogonal to any function in 
$\CV_{\CI_1}$ such that $\CI_1 \ne \CI_0$.  It follows that 
$\lp\vG_c, \CX_{\io_0}\rp_{G^F} = 
      \lp\ (\vG_c)_{\CI_0}, \CX_{\io_0}\rp_{G^F}$.
But since $\CI_0 = \{ \io_0\}$, 
we have $L = G, \CW = \{ 1\}$ and $\supp (\io_0) = C$.
Hence by Lemma 5.11, we have 
\begin{equation*}
\lp\ (\vG_c)_{\CI_0}, \CX_{\io_0}\rp_{G^F} = 
               q^{(-\codim C + \dim Z^0_G)/2}\z_{\CI_0}\iv\e(cc_0)\iv.
\end{equation*}
Substituting this into (7.7.2), we obtain the lemma. 
\end{proof}
\para{7.8.}
Returning to the original setting in 1.1, we consider 
an $Fw_{\d}$-stable cuspidal character $\d$ of $L^{F^m}$, and 
$\CP_{\d} = \Ind_{P^{F^m}}^{G^{F^m}}\d$ as in 4.9.
We shall describe $H(\d) = \End \CP_{\d}$ more precisely.
By Corollary 7.4, $(L, \d)$ is given as follows;  
$L = \wt L \cap G$, where 
$\wt L$ is a Levi subgroup of $\wt G$ such that 
$\wt L = \wt L_1 \times \cdots \times \wt L_r$ with 
$\wt L_i \simeq GL_t\times \cdots\times GL_t$ 
($n_i/t$-times) for a fixed integer $t$. 
There exists a (unique) cuspidal character $\wt\d$ of $\wt L^{F^m}$ 
belonging to $\CE(\wt L^{F^m}, \{ \ds \})$, where $\ds$ is as in 
(7.1.1) (by replacing $G$ by $L)$. 
Then the restriction of $\wt\d$ is a sum of cuspidal characters of 
$L^{F^m}$, and $\d$ is obtained in the form  
$\d = \th'\otimes\d_0$, where $\d_0$ is an irreducible consistent of 
$\wt\d|_{L^{F^m}}$ and $\th'$ is a linear character of 
$L^{F^m}$ corresponding to $z \in Z_{L^*}^{F^m}$. 
By Lehrer [Le, Theorem 10], 
$\CW_{\d} \simeq  \CW_{\d}^0\rtimes \Om_{\d}$, where
$\Om_{\d} \simeq \BZ/t_1\BZ$ for some integer $t_1 > 0$, and $\CW_{\d}$
acts on $L^{F^m}$ as a permutation of factors in the direct
product.  Moreover, $\CW_{\d}^0$ is isomorphic to the ramification 
group $\CW_{\wt\d}$ of $\wt\d$ in $\wt L^{F^m}$.
In our case, it is easy to see that $\CW_{\wt\d}$ is isomorphic to
$\CW$.  
\par
Now assume that $\d$ is a constituent of $\wt\d|_{L^{F^m}}$.  
Then $\d \in \CE(L^{F^m}, \{ s\})$.  Moreover we have 
$\CW_{\d} = \CW$ since $\CW_{\d}$ is a subgroup of $\CW$.
Let $\th$ be a linear character of $L^{F^m}$ as above. 
Then $\th\otimes\d$ is a cuspidal character belonging to 
$\CE(L^{F^m}, \{ zs\})$.  
\par
We can now prove Lusztig's conjecture for $G^F$. 
\begin{thm} 
Lusztig's conjecture holds for $G^F$.  More precisely, there 
exists a bijection between the set of $F$-stable 
character sheaves and the set
of almost characters of $G^F$ satisfying the following:
\begin{enumerate}
\item 
For each almost character $R_x$ of $G^F$, we denote by 
$A_x$ the corresponding character sheaf of $G$, and by $\f_x$ the 
isomorphism $F^*A_x \isom A_x$ as given in 6.8.  Then 
\begin{equation*}
\x_{A_x, \f_x} = \nu_xR_x
\end{equation*}
for a certain constant $\nu_x \in \Ql^*$.  Here $\nu_x$ is a root of 
unity contained in a fixed cyclotomic field independent of $q$. 
\item
Let $\x_{\CE, z, \e}$ be the characteristic function of 
the $F$-stable cuspidal character sheaf $\CE\otimes A_{z,\e}$ as given in 
(6.5.3).  Let $z_1 \in Z_{G^*}$ be the element corresponding to 
$\CE' \in \CS(G/G\der)$. 
Let $s$ be as in (7.1.1), and put $s' = z_1s$.  Let 
$R_{z,\e\iv}$ be the almost character
of $G^F$ corresponding to $(z, \e\iv) \in \CM_{s', N}$.
Then we have 
\begin{equation*}
\x_{\CE, z,\e} = \z_{\CI_0}\iv\e(c_0)\iv R_{z,\e\iv}. 
\end{equation*}
\end{enumerate}
\end{thm}
\begin{proof}
Let $L$ be the $F$-stable Levi subgroup containing $T$ 
of a proper standard parabolic subgroup $P$ of $G$.
By induction on $\dim G$, we may assume that Lusztig's conjecture
holds for any $L^{F'}$, where $F' = F\dw_1$ for some 
$\dw_1 \in N_G(L)$.  Assume that   
$A_0$ is an $F'$-stable  cuspidal character sheaf on $L$.
 Then $A_0$ can be
written as $A_0 = \CE\otimes A_{z,\e}$, where $z \in (Z_L/Z_L^0)^{F'}$, 
$u_1 \in C_0^{F'}$ ($C_0$ is the regular unipotent class in $L$),
 $\e \in A_L(\dz u_1)\wg_{F'}$ and 
$\CE$ is a local system on $L$ such that $\CE|_{\dz Z_L^0}$ is $F'$-stable.
By (ii), there exists an $F'$-stable cuspidal character 
$\d$ of $L^{F^m}$ associated to $(z, \e, \CE)$ such that 
the corresponding almost character of $L^{F'}$
is given by $R_{z,\e\iv}$.  
This holds for any $F' = Fw$ for $w \in N_W(W_L)$ such that 
$(Fw)^*A_0 \simeq A_0$.
\par
On the other hand, let $\CE_1$ be the tame local system on $Z_L^0$
obtained by restricting $\CE$ to $Z_L^0$,  and let $\CZ_{\CE_1}$ be as
in 6.8.
The above discussion then shows that $\CZ_{\CE_1} = \CZ_{\d}$.
Since $\CZ_{\d}$ is a coset of $\CW_{\d}$ 
and $\CZ_{\CE_1}$ is a coset of $\CW_{\CE_1}$, we see that 
\begin{equation*}
\tag{7.9.1}
\CW_{\d} = \CW_{\CE_1}.
\end{equation*}
\par
Now $w_{\CE_1}$ given in 6.8 coincides with $w_{\d}$ given in 
4.9, and so $\g_{\CE_1} = \g_{\d}$.  
Let $K^w, \vf^w$ etc. be as in 6.8.  Then by (6.8.1), we have
\begin{equation*}
\tag{7.9.2}
\x_{K^w, \vf^w} = q^{-b_0}\sum_{E \in (\CW_{\d}\wg)^{\g_{\d}}}
\Tr(\g_{\d} y, \wt E)\x_{A_E,\f_{A_E}},
\end{equation*}
where $A_E$ is the simple component of $K^{w}$ corresponding to
$E \in \CW\wg_{\CE_1}$.  On the other hand, 
under the isomorphism $(L^w)^F \simeq L^{F\dw}$, 
$\x_{A_0^w,\f_0^w} \in C((L^w)^F\ssim)$ is regarded as an element in 
$C(L^{F\dw}\ssim)$. Then it is known by 
[L6, Prop. 9.2] that 
\begin{equation*}
\tag{7.9.3}
\x_{K^{w}, \vf^w} = (-1)^{\dim C_0 + \dim Z^0_L}
          R_{L(w)}^G(q^{-b_0}\x_{A^w_0, \f^w_0}).
\end{equation*}
(Note that $\x_{\CE^{\sharp}}$ in [loc. cit.] coincides with 
$\x_{A_0^w, \vf_0^w} = q^{-b_0}\x_{A_0^w}, \f_0^w$. 
Also we note that (7.9.3) holds only under some restriction that
$q > q_0$ for some constant $q_0$ (see [loc. cit.]).  However, 
this restriction on $q$ can be removed by using a similar method
as in [S1] based on the Shintani descent identity of character
sheaves.  Since this argument will appear in the proof of 
Proposition 9.12 in a more extended form, we omit the details here.)  
\par
Since $\x_{A_0^w, \f^w_0}$ coincides with the almost character 
$R_{z,\e\iv}$ up to
scalar, we see by Proposition 4.10 together with (7.9.2) and (7.9.3),
that 
\par\medskip\noindent
(7.9.4) \ $\x_{A_{E}, \f_{A_{E}}}$ coincides with $R_{E}$ 
up to scalar.
\par\medskip
The above argument implies that $\x_{A, \f_A}$ is  
identified with some almost character of $G^F$ up to scalar unless
$A$ is cuspidal.  So we assume that $\wh G^F_0 \ne \emptyset$.  
Let $\CV_0$ be the subspace of $\CV_{G}$ spanned by 
$\x_{A, \f_A} \in \wh G^F_0$.  
Then in view of 
the previous discussion, $\CV_0$ coincides with the subspace of 
$\CV_{G}$ spanned
by almost characters of ${G}^F$ obtained from $F$-stable cuspidal
irreducible characters of ${G}^{F^m}$ by Shintani descent.
Now by Lemma 7.5, $R_{z,\e}$ is characterized as the
unique function in $\CV_0$ satisfying the property of inner product 
with various $\vG_{c,\xi,\th}$ for $(c,\xi) \in \ol\CM_{s',N}$ with 
$\th = \vD(\wt\r_{\ds',E})|_{Z_G^F}$. 
By Lemma 7.7, $\x_{\CE, z,\e}$ is also characterized by the inner product
with $\vG_{c,\xi,\th}$. Hence by comparing the formulas in Lemma 7.5
and Lemma 7.7, we see that 
\begin{equation*}
\x_{\CE,z,\e} = \z_{\CI_0}\iv\e(c_0)\iv \p(\hat z\iv)\th(\dz)R_{z,\e\iv}.
\end{equation*}
Here we claim that $\p(\hat z)= \th(\dz)$.  In fact, let 
$z_0$ be an element in $Z_{G^*}$ corresponding to 
$\CE' \in \CS(G/G\der)$, and $\dz_0 \in Z_{\wt G^*}$ a representative 
of $z_0$.  Let $\wt\p'$ be the linear character of 
$\wt G^{F^m}$ corresponding to $\dz_0$ for a large $m$.  
The restriction of $\wt\p'$ on $T^{F^m}$ gives the
character $\p$.  But it follows from the discussion on 
$\vD(\wt\r_{\ds',E})$ (see 2.4) that 
the restriction $\wt\p$ of $\wt\p'$ on 
$Z_{\wt G}^{F^m}$ is also $F$-stable, and the restriction of 
$Sh_{F^m/F}(\wt\p)$ on $Z_G^F$ gives the character $\th$.
This shows that $\p(\hat z) = \wt\p(\hat z) = \th(\dz)$, and 
the claim follows.  
Thus we have proved 
\begin{equation*}
\tag{7.9.5}
\x_{\CE,z,\e} = \z_{\CI_0}\iv\e(c_0)\iv R_{z,\e\iv},
\end{equation*}
and the assertion (ii) follows. 
\par
In order to complete the proof, we have only to show 
the assertion on the scalars $\nu_x$.
The assertion is certainly true for the case of cuspidal character
sheaves by (ii).  In the general case, this scalar
$\nu_x$ is given by Proposition 4.10 as 
$\nu_x = \nu_0\mu_{\wt\d,y}\iv\mu_{\wt E}$ under the notation of 
Proposition 4.10,  where $\x_{A_0^{w_{\d}y}} = \nu_0R_{\d,y}$, and 
$\x_{A_x} = \nu_xR_x$ (here $R_x = R_E$).  
Note that $\mu_{\wt\d,y}$ and $\mu_{\wt E}$ are
determined as in Theorem 4.7, by the choice a representative 
$\dc \in T^F$ of $c \in A^F_{\la}$ (for a fixed $F$) 
and of an extension field 
$\BF_{q^m}$; $m$ is chosen sufficiently divisible so that the 
Shintani descent gives the almost characters.  But the argument in 
[S2] shows that the requirement for $m$ only comes from the Shintani 
descent of modified generalized Gelfand-Graev characters, and so 
we can choose $m$ independent of the base field $\Fq$.  Hence 
$\mu_{\wt\d,y}$ and $\mu_{\wt E}$ are root of unities contained in a fixed
cyclotomic field (cf. Theorem 4.7).  This proves the assertion (i), 
and the theorem follows. 
\end{proof}
\par\medskip
\section{Parametrization of almost characters}
\para{8.1.}
Theorem 7.9 is based on the parametrization 
in terms of the induction from cuspidal character
sheaves, and its counter part for almost characters.
However, in order to decompose almost characters into
irreducible characters, one needs the parametrization
of almost characters given in 4.5. 
In this section, we discuss the relationship between 
these two parametrizations. 
In the remainder of this paper, we assume that  
$\wt G = GL_n$ and $G = SL_n$, for simplicity.  The general 
case is dealt with similarly.
\par
We consider the following semisimple element $\ds$ in $\wt G^*$,  
which is a more general type than (7.1.1).  
\par\medskip\noindent
(8.1.1) \ Let $t$ be a divisor of $n$ prime to $p$. 
Take $\ds \in \wt T^*$ 
such that 
\begin{equation*}
\ds = \Diag (\underbrace{1,\dots, 1}_{n/t\text{-times}}, 
               \underbrace{\z, \dots, \z}_{n/t\text{-times}},
                \dots, 
   \underbrace{\z^{t-1}, \dots, \z^{t-1}}_{n/t\text{-times}}),
\end{equation*}
where $\z$ is a primitive $t$-th root of unity in $k$. 
Put $s = \pi(\ds)$.  
\par\medskip
Then 
$W_{\ds} \simeq \FS_{n/t}\times\cdots\times \FS_{n/t}$
($t$-times), and 
$\Om_s$ is a cyclic group of order $t$ generated by 
$w_0 \in W$ which permutes the factors of $W_{\ds}$ cyclicly. 
Hence $W_s$ is of the form as given in (2.5.1).  
We note that the class $\{ s\}$ is the unique class in 
$G^*$ satisfying (2.5.1) for a fixed $t$. 
We now assume that $F^m$ acts trivially on 
$Z_G$.  Since $F^m(\z) = \z$, we have
$F^m(\ds) = \ds$. 
Then for $x = w_0^i \in \Om_s$,
$\ds_x \in \wt T^{*xF^m}$ is defined as $\ds_x = \ds z_x$ with 
$z_x \in Z_{\wt G^*}$ such that  
$z_x^{q^m-1} = \ds\iv x\iv\ds x = 
    \Diag(\z^{-i},\dots,\z^{-i}) \in Z_{\wt G^*}$. 
Let $d$ be a divisor of $t$, and consider the 
Levi subgroup $\wt L$ in  $\wt G$ such that 
$\wt L \simeq GL_{d} \times\cdots\times GL_{d}$ ($n/d$-times).
Let $\ds_L$ be a regular semisimple element in $\wt L^*$ defined 
as follows; 
under the isomorphism $\wt L^* \simeq GL_{d}\times\cdots\times GL_d$, 
$\ds_L$ is a product of $\Diag(1, \z_0, \dots, \z_0^{d-1})$, where
$\z_0 = \z^{t/d}$.  Put $s_L = \pi(\ds_L)$ under the natural 
map $\pi: \wt L^* \to L^*$.  Then one sees easily that there exists 
$\dz_L \in Z_{\wt L^*}$ such that $\ds_L\dz_L$ is $W$-conjugate to $\ds$.
Here $\Om_{s_L}$ is a cyclic group of order $d$ generated by 
$w_{0,L} \in W_L$, and there exists
an injective homomorphism $\Om_{s_L} \to \Om_s$ such that
the image of $w_{0,L}$ coincides with $w_0^{t/d}$.
For any $y \in \Om_{s_L}$, one can define
$(\ds_L)_y = \ds_Lz'_y$ for some $z'_y \in Z_{\wt G^*}$ as in the case of $G$.
It follows that $(\ds_L)_yz_L$ is $W$-conjugate to $\ds z'_y$.
If $y = (w_{0,L})^i$, then 
$(z_y')^{q^m-1} = \Diag(\z_0^i, \dots, \z_0^i)$, and we may 
choose $z_y' = z_{y^{t/d}}$.
Summing up the above argument, we have
\par\medskip\noindent
(8.1.2) \ For each generator $y$ of $\Om_{s_L}$, the class 
$\{ (\ds_L)_y\dz_L\}$ in $\wt G^*$ gives rise to a class $\{ \ds_x\}$ for 
some $x \in \Om_s$ such that the order of $x$ is $d$.  The
correspondence $\{ (\ds_L)_y\} \to \{\ds_x\}$ gives a bijection
between the classes in $\wt L^*$ associated to $y \in \Om_{s_L}$ of
order $d$,  and the classes in $\wt G^*$ associated to 
$x \in \Om_s$ of order $d$.
\para{8.2.}
We write $\ds_L\dz_L$ in a more explicit way.
$\dz_L\in Z_{\wt L^*}$ can be written as 
\begin{equation*}
\dz_L = (a, \dots, a) \quad\text{ with } 
    a = (1, \z, \dots, \z^{t/d-1}) \in (k^*)^{t/d}
\end{equation*}
under the isomorphism 
$Z_{\wt L^*} \simeq (k^*)^{t/d} \times\cdots\times (k^*)^{t/d}$
($n/t$-times).
Then we have 
\begin{equation*}
\ds_L\dz_L = (a_1, \dots, a_{t/d-1}, a_1, \dots, a_{t/d-1}, 
            \dots, a_1, \dots, a_{t/d-1})
\end{equation*}
with 
\begin{equation*}
a_j = \Diag(\z^{j-1}, \z^{t/d + j-1}, \z^{2t/d + j-1}, 
            \dots, \z^{(d-1)t/d + j-1}) \in GL_d.
\end{equation*}
From this, we see easily that there exists $w_L \in N_W(W_L)$
such that $Fw_L(\ds_L\dz_L) = \ds_L\dz_L$.
Then $\dz_L$ is also $Fw_L$-stable modulo $(Z_{\wt L^*})_d$, where
$(Z_{\wt L^*})_d = \{ z \in Z_{\wt L^*} \mid z^d = 1\}$.
We choose $w_L$ in a standard way. 
Hence $w_L$ gives the permutation
of factors in each block $a_1, \dots, a_{t/d-1}$.  We put $F'' = Fw_L$. 
Take $m$ large enough so that $\dz_L \in Z_{\wt L^*}^{F^m}$.
Let $\wt\th$ be the linear character of $\wt L^{F^m}$ corresponding to 
$\dz_L \in Z^{F^m}_{\wt L^*}$, and $\th$ the restriction of $\wt\th$ to
$L^{F^m}$.  Then by Lemma 1.2, $\wt\th|_{Z^{F^m}_{\wt L}}$ 
is $F''$-stable, and so $\th_1 = \th|_{Z^{0F^m}_L}$ is also $F''$-stable.
We put $\wt\th_1 = \wt\th|_{Z_{\wt L}^{F^m}}$.
\par
Let us denote by $z_L \in L^*$ the image of $\dz_L$ under the natural map
$\wt L^* \to L^*$. 
We consider the parameter set $\CM^L_{s_Lz_L, N_0}$ (a similar set as 
$\CM_{s,N}$ in 4.5, defined by replacing $G, F$ by $L, F''$) with respect 
to the $F''$-stable regular semisimple class $\{s_Lz_L\}$ in $L^*$, 
where $N_0$ is a regular nilpotent element in $\Lie L$.
Let $\d_{z,\e} (= \r_{z,\e})$ be a cuspidal irreducible character 
of $L^{F^m}$, stable by $F''$ corresponding to 
$(z,\e) \in \CM^L_{s_Lz_L, N_0}$.
Then there exists a cuspidal irreducible character 
$\wt\d = \wt\r_{(\ds_L)_y\dz_L,1}$ of $\wt L^{F^m}$ such that 
$\d_{z,\e}$ is an irreducible constituent of the 
restriction of $\wt\d$ to $L^{F^m}$.
Note that $\wt\d$ can be written as $\wt\d = \wt\th \otimes \wt\d'$
for $\wt\d' = \wt\r_{(\ds_L)_y,1}$.  The class $\{s_L\}$ is 
$F''$-stable, and there exists a cuspidal irreducible character
$\d'_{z,\e}$ parametrized by $(z,\e) \in \CM^L_{s_L, N_0}$,
 which is a constituent of $\wt\d'|_{L^{F^m}}$, 
such that $\d_{z,\e} = \th\otimes \d'_{z,\e}$.   
We consider the Harish-Chandra induction 
$\wt I = \Ind_{\wt P^{F^m}}^{\wt G^{F^m}}(\wt\d)$
and its restriction $I$ on $G^{F^m}$.
By (8.1.1), irreducible components of $I$ belong to 
$\CE(G^{F^m}, \{ s\})$. 
\par
By (6.7.2) and (7.9.1) we have $\CW_{\d} \simeq \CW_{\th_1}$. 
Since $\CW_{\wt\d} \simeq \CW_{\d}^0$, we see that 
$\CW_{\wt\d} \simeq \CW_{\th_1}^0$.  
We have $\CW_{\th_1} \simeq \CW_{\th_1}^0\rtimes \Om_{\th_1}$, 
where $\CW_{\th_1}^0$ is given as in 
\begin{equation*}
\tag{8.2.1}
\CW_{\th_1}^0 \simeq \FS_{n/t}\times\cdots\times \FS_{n/t} \quad
                  (t/d\text{-times}).
\end{equation*}
and $\Om_{\th_1}$ is the cyclic group of order $t/d$ 
generated by an element $y_0 \in \CW$ permuting 
the factors $\FS_{n/t}$ cyclicly.
Now $\wt I$ is decomposed as 
\begin{equation*}
\tag{8.2.2}
\wt I = \sum_{E_1 \in (\CW_{\th_1}^0)\wg}(\dim E_1)\wt\r_{E_1}, 
\end{equation*}
where $\wt\r_{E_1}$ is the irreducible characters of $\wt G^{F^m}$.
But since $\wt\r_{E_1}$ is contained in 
$\CE(\wt G^{F^m}, \{\ds_x \})$,  they are expressed as  
$\wt\r_{\ds_x, E_1'}$ with $E_1' \in W_{\ds}$.  The relationship 
of these two parametrizations are given as follows. 
We have $W_{\ds_x} \simeq (\FS_{n/t})^t$, and there exists a natural
embedding $\CW_{\th_1}^0 \to W_{\ds_x}$ via the diagonal embedding 
\begin{equation*}
(\FS_{n/t})^{t/d} \hra (\FS_{n/t})^t \simeq 
       (\FS_{n/t})^{t/d}\times\cdots\times (\FS_{n/t})^{t/d}.
\end{equation*}
Then one can define a map $f: (\CW_{\th_1}^0)\wg \to W_{\ds}\wg$ through
$E_1 \mapsto E_1\boxtimes\cdots\boxtimes E_1$.
We note that 
\par\medskip\noindent
(8.2.3) \ $\wt\r_{E_1} = \wt\r_{\ds_x, E_1'}$ with 
 $E_1' = f(E_1) \in W_{\ds_x}\wg$.
\par\medskip
In fact, let $w_1$ be the Coxeter element in $W_L$, and $\p_0$ a regular
character of $\wt T^{F}_{w_1} \simeq \wt T^{Fw_1}$ obtained by the Shintani descent 
$Sh_{F^m/Fw_1}(\th|_{\wt T^{F^m}})$.
The cuspidal character $\wt\d$ can be expressend as 
$\pm R^{\wt L}_{\wt T_{w_1}}(\p_0)$. 
Thus $\wt I$  coincides with 
$\pm R_{\wt L}^{\wt G}(R_{\wt T_{w_1}}^{\wt L}(\p_0))$, and 
\begin{equation*}
\tag{8.2.4}
R_{\wt L}^{\wt G}(R_{\wt T_{w_1}}^{\wt L}(\p_0)) = 
       R_{\wt T_{w_1}}^{\wt G}(\p_0)
    = \sum_{ E_1' \in W_{\ds}\wg}\Tr(w_1, E'_1)\wt\r_{E_1'}.
\end{equation*}
Then  we have $\Tr(w_1, E_1') = \dim E_1$ 
if $E_1' = E_1\boxtimes\cdots\boxtimes E_1$ for some 
$E_1 \in (\CW_{\th_1}^0)\wg$, and $\Tr(w_1, E_1') = 0$ otherwise.
Comparing (8.2.2) with (8.2.4), we see that 
$\wt\r_{E_1} = \wt\r_{E_1'}$ for $E_1 \in (\CW_{\th_1}^0)\wg$.
This proves (8.2.3).
\par
Put $\d = \d_{z,\e}$.  Then $\CW_{\d} \simeq \CW_{\th_1}$.
Take $E \in (\CW_{\d}\wg)^{F''}$, and put $E_1 = E|_{\CW_{\d}^0}$.
We assume that $E_1$ is of the form that 
\begin{equation*}
\tag{8.2.5}
E_1 = E_2\boxtimes\cdots\boxtimes E_2 \in (\CW_{\th_1}^0)\wg 
\end{equation*}
with $E_2 \in (\FS_{n/t})\wg$. 
Then $E_1' = f(E_1) \in W_{\ds}\wg$ satisfies the condition in 
(2.5.1).
Let $\r_E$ be the irreducible component of 
$\Ind_{P^{F^m}}^{G^{F^m}}\d$ corresponding to $E \in (\CW_{\d}\wg)^F$.
We denote by $R_E$ the almost character of $G^F$ corresponding to
$\r_E$ through the Shintani descent.
Let $\CO_{\wt\r}$ be the nilpotent orbit in $\Fg$ associated to 
$\wt\r = \wt\r_{\ds_x, E_1'}$.  
We choose a nilpotent element $N$ such that 
$\CO_{\wt\r} = \CO_N$.  Let $\wt P_N$ and $\wt L_N$ be the 
associated parabolic subgroup and its Levi subgroup in $\wt G$.  
We put $\wt M = \wt L_N$ and $M = \wt M \cap G$. 
We consider the modified generalized Gelfand-Graev characters 
$\vG_{c,\xi}$ of $G^F$ associated to $(c,\xi) \in \ol\CM_{s, N}$.
A similar formula as in Lemma 7.5 holds also for the general 
$\ol\CM_{s,N}$, which implies the following lemma.  
(Note that we don't need to consider $\vG_{c,\xi,\th}$ since 
$\vD(\wt\r)|_{Z_M(\la)^F}$ is trivial for any 
$\wt\r = \wt\r_{\ds_x,E_1'}$ 
in the case of $G = SL_n$ by [S2, Lemma 2.18].)
\begin{lem}  
The almost character $R_E$ is written as $R_{z',\e'}$ with 
$(z',\e') \in \CM_{s,N}$.  The correspondence 
$E \lra (z',\e')$ is determined by computing the inner
product $\lp \vG_{c,\xi}, R_E\rp$ for various $\vG_{c,\xi}$. 
\end{lem}
\para{8.4.}
Let $C_0$ be the regular unipotent class in $L$, and 
$\CE$ the local system of $L$ corresponding to the linear character 
$\th$ of $L^{F^m}$ in 8.2. 
Let $A_0 = A_{\CE,z,\e}$ be the $F''$-stable 
cuspidal character sheaf on $L$ associated to 
$(z,\e)\in \CM^L_{s_Lz_L,N_0}$ and $\CE$, defined as in 6.5, 
by replacing $G$ by $L$.  Put $K = \ind_P^G A_0$.
Then $\End K \simeq \Ql[\CW_{\CE_1}]$ with 
$\CW_{\CE_1} \simeq \CW_{\th_1}$.  We denote by 
$A_E$ the character sheaf of $G$ which is a direct summand  
of $K$ corresponding to $E \in \CW_{\th_1}\wg$.   
Put $\CZ_{\th_1} = \{ w \in \CW \mid {}^{Fw}\th_1 = \th_1\}$.
Then $\CZ_{\th_1} = w_L\CW_{\th_1}$, and $F'' = Fw_L$ acts
on $\CW_{\th_1}$.
For each $w \in \CZ_{\th_1}$, we denote by $\th^w_0$ the linear 
character of $Z_L^{0Fw}$ given by 
$\th^w_0 = Sh_{F^m/Fw}(\th_1)$. 
\par
Let $E \in \CW_{\th_1}\wg$ and $\CO_N$ be as in 8.2. 
By Theorem 7.9 (see (7.9.4)), 
$\x_{A_E}$ coincides with $R_{E}$ up to scalar.
Then we see that $N = N_{\mu} \in \Fg$ is of the form 
such that $\mu = (\mu_1 \ge \mu_2\ge \cdots \ge \mu_r)$ with
$\mu_j$ divisible by $t$.
We choose $N_{\mu}$ in Jordan's normal form corresponding to the 
partition $\mu$.  
Then $A_{\la} = Z_M(\la)/Z_M^0(\la)$ is a cyclic group of order 
$t'$, where $t'$ is the largest common divisor of 
$\mu_1, \dots, \mu_k$, coprime to $p$. Note that $Z_M(\la)$ is 
$F$-stable.  By our choice of $\ds z_L$ and of $w_L$, $w_L$ 
acts on $Z_M(\la)$, and so we have an action of $F'' = Fw_L$ on 
$Z_M(\la)$.  
\par
There exists an $F''$-stable subgroup $Z_M^1(\la)$ 
of $Z_M(\la)$ such that $\bar A_{\la} = Z_M(\la)/Z_M^1(\la)$ is  a
cyclic group of order $t$ (see [S2, 2.19]), which is given as follows; 
if we write the partition $\mu$ as $\mu = (1^{m_1}, 2^{m_2}, \dots)$.
then we have 
\begin{equation*}
Z_M^1(\la) \simeq 
  \{ (\dots, x_i, \dots) \in \prod_{m_i>0}GL_{m_i} 
          \mid \prod_i(\det x_i)^{i/t} = 1\}.
\end{equation*}  
(Note that $m_i = 0$ if $i$ is not divisible by $t$.)
We may choose $T \subset M \cap L$.  Then we have 
\begin{equation*}
\bar A_{\la} \simeq 
   Z_M(\la)/Z_M^1(\la) \simeq (Z_M(\la) \cap T)/(Z_M^1(\la) \cap T).
\end{equation*}
By our choice of $N$, we have $Z_M(\la) \cap T \subset Z_L$.
On the other hand, the above description of $Z_M^1(\la)$ implies 
that $Z_M^1(\la) \cap T \subset Z_L^0$.
It follows that we have a natural homomorphism 
\begin{equation*}
\tag{8.4.1}
\bar A_{\la} \to Z_L/Z_L^0
\end{equation*}
compatible with the action of $F$.
Note that $\e$ is an $F$-stable linear character of $Z_L/Z_L^0$.
Thus one can define an $F$-stable linear character 
$\e_1 \in (\bar A_{\la}\wg)^F$ as the pull back of $\e$ under the
above homomorphism..
\par
Since $T \subset M$, we have $Z_L \subset M$.  
One can write 
\begin{equation*}
\tag{8.4.2}
Z_L^0 = 
 \{ (a_1, \dots, a_{n/d}) \in (k^*)^{n/d} \mid \prod_ia_i = 1\}
\end{equation*}
under the identification 
\begin{equation*}
Z_L \simeq 
  \{ (a_1, \dots, a_{n/d}) \in (k^*)^{n/d} \mid \prod_i a_i^d = 1\}.
\end{equation*}
Let $zZ_L^0$ be a coset in $Z_L/Z_L^0$.  We may take a representative 
$z$ in $Z_G$.  
Then we have $Z_M(\la) \cap zZ_L^0 = z(Z_M(\la) \cap Z_L^0)$, and   
\begin{align*}
\tag{8.4.3}
Z_M(\la)\cap Z_L^0 &= 
  \{ (\underbrace{b_1, \dots, b_1}_{\mu_1/d\text{-times}},
      \underbrace{b_2,\dots, b_2}_{\mu_2/d\text{-times}}, \dots, 
  \underbrace{b_r,\dots, b_r}_{\mu_r/d\text{-times}}) \in Z_L^0 \mid 
                                b_i \in k^* \}  \\
           &\simeq \{ (b_1, \dots, b_r) \in (k^*)^r \mid 
                    \prod_i  b_i^{\mu_i/d} = 1 \}, 
 \end{align*}
and 
\begin{equation*}
Z^1_M(\la)\cap Z_L^0 \simeq \{ (b_1, \dots, b_r) \in (k^*)^r \mid 
                    \prod_i  b_i^{\mu_i/t} = 1 \}.  
\end{equation*}
Now $Z_M^1(\la) \cap Z_L^0$ acts on $z(Z_M(\la) \cap Z_L^0)$ by a 
left multiplication.  Put 
$X_M = (Z_M(\la)\cap zZ_L^0)/(Z_M^1(\la) \cap Z_L^0)$.  
It follows from the above argument that there exists a bijection 
\begin{equation*}
\tag{8.4.4}
  X_M \simeq \BZ/(t/d)\BZ.
\end{equation*}
If $zZ_L^0$ is $F''$-stable, then $F''$ acts naturally on $X_M$, 
which is compatible  
with the natural action of $F$ on $\BZ/(t/d)\BZ \simeq \lp \z^d\rp$.
\par
Assume that the coset $zZ_L^0$ is $Fw$-stable for $w \in \CW$.  Then 
we choose a representative $\dz_w \in Z_L^{Fw}$ of $zZ_L^0$.
Since $\dz_w \in T^{Fw}$, there exists $\a_w \in T$ such that 
$F^m(\a_w)\a_w\iv = \dz_w$.  We put $\hat z_w = \a_w\iv {}^{Fw}\a$.
Then we have $\hat z_w \in T^{F^m}$. 
Let $\p$ be a linear character of $T^{F^m}$ obtained by restricting 
$\th$ to $T^{F^m}$.  The value $\p(\hat z_w)$ does not depend on the
choice of $\a_w \in T$.
\par
We show the following lemma.
\begin{lem} 
Let the notations be as before.
\begin{enumerate}
\item
Let $t_1 \in Z_G \cap zZ_L^0$.  There exists a linear character 
$\Psi_{t_1}$ of $\CW_{\th_1}$ satisfying 
the following properties; $\Psi_{t_1}$ is  
trivial on $\CW^0_{\th_1}$, and is 
regarded as a character of $\Om_{\th_1}$.  
If $t_1 \in (zZ_L^0)^{F''}$, then $\Psi_{t_1}$ is $F''$-stable, and 
in that case, we have $t_1 \in (zZ_L^0)^{Fw}$ for $w = w_Ly$  
with $y \in \CW_{\th_1}$, 
 and 
\begin{equation*}
\p(\hat z_w)\th_0^{w}(t_1\dz_w\iv) = 
    \p(\hat z_{w_L})\th_0^{w_L}(t_1\dz_{w_L}\iv)\Psi_{t_1}(y).
\end{equation*}
\item
$\Psi_{t_1}$ depends only on the coset in $X_M$ to which $t_1$ belongs,
and we have 
\begin{equation*}
\{ \Psi_{t_1} \mid t_1 \in X_M^{F''} \} = (\Om\wg_{\th_1})^{F''}. 
\end{equation*}
\end{enumerate}
\end{lem}
\begin{proof}
Put $w = w_Ly$ for $y \in \CW_{\th_1}$.  
For $t_1 \in  Z_G \cap zZ_L^0$, choose $\a \in Z_{\wt L}$ such that
$t_1 = F^m(\a)\a\iv$.  
We choose $m' $,  a multiple of $m$, such that $\a$ is
$F^{m'}$-stable.  One can find a linear character 
$\wt\th'_1$ on $Z_{\wt L}^{F^{m'}}$, which is an extension of 
$\wt\th_1$ on $Z_{\wt L}^{F^m}$ (see 8.2), stable by $F''$, 
such that $\CW^0_{\th_1} \subseteq \CW_{\wt\th'_1}$.
We put $\Psi_{t_1}(y) = \wt\th_1'(\a\iv{}^y\a)$ 
(for a fixed $\a$ and $\wt\th_1'$). 
We have $\wt\th'_1(\a\iv{}^y\a) = \wt\th'_1(\a\iv)\wt\th_1'({}^y\a) = 1$
for $y \in \CW_{\th_1}^0$, and so $\Psi_{t_1}$ is trivial on $\CW^0_{\th_1}$. 
It follows that $\Psi_{t_1}(y) = \Psi_{t_1}(y_2)$ if $y = y_1y_2$ with 
$y_1 \in \CW_{\th_1}^0, y_2 \in \Om_{\th_1}$.  Hence 
in order to show that $\Psi_{t_1}$ is a homomorphism 
$\CW_{\th_1} \to \Ql^*$, it is enough to see that 
(*) $\Psi_{t_1}$ is a homomorphism on $\Om_{\th_1}$.  
This will be shown soon later.
\par
Now assume that $t_1 \in (zZ_L^0)^{F''}$.  
Since $z \in Z_G$, we have $t_1 \in (zZ_L^0)^{Fw} = \dz_wZ_L^{0Fw}$.
Choose $\b_w \in Z_L^0$ such that 
$F^m(\b_w)\b_w\iv = t_1\dz_w\iv \in Z_L^{0Fw}$.
Then $\th_0^w(t_1\dz_w\iv) = \th_1(\b_w\iv {}^{Fw}\b_w)$.
Since $t_1 = F^m(\a)\a\iv$, we may assume that 
$\a_w = \a\b_w\iv$.
In particular, we have $\a\iv {}^{Fw}\a \in T^{F^m}$ and 
$\p(\hat z_w)\th_0^w(t_1\dz_w\iv) = \p(\a\iv {}^{Fw}\a)$.
A similar formula also holds for $w_L$.
Since 
$\a\iv{}^{Fw}\a = \a\iv{}^{Fw_L}\a\cdot{}^{Fw_L}(\a\iv{}^y\a)$, 
we see that $\a\iv {}^y\a \in T^{F^m}$.
Since $\wt\th_1$ and $\p$ coincides with 
each other on $T^{F^m} \cap Z_{\wt L}^{F^m}$, 
we have 
\begin{align*}
\p(\hat z_w)\th_0^w(t_1\dz_w\iv) &= \p(\a\iv{}^{Fw}\a) \\ 
     &= \p(\a\iv{}^{Fw_L}\a)\p({}^{Fw_L}(\a\iv{}^{y}\a))  \\
     &= \p(\hat z_{w_L})\th_0^{w_L}(t_1\dz_{w_L}\iv)\Psi_{t_1}(y). 
    \end{align*}
Hence the formula in (i) holds.  In particular, $\Psi_{t_1}(y)$
does not depend on the choice of $\a \in Z_{\wt L}$.  
Put $\a' = F''(\a)$.  We have  
$t_1 = F^m(\a)\a\iv = F^m(\a')(\a')\iv$ 
since $t_1$ is $F''$-stable, and one can replace $\a$ by $\a'$ in 
defining $\Psi_{t_1}$.
Since $\wt\th_1'$ is $F''$-stable, we have
\begin{align*}
\Psi_{t_1}(y) = \wt\th'_1({\a}\iv{}^y\a) = 
    \wt\th'_1({\a'}\iv{}^{F''(y)}\a') = \Psi_{t_1}(F''(y)).
\end{align*}
This shows that $\Psi_{t_1}$ is $F''$-stable.
Thus the assertion in (i) was proved modulo (*).
\par
Next we show (ii).    
Take $t_1 \in Z_G \cap (Z_M(\la) \cap zZ_L^0)$.  Then as in (8.4.2)  
$t_1$ can be written as $t_1 = zb$ with $b \in Z_M(\la) \cap Z_L^0$. 
Hence $b = (b_1, \dots, b_r)$ such that 
$\prod_ib_i^{\mu_i/d} = 1$ with $b_i \in k^*$.
One can choose $\a \in T$ such that $F^m(\a)\a\iv = t_1$ 
as follows; $\a = \b\g$ with $\b \in Z_L^0$ and $\g \in Z_{\wt G}$
such that $F^m(\g)\g\iv = z$, $F^m(\b)\b\iv = b$.  More precisely, we 
can choose $\b$ as  
\begin{equation*}
\tag{8.5.1}
\b = (\nu(t_1)\iv\b_1, 
       \underbrace{\b_1, \dots, \b_1}_{\mu_1/d-1\text{-times}},
       \underbrace{\b_2, \dots, \b_2}_{\mu_2/d\text{-times}},
       \dots,
       \underbrace{\b_r, \dots, \b_r}_{\mu_r/d\text{-times}}),
\end{equation*}
where $\b_i^{q^m-1} = b_i$ for $i = 1, \dots, r$ and 
$\nu(t_1) = \prod_{i=1}^r\b_i^{\mu_i/d} \in \BF^*_{q^m}$.
Let $y_0$ be the generator of $\Om_{\th}$ as in 8.2.  
Since $y_0$ is a cyclic permutation of order $t/d$ of 
consecutive factors,  
it makes no change except the part 
$\nu(t_1)\iv\b_1, \b_1, \dots \b_1$.  Hence we see that 
\begin{equation*}
\tag{8.5.2}
\b\iv {}^{y^j_0}\b = 
     (\nu(t_1), 1, \dots, 1, \nu(t_1)\iv, 1, \dots, 1) \in Z_L^{0F^m}
\end{equation*}
for $1 \le j \le t/d-1$, where $\nu(t_1)\iv$ occurs in the $(j+1)$th
factors.  Moreover, since $\g \in Z_{\wt G}$, we have ${}^{y_0}\g = \g$. 
It follows that 
\begin{equation*}
\a\iv{}^{y_0^j}\a = \b\iv {}^{y_0^j}\b.
\end{equation*}
On the other hand, $\th$ is the restriction of the linear character 
$\wt\th$ of $(\wt L/\wt L\der)^{F^m}$ such that 
$\wt\th = 
  (1\boxtimes\wt\vT\boxtimes\cdots\boxtimes
           \wt\vT^{t/d-1})^{\boxtimes n/t}$, where $\wt\vT$ is a linear 
character of $GL_d^{F^m}$ of order $t$.  
Since $\wt L/\wt L\der \simeq Z_{\wt L}/(Z_{\wt L} \cap \wt L\der)$,
and $\wt L\der = L\der = SL_d \times \cdots \times SL_d$, 
the linear character $\wt\th_1$ on $Z_{\wt L}^{F^m}$ can be written as
\begin{equation*}
\tag{8.5.3}
\wt\th_1 = 
   (1\boxtimes \vT\boxtimes\cdots\boxtimes 
          \vT^{t/d-1})^{\boxtimes n/t},
\end{equation*} 
where $\vT$ is a linear character of 
$Z_{GL_d}^{F^m}$ which is trivial on $Z_{SL_d}^{F^m}$.  
Hence $\vT$ is identified with a linear character of 
$\BF_{q^m}^*$ of order $t/d$. 
If we replace $m$ by $m'$, $\vT$ can be extended to a linear character
$\vT'$ of $\BF^*_{q^{m'}}$, and 
$(1\boxtimes \vT'\boxtimes\cdots\boxtimes 
          {\vT'}^{t/d-1})^{\boxtimes n/t}$ gives rise to 
a linear character of 
$Z_{\wt L}^{F^{m'}}$, which gives $\wt\th_1'$.
It follows that 
\begin{equation*}
\tag{8.5.4}
\Psi_{t_1}(y^j_0) = \wt\th_1'(\a\iv{}^{y^j_0}\a) = \vT(\nu(t_1))^{-j}.
\end{equation*}
This proves (*) since if $t_1 \in Z_G \cap zZ_L^0$, then 
$t_1 \in Z_G \cap (Z_M(\la) \cap zZ_L^0)$.
\par
Now assume that $t_1 \in Z_M^1(\la)\cap Z_L^0$.
Since $\prod_ib_i^{\mu_i/t} = 1$, we have 
$\nu(t_1)^{(q^m-1)d/t} = 1$. It follows that there exists 
$\nu_1 \in \BF_{q^m}$ such that $\nu(t_1) = \nu_1^{t/d}$, 
and we have 
$\Psi_{t_1}(y^j_0) = 1$. 
This implies that $\Psi_{t_1}$ depends only on 
$t_1 \in X_M$.  
Now $X_M$ is in bijection with a cyclic group 
of order $t/d$, and we can choose a representative $x_0$
of a generator of $X_M$ as $x_0 = zc_0$ with   
$c_0 = (c, \dots, c) \in Z_G \cap Z_L^0$ such that $c^{n/t}$ 
is a primitive $t/d$-th root of unity in $k$.  
Put $\nu_0 = \nu(x_0)$.  Then $\nu_0$ is a generator of the 
cyclic group $\BF_{q^m}^*$, and we have $\nu(x_i) = \nu_0^i$ for
$x_i = zc_0^i$.
It follows from (8.5.1) that we have 
$\Psi_{x_i}(y_0) = \vT(\nu_0)^{-i}$ for $i = 0, \dots, t/d-1$.
Since $\nu_0$ is a generator of $\BF_{q^m}^*$, $\vT(\nu_0)$ is of 
order $t/d$ in $\Ql^*$.  
Hence $\Psi_{x}(y_0)$ are all distinct for $x \in X_M$.  Since 
$\Om_{\th_1}$ is a cyclic group of order $t/d$, we see that
there exists a bijection between $\{ \Psi_{x} \mid x \in X_M\}$ and
$\Om_{\th_1}\wg$.  
We note that $\Psi_x$ is $F''$-stable for $x \in X_M^{F''}$.
In fact,  take a representative $t_1 \in Z_M(\la) \cap zZ_L^0$ of 
$x \in X_M^{F''}$.  Then we have $F''(t_1) = t_1t_2$ with 
$t_2 \in Z_M^1(\la) \cap zZ_L^0 $.  
It follows that $\Psi_{F''(t_1)} = \Psi_{t_1t_2} = \Psi_{t_1}$.
Now take $\a \in Z_{\wt L}$ such that $F^m(\a)\a\iv = t_1$, and put
$\a' = F''(\a)$.  
Since $\wt\th_1'$ is $F''$-stable, we have
\begin{align*}
\Psi_{t_1}(y) = \wt\th_1'(\a\iv{}^y\a)
              = \wt\th_1'({\a'}\iv{}^{F''(y)}\a').
\end{align*}
But since $F''(t_1) = F^m(\a'){\a'}\iv$, we have
\begin{equation*}
 \wt\th_1'({\a'}\iv{}^{F''(y)}\a') = 
         \Psi_{F''(t_1)}(F''(y)) = \Psi_{t_1t_2}(F''(y)) 
                                 = \Psi_{t_1}(F''(y)). 
\end{equation*}
This implies that $\Psi_{t_1}(y) = \Psi_{t_1}(F''(y))$, and so 
$\Psi_x$ is $F''$-invariant. 
Actually this argument shows that $\Psi_x$ is $F''$-invariant
if and only if $x \in X_M^{F''}$.  Thus we have
$\{ \Psi_x \mid x \in X_M^{F''}\} = (\Om_{\th_1}\wg)^{F''}$,
and (ii) is proved. 
\end{proof} 
\begin{thm}  
Let $G, s$ and $L, s_L, z_L$ be as in 8.1, 8.2. 
In particular, $L$ is a Levi subgroup of $G$ such that
$L/Z_L$ is a product of $PGL_d$ with $d|t$.  Let 
$\th$ be the linear character of $L^{F^m}$ corresponding to
$z_L \in Z^{F^m}_{L^*}$, and let $\th_1 = \th|_{Z_L^{0F^m}}$ be 
the $F''$-stable linear character.  Let $\CE$ be the local system 
on $L$ corresponding to $\th$.   
\begin{enumerate}
\item
Let $A_{\CE, z,\e} = \CE\otimes A_{z,\e}$ be an $F''$-stable 
 cuspidal character sheaf of $L$ as in 7.6 for 
$(z, \e) \in  \CM^L_{z_Ls_L, N_0}$, and $A_E$ be the character sheaf 
corresponding to $E \in (\CW_{\th_1}\wg)^{F''}$, where 
$E_1 = E|_{\CW_{\th_1}^0}$ satisfies the condition in (8.2.5). 
Then there exists $z_E \in \bar A_{\la}^F$
satisfying the following.
\begin{equation*}
\tag{8.6.1}
\x_{A_{E}} = \nu_E R_{z_E, \e_1\iv},  
\end{equation*}
where $(z_E, \e_1\iv) \in \CM_{s, N}$ is given as follows; 
$ \e_1$ is the $F$-stable character of $\bar A_{\la}$ defined as 
the pull back of $\e$ in 8.4.
$z_E$ is determined uniquely by the following condition;
let $\io$ be the unique element of $\CI_0$ such that 
$\supp(\io) = \CO_N$, and let $E_{\io} \in \CW\wg$ the 
corresponding character under the 
generalized Springer correspondence.
Then there exists a unique class $x_E \in X_M^{F''}$ such that
\begin{equation*}
\lp E\otimes\Psi_{t\iv}, E_{\io}\rp_{\CW_{\th_1}} = 
       \begin{cases}
            1 &\quad\text{ if } t \in x_E, \\
            0 &\quad\text{ otherwise.}
       \end{cases} 
\end{equation*}
$z_E$ is the image of $x_E$ under the natural map 
$X_M^{F''} \to \bar A_{\la}^F$.
$\nu_E$ is a root of unity, which is determined explicitly
(see (9.10.6)). 
\par
All the almost characters $R_{z',\e'}$ associated to 
$(z' \e') \in \CM_{s,N}$  
such that $\e'$ is of order $d$ is obtained in this way from 
some $A_E$.
\item
Let $\d = \d_{z,\e} = \th\otimes\d'_{z,\e} \in \CT^{(m)}_{s_Lz_L, N_0}$ 
be the cuspidal irreducible character of $L^{F^m}$ as in 8.2.  
Let $\r_{E}$
be the $F$-stable irreducible character of $G^{F^m}$ contained in 
$\Ind_{P^{F^m}}^{G^{F^m}}\d$ corresponding to 
$E \in \CW_{\d}\wg$.  Then we have $\r_{E} = \r_{z_{E}, \e_1\iv}$
with $(z_{E}, \e_1\iv) \in \CM_{s, N}$.  
In particular, we have $R_{E} = R_{z_E, \e_1\iv}$. 
\end{enumerate}
\end{thm}
\remark{8.7.}
Theorem 8.6 gives an identification of two parametrizations
of almost characters of $G^F$, one given by the parameter set
$\CM_{s,N}$ and the other by the Harish-Chandra induction, 
 in the case where $(s, E)$ satisfies the property (2.5.1). 
The parametrizations of irreducible characters, and of almost 
characters, are described in 2.8, (a) - (c).  The above result
covers the case (a). The case (b) is obtained by applying 
the Harish-Chandra induction for the case (a).  Thus the above
parametrization can be extended also for this case.  The case
(c) is obtained by applying the twisted induction for the 
cases (a) or (b).  Since the twisted induction of almost 
characters corresponds to the induction of character sheaves, 
our result gives an enough information for decomposing the 
characteristic functions of character sheaves in terms of 
irreducible characters of $G^F$ based on our parametrization
$\ol\CM_{s,E}$.  
\para{8.8.}
The proof of the theorem will be done in the next section.
Here we prove some preliminary results. 
As discussed in 5.5, we identify 
the set $\CI_G$ with the set of pairs $(\CO, \CE)$, where 
$\CO$ is a nilpotent orbit in $\Fg$ and $\CE$ is a simple
$G$-equivariant local system on $\CO$.  For $\io \in \CI_G$ 
belonging to $(L, C_0, \CE_0) \in \CM_G$, we 
define two integers $b(\io)$ and $b_0$ as follows (cf. 5.1 and 6.8).
\begin{align*}
\tag{8.8.1}
b(\io) &= (a_0 + r)/2 = (\codim_G\supp(\io) - \codim_L C_0)/2, \\
b_0    &= (\dim G - \dim\supp K_{\io})/2 = (\codim_L C_0 - \dim Z_L)/2.
\end{align*}
\par
Suppose that $\io$ belongs to the triple $(L, C_0, \CE_0) \in \CM^{F}_G$. 
Let $z \in (Z_L/Z_L^0)^{F}$, and we choose a representative 
$\dz \in Z_L^{F}$. By the translation $C_0 \isom \dz C_0$, we may
regard $\CE_0$ the $F$-stable local system on $\dz C_0$.
Let $\vS = zZ_L^0\times C_0 = Z_L^0 \times \dz C_0$, and we follow
the setting in 6.7 and 6.8.  In particular, $A_0 = \CE\otimes A_{z,\e}$ is the 
cuspidal character on $L$ associated to $\CE_0 = \CE_{\e}$, and 
to the local system $\CE$ on $L$ which is the pull back of 
$\CE' \in \CS(L/L\der)$.  
We put $\CE_1 = \CE|_{zZ_L^0}$ as before,
and consider the corresponding character $\th_1$ of $Z_L^{0F^m}$, etc.
We identify $\CW_{\CE_1}$ with $\CW_{\th_1}$ and with $\CW_{\d}$, 
and similarly for $\CZ_{\CE_1}$ with $\CZ_{\th_1}$, $\CZ_{\d}$.
We write the automorphism $\g_{\CE_1}$ on $\CW_{\CE_1}$ as $\g_{\d}$.  
Let $K = \Ind_P^GA_0$ be the induced complex on $G$, and we consider 
$\vf^w: (Fw)^*K^w \isom K^w$ for $w \in \CZ_{\CE_1}$.  
Then by (7.9.2), we have
\begin{equation*}
\tag{8.8.2}
\x_{A_{E}} =   q^{b_0}|\CW_{\th_1}|\iv\sum_{y \in \CW_{\d}}
     \Tr((\g_{\d}y)\iv, \wt E)\x_{K^w,\vf^w}
\end{equation*}
with $w = w_{\d}y$. 
\par
Let $\CE_0' = \Ql\boxtimes \CE_0$ be the local system on
 $\vS = Z_L^0\times C_0$.
Let $\CF = \CE\otimes \CE_0'$ be the local system on $\vS$.  Then 
$A_0$ coincides with $\IC(\ol\vS, \CF)[\dim \vS]$.  
For each $w \in \CZ_{\CE_1}$, take $\a \in G$ such that 
$\a\iv F(\a) = \dw$, where $\dw$ is a representative of $w$ in
 $N_G(L)$.
Then $A_0^w$ is constricted from the twisted data 
$(L^w, \vS^w, \CF^w)$, 
where $L^w = \a L\a\iv, \vS^w = \a \vS\a\iv$ and 
$\CF^w = \ad(\a\iv)^*\CF$.  
By applying the argument in 6.5 (see also 6.8), 
we can construct an isomorphism 
$F^*\CF^w \isom \CF^w$ which induces $\vf_0^w: F^*A_0^w \isom A_0^w$. 
We denote this map also by $\vf_0^w$. 
We put $C_0^w = \a C_0\a\iv$.
The set $zZ_L^0$ is $Fw$-stable for $w \in \CZ_{\CE_1}$ and one can 
choose a representative $\dz_w \in Z_L^{Fw}$ of the coset $zZ_L^0$.
We have an $F$-stable set $\vS^w = \dz^wZ^0_{L^w}\times C_0^w$ 
with $\dz^w = \a \dz_w \a\iv \in Z^F_{L^w}$.
For each $w \in \CZ_{\th_1}$, let $\th_0^w = Sh_{F^m/Fw}(\th_1)$
be the linear character of $Z_L^{0Fw}$.  Under the isomorphism 
$Z_L^{0Fw} \simeq Z_{L^w}^{0F}$, we regard $\th_0^w$ as the character
of $Z_{L^w}^{0F}$, which we denote by $\bar \th_0^w$.
\par
Now take $t \in G^F$ and fix it.  Let $x \in G^F$ be an element such that 
$x\iv tx \in \dz^w Z^0_{L^w}$.  Then $L^w_x = xL^wx\iv$
is a Levi subgroup of some parabolic subgroup of $Z_G(t)$.
Let $C_x^w = xC_0^wx\iv$ be the unipotent class in $L_x^w$.
We denote by $\CF^w_x$ the local system on $C_x^w$ obtained by 
the pull back of $\CE^w_0 = \ad(\a\iv)^*\CE_0$ 
by $\ad x\iv: C_x^w \to C_0^w$. 
By the map
$\b: C^w_x \to \vS^w$, $v \mapsto x\iv tvx$,  
$\vf_0^w: F^*\CF^w \isom \CF^w$ induces an isomorphism
$\vf'_x = \b^*\vf_0^w: F^*\CF^w_x \isom \CF^w_x$.
Let $\vf_x: F^*\CF^w_x \isom \CF^w_x$ be the 
isomorphism defined by $\vf_x = (\ad x\iv)^*\vf^w_0$.  Then we have 
\begin{equation*}
\vf'_x = \bar\th_0^w(x\iv tx(\dz^{w})\iv)\p(\hat z_w)\vf_x
\end{equation*}
by (6.5.6), where $\p$ is a character of $T^{F^m}$ as given in
(6.5.6),
and $\hat z_w$ is as in 8.4.
Assume that $t, v \in G^F$ such that $tv = vt$, where $t$
is semisimple and $v$ is unipotent.  Then by 
the character formula [L3, II, Theorem 8.5] 
for the function $\x_{K^w,\vf^w}$, we have 
\begin{equation*}
\tag{8.8.3}
\x_{K^w, \vf^w}(tv) = |Z_G(t)^F|\iv
    \sum_{\substack{ x \in G^F \\ x\iv tx \in \dz^wZ_{L^w}^0}}
       Q^{Z_G(t)}_{L_x^w, C^w_x, \CF^w_x, \vf_x}(v)
               \bar\th^w_0(x\iv tx(\dz^w)\iv)\p(\hat z_w),
\end{equation*}
where $Q^{Z_G(t)}_{L_x^w, C_x^w, \CF_x^w, \vf_x}$ is the 
generalized Green function of $Z_G(t)^F$ 
(note that $Z_G(t)$ is connected).
\para{8.9.}
Recall that $\CX_{\io}$ is the function on $G^F\uni$ associated to
$\io$ given in 5.2.
It is known by [L3, V, (24.2.8)]
that $\CX_{\io}$ is expressed as 
\begin{equation*}
\tag{8.9.1}
\CX_{\io} = |\CW|\iv q^{-b(\io)}\sum_{w \in \CW}
           \Tr(w\iv, E_{\io})Q^G_{L^w, C_0^w, \CE_0^w, \vf_0^w},
\end{equation*}
where $\vf_0^w: F^*\CE_0^w\isom \CE_0^w$ is given by 
$\vf_0^w = \ad(\a\iv)^*\vf_0$ from $\vf_0: F^*\CE_0 \isom \CE_0$.
For each 
linear character $\th$ of $Z_{L^w}^{0F}$, we denote by $K^w_{\th}$
the complex $K^w$ given in 8.8 (subject to the condition that 
$z = 1$, i.e., $\vS = Z_L^0 \times C_0$) such that $\th_0^w = \th$, and 
denote by $\x_{K_{\th}^w}$ the characteristic function 
$\x_{K^w,\vf^w}$
Then by the character formula (8.8.3), we see that 
\begin{equation*}
\tag{8.9.2}
Q^G_{L^w, C_0^w, \CE_0^w, \vf^w_0} = |Z_{L^w}^{0F}|\iv
            \sum_{\th \in (Z_{L^w}^{0F})\wg}\x_{K^w_{\th}},
\end{equation*}
where we regard the left hand side as the class function on $G^F$
by extending by 0 outside of $G^F\uni$.
Under the isomorphism $Z_{L^w}^{0F} \simeq Z_L^{0Fw}$, $\th$
determines an $Fw$-stable  character $\th_1 \in (Z_L^{0F^m})\wg$ such
that $Sh_{F^m/Fw}(\th_1) = \th$.  Put 
$\CZ_{\th_1} = \{ w' \in \CW \mid Fw'(\th_1) = \th_1\}$.  
Then there exists a  $w_1 \in \CW$ such that 
$\CZ_{\th_1} = w_1\CW_{\th_1}$, and we define 
$\g_1: \CW_{\th_1} \to \CW_{\th_1}$ by $\g_1 = \ad w_1$.
Hence, as in the case of $K^w$, one can decompose $K^w_{\th}$ by 
\begin{equation*}
\tag{8.9.3}
\x_{K^w_{\th}} = q^{-b_0}\sum_{ E' \in (\CW_{\th_1}\wg)^{\g_1}}
                      \Tr(\g_1y, \wt E')\x_{A_{E',\th}}, 
\end{equation*} 
for $y \in \CW_{\th_1}$ such that $w = w_1y$, 
where $\x_{A_{E',\th}}$ is the (normalized) characteristic function of 
the character sheaf $A_{E',\th}$.
\para{8.10.}
For each irreducible character $\e$ of $Z_M(\la_c)^F$, put 
$\vG_{c,\e} = \Ind_{Z_M(\la_c)^F}^{G^F}(\e\otimes\vL_c)$.  Then we
have 
\begin{equation*}
\vG_c = \sum_{\e}(\deg\e)\vG_{c,\e},
\end{equation*}
where $\e$ runs over all the irreducible characters in $Z_M(\la_c)^F$.
We denote 
by $\vG_c'$ the sum of $\vG_{c,\e}$, where $\e$ runs over all the
linear characters of $Z_M(\la_c)^F$, and by $\vG_c''$ the complement
of $\vG_c'$ in $\vG_c$ so that $\vG_c = \vG_c' + \vG_c''$.
Then we have the following lemma.
\begin{lem}  
The inner product 
$\lp \vG_c'', Q^G_{L^w, C_0^w, \CE_0^w, \vf_0^w}\rp_{G^F}$ can be
expressed as 
\begin{equation*}
\lp \vG_c'', Q^G_{L^w, C_0^w, \CE_0^w, \vf_0^w}\rp_{G^F} =
    q^{-b_0+1}|Z_{L^w}^{0F}|\iv n_0\iv\b,
\end{equation*}
where $\b$ is an algebraic integer contained in a fixed cyclotomic 
field $\CA$ independent of $q$, and $n_0$ is an integer 
independent of $q$.
\end{lem}
\begin{proof}
By (8.9.2) and (8.9.3), we have
\begin{equation*}
\tag{8.11.1}
Q^G_{L^w,C_0^w,\CE_0^w, \f_0^w} = 
   |Z_{L^w}^{0F}|\iv q^{-b_0}\sum_{\th \in (Z_{L^w}^{0F})\wg}
      \sum_{E' \in (\CW_{\th_1}\wg)^{\g_1}}\Tr(\g_1y,\wt E')\x_{A_{E',\th}},
\end{equation*}
where $A_{E',\th}$ is the character sheaf which is a direct summand of 
$K^w_{\th}$.
Here we may assume that $\Tr(\g_1y, \wt E') \in \CA$.
Moreover, by Theorem 7.9, $\x_{A_{E',\th}}$ coincides with 
the almost character $R_x$ up to a scalar $\nu_x$, and we may assume that
$\nu_x$ is a unit in $\CA$.  On the other hand, it is known that 
$\deg\e$ is a
polynomial in $q$, and that $\deg\e$ is divisible by $q$ if $\e$ 
is not a linear character.  
It follows that $\lp \vG_c'', \r\rp_{G^F} \in n_1\iv q\BZ$ for any 
$\r \in \Irr G^F$.  In particular we have 
$\lp\vG_c'', R_x\rp_{G^F} \in n_2\iv q\CA$ for any almost character of 
$G^F$, where $n_1, n_2$ are integers independent of $q$. 
Then the lemma follows from (8.11.1).
\end{proof}
\par\medskip
\section{Proof of Theorem 8.6}
\para{9.1.}
In this section, we prove Theorem 8.6. 
First we note that (ii) follows from (i).  In fact, since 
$\r_{E}$ and $A_{E}$ have the 
same parametrizaition via the decomposition of 
$\Ind_{P^{F^m}}^{G^{F^m}}\d$ and of $\ind_P^GA_{\CE, z,\e}$,
we see that $R_E$ coincides with $\x_{A_E}$ up to scalar.  
Hence $R_E$ coincides with $R_{zz_E, \e_1\iv}$ by Theorem 8.6 (i).
This also shows that $\r_E = \r_{zz_E, \e_1\iv}$, and 
(ii) follows. 
\par  
In order to prove (i), first we show the following.
\begin{prop} 
Suppose that $q$ is large enough (but we don't assume that 
$q$ is sufficiently divisible).  Then the statement of 
(i) in Theorem 8.6 holds. 
\end{prop}
\para{9.3.}
The proof of the proposition will be done through 9.3 to 9.10.
We shall prove it by computing the inner product
$\lp\vG_{c,\xi}, \x_{A_E}\rp$ under the assumption that 
$q$ is large enough.  By Lemma 8.3, we have only to compare the 
inner products of $R_{z',\e'}$ and of $\x_{A_E}$ 
with various $\vG_{c,\xi}$ associated to the nilpotent element 
$N$ such that $N \in \CO_{\wt\r_{\ds_x, E'}}$.
\par
First we shall compute the inner product 
$\lp \vG_{c,\xi}, \x_{K^w,\vf^w}\rp_{G^F}$.  
By (8.8.3) and Proposition 3.5, we have 
\begin{align*}
\lp \vG_{c,\xi}, &\x_{K^w,\vf^w}\rp_{G^F} = 
   \frac{1} {|G^F|}\sum_{t', v' \in G^F} \frac{1}{|Z_G(t')^F|^2} 
      \sum_{\substack{g \in G^F \\ g\iv t'g \in Z_M(\la_c)^F}}
  \frac{|Z_{M}(\la_c)^F \cap Z_G(g\iv t'g)^F|}{|Z_M(\la_c)^F|} \\
&\times \sum_{\substack{x \in G^F \\ x\iv t'x \in \dz^wZ_{L^w}^{0F} }}
  \xi\nat(g\iv t'g)\ol{\bar\th_0^w(x\iv t'x(\dz^w)\iv)\p(\hat z_w)}
           \vG^{Z_G(t')}_{N_g, 1}(v')
             Q^{Z_G(t')}_{L_x^w, C_x^w,\CE_x^w, \vf_x}(v'),
\end{align*}
where in the first sum, $t', v' \in G^F$ runs over semisimple elements
$t'$, unipotent elements $v'$ such that $t'v' = v't'$.
Then the right hand side of the above expression can be written as 
\begin{align*}
 &\frac {1}{|G^F|}\sum_{t \in \dz^wZ_{L^w}^{0F}}
\frac 1 {|Z_G(t)^F|^2}\!\!\!\sum_{\substack{ g, x \in G^F\\ 
        g\iv xtx\iv g \in Z_M(\la_c)^F}}
 \!\!\!\frac{|Z_M(\la_c)^F \cap Z_G({}^{g\iv x}t)^F|}{|Z_M(\la_c)^F|}
\xi\nat({}^{g\iv x}t)\bar\th_0^w(t\iv\dz^w)\p(\hat z_w\iv) \\ 
  &\times\sum_{v \in Z_G(t)^F\uni}
      \vG^{Z_G(xtx\iv)}_{N_g,1}(xvx\iv)
    Q^{Z_G(xtx\iv)}_{L_x^w, C_x^w, \CE_x^w,\vf_x}(xvx\iv).
\end{align*}
Note that $\vG^{Z_G(xtx\iv)}_{N_g, 1}(xvx\iv) 
   = \vG^{Z_G(t)}_{N_{x\iv g}, 1}(v)$ and that
\begin{equation*}
Q^{Z_G(xtx\iv)}_{L_x^w, C_x^w, \CE_x^w,\vf_x}(xvx\iv) = 
 Q^{Z_G(t)}_{L^w, C_0^w, \CE_0^w, \vf^w_0}(v).
\end{equation*}
Hence by replacing $x\iv g$ by $g$ in the previous expression, we have
\begin{align*}
\tag{9.3.1}
\lp\vG_{c,\xi}, \x_{K^w,\vf^w}\rp_{G^F} 
   = &\sum_{t \in \dz^wZ^{0F}_{L^w}}
      \sum_{\substack{ g \in G^F \\ g\iv tg \in Z_M(\la_c)^F}}
  \frac{|Z_M(\la_c)^F \cap
    Z_G({}^{g\iv}t)^F|}{|Z_G(t)^F||Z_M(\la_c)^F|} \ \times \\
&\times \xi\nat_g(t)\bar\th_0^w(t\iv\dz^w)\p(\hat z_w\iv)
   \lp\vG^{Z_G(t)}_{N_g,1}, 
         Q^{Z_G(t)}_{L^w,C_0^w, \CE_0^w,\vf^w_0}\rp_{Z_G(t)^F},
\end{align*}
where $\xi_g\nat(t) = \xi\nat(g\iv tg)$.
\para{9.4.}
Returning to the original setting, we consider $A_E$ as in the
theorem, i.e., $A_E$ is a direct
summand of $K_{\th}^w$, where $\th \in (Z^{0F}_{L^w})\wg$ corresponds to 
$\th_0^w$ under the isomorphism $Z_{L^w}^{0F} \simeq Z_L^{0Fw}$. 
Hence $\CW_{\th_1}, \CZ_{\th_1}$, etc. are nothing but the objects 
discussed in 8.2.  So we write $w_{\d} = w_L$ and $\g_{\d} = \g$.
We continue the computation of
$\lp\vG_{c,\xi}, \x_{A_E}\rp$ under this setting.  
By (8.8.2) and (9.3.1), we have
\begin{align*}
\lp\vG_{c,\xi}, \x_{A_E}\rp &= |\CW_{\th_1}|\iv q^{b_0} 
       \sum_{ y \in \CW_{\th_1}}\Tr((\g y)\iv, \wt E) 
    \sum_{t \in \dz^wZ_{L^w}^{0F}}
   \sum_{\substack{g \in G^F \\ g\iv tg \in Z_M(\la_c)^F}}
  \!\!\!\frac{|Z_{M}(\la_c)^F \cap
       Z_G({}^{g\iv}t)^F|}{|Z_G(t)^F||Z_M(\la_c)^F|} \\
 &\times \xi_g\nat(t)\bar\th_0^w(t\iv\dz^w)\p(\hat z_w\iv)
   \lp\vG^{Z_G(t)}_{N_g,1}, 
         Q^{Z_G(t)}_{L^w,C_0^w, \CE_0^w,\vf^w_0}\rp_{Z_G(t)^F}.
\end{align*}
Let $H$ be an $F$-stable reductive subgroup of $G$ containing 
$L$.  Put $\CW_H = N_H(L)/L \subset \CW$. 
Recall that $\th_1 \in (Z_L^{0F^m})\wg$, and 
$\CZ_{\th_1} = \{ w \in \CW \mid {}^{Fw}\th_1 = \th_1\} = w_L\CW_{\th_1}$.  
Then $\CW_H \cap \CZ_{\th_1} =w_L'\CW_{H,\th_1}$, where 
$w_L' = w_Lw_H \in \CW_H$ with $w_H \in \CW_{\th_1}$ and $\CW_{H,\th_1}$ is
the stabilizer of $\th_1$ in $\CW_H$. 
Then $H$ contains $L^{w_L'y}$ for any  $y \in \CW_{H,\th_1}$.
Put 
\begin{align*}
\tag{9.4.1}
\vD_H^{(g)} &= |\CW_{\th_1}|\iv q^{b_0}\sum_{y \in \CW_{H,\th_1}}
       \Tr((\g w_Hy)\iv, \wt E)\ \times \\ 
   &\times \sum_{t}
\xi\nat_g(t)\bar\th_0^w(t\iv\dz^w)\p(\hat z_w\iv)
   \lp\vG^{H}_{N_g,1}, 
         Q^{H}_{L^w,C_0^w, \CE_0^w,\vf^w_0}\rp_{H^F},   
\end{align*}
where in the second sum, $t$ runs over all the elements in 
$\dz^wZ_{L^w}^{0F}$ with $w = w_L'y = w_Lw_Hy$ 
such that $Z_G(t) = H$ and that $g\iv tg \in Z_M(\la_c)$ 
for $g \in G^F$.
Then we have 
\begin{equation*}
\tag{9.4.2}
\lp\vG_{c,\xi}, \x_{A_E}\rp = 
    \sum_{H}\sum_{g}\vD^{(g)}_H, 
\end{equation*}
where in the first sum, $H$ runs over all the $F$-stable reductive subgroups of 
$G$ containing $L$ such that $H = Z_G(t)$ for a fixed $t \in Z_L$, 
and in the second sum  
$g$ runs over all the elements in the double cosets 
$H^F\backslash G^F/Z_M(\la_c)^F$ such
that $g\iv tg \in Z_M(\la_c)^F$.
Here we note the following lemma, which is a stronger version of 
Lemma 8.11.
\begin{lem}
Assume that $q$ is large enough. 
Let $\vG_c''$ be as in 8.10. Then we have 
\begin{equation*}
\lp \vG_c'', Q^G_{L^w, C_0^w, \CE_0^w, \vf_0^w}\rp_{G^F} 
       = q^{-b_0+1}n_0\iv\b,
\end{equation*}
where $\b \in \CA$, and $n_0$ is an integer independent of $q$.
\end{lem}
\begin{proof}
First we note that 
\begin{equation*}
\tag{9.5.1}
\lp \vG_{c,1}, Q^G_{L^w, C_0^w, \CE_0^w, \vf_0^w}\rp_{G^F} =
    q^{-b_0}n_1\iv \b',
\end{equation*}
where $\b'$ is an algebraic integer contained in a fixed cyclotomic
field $\CA$ independent of $q$, and $n_1$ is an integer independent of $q$.
We show (9.5.1).  We apply (9.3.1) to the situation that $\xi = 1$ 
and $\th_1 = 1$, i.e, $\CW_{\th_1} = \CW$.  Then by a similar argument
as in the proof of (9.4.2), we have
\begin{equation*}
\tag{9.5.2}
\lp\vG_{c,1}, \x_{K^w, \vf^w}\rp_{G^F} = 
   |Z_G^F|\lp\vG_{c,1}, Q^G_{L^w,C_0^w, \CE_0^w,
          \vf_0^w}\rp_{G^F} + \sum_{H \ne G}\sum_g \Xi_H^{(g)},
\end{equation*} 
where $\Xi_H^{(g)} = 
   \sum_t\lp\vG^{H}_{N_g,1},Q^{H}_{L^w,C_0^w,\CE_0^w,\vf^w_0}\rp_{H^F}$.
By induction on the rank of $G$, we may assume that 
$\Xi_H^{(g)}$ can be expressed as 
$\Xi_H^{(g)} = q^{-b_0}n_H\iv\b'_H$, where 
$n_H, \b'_H$ are similar elements as $n_1, \b'$ in (9.5.1). 
(Note that $b_0$ has common value for all $H$ containing $L$.) 
On the other hand, by Theorem 7.9, $\x_{A_E}$ coincides with 
an almost character $R_x$ of $G^F$ up to a scalar $\nu_x$ which
is a root of unity in $\CA$.   
It follows, by (7.9.2), that 
\begin{equation*}
\lp\vG_{c,\xi}, \x_{K^w,\vf^w}\rp_{G^F} 
     \in q^{-b_0}n_2\iv\CA
\end{equation*}
with some 
$n_2 \in \BZ$ independent of $q$.  (9.5.1) now follows from (9.5.2). 
\par
Next we note that  
\begin{equation*}
\tag{9.5.3}
\lp\vG_{c}, Q^G_{L^w,C_0^w, \CE_0^w, \vf_0^w}\rp_{G^F} \in q^{-b_0}\CA.
\end{equation*}
In fact, by Lemma 5.12 (iii), we have 
\begin{equation*}
\lp\vG_{c}, \CX_{\io}\rp_{G^F} \in q^{-b_0-b(\io)}\CA.
\end{equation*}
Since
\begin{equation*}
Q^G_{L^w, C_0^w, \CE_0^w, \vf_0^w} = 
\sum_{\io \in \CI_0^F}\Tr(w, E_{\io})q^{b(\io)}\CX_{\io}
\end{equation*}
by (8.9.1), we obtain (9.5.3).
\par
Now we have
\begin{equation*}
\lp\vG_c'', Q^G_{L^w,C_0^w, \CE_0^w, \vf_0^w}\rp_{G^F} = 
\lp\vG_{c}, Q^G_{L^w,C_0^w, \CE_0^w, \vf_0^w}\rp_{G^F} - 
      f_{N_c}\lp\vG_{c,1}, Q^G_{L^w,C_0^w, \CE_0^w, \vf_0^w}\rp_{G^F},
\end{equation*}
where $f_{N_c}$ is the number of linear characters of $Z_M(\la_c)^F$.
Hence by (9.5.1) and (9.5.3), we have 
\begin{equation*}
\tag{9.5.4}
\lp \vG_c'', Q^G_{L^w, C_0^w, \CE_0^w, \vf_0^w}\rp_{G^F} 
             = q^{-b_0}n_1\iv\b'
\end{equation*}
with $n_1 \in \BZ$ independent of $q$, and $\b' \in \CA$.
\par
On the other hand, by Lemma 8.11, we have 
\begin{equation*}
\tag{9.5.5}
\lp \vG_c'', Q^G_{L^w, C_0^w, \CE_0^w, \vf_0^w}\rp_{G^F} =
    q^{-b_0+1}|Z_{L^w}^{0F}|\iv n_2\iv\b''
\end{equation*}
with $\b'' \in \CA$.
We may assume that 
$\lp \vG_c'', Q^G_{L^w, C_0^w, \CE_0^w, \vf_0^w}\rp_{G^F} \ne 0$.
Then by (9.5.4) and (9.5.5), we see that 
\begin{equation*}
n_2\b' = q|Z_{L^w}^{0F}|\iv n_1\b'' \in \CA.
\end{equation*}
Since $q$ and $|Z_{L^w}^{0F}|$ are prime to each other, 
$n_1\b''$ is divisible by $|Z_{L^w}^{0F}|$.
Hence $\b'$ can be written as $\b' = qn_2\iv \b$ with 
$\b \in \CA$, and we have 
\begin{equation*}
\lp \vG_c'', Q^G_{L^w, C_0^w, \CE_0^w, \vf_0^w}\rp_{G^F} =
                q^{-b_0+1}(n_1n_2)\iv\b.
\end{equation*}  
This proves the lemma.
\end{proof}
\para{9.6.}
We continue the computation in 9.4.
By (8.9.1) applied to the 
reductive group $H$, we have
\begin{equation*}
Q^{H}_{L^w, C^w_0, \CE_0^w,\f^w_0} = \sum_{\io \in (\CI'_0)^F}
     \Tr(w, E_{\io})q^{b_{H}(\io)}\CX^{H}_{\io}
\end{equation*}
for any $w \in \CW_H$, where
$b_{H}(\io)$ is given as in (8.8.1) by replacing $G$ by $H$, 
and $(\CI_H)_0$ is a block in $\CI_H$ corresponding to $(L, C_0,\CE_0)$.  
Substituting this into the previous formula, we see that 
\begin{equation*}
\tag{9.6.1}
\vD^{(g)}_{H} = q^{b_0}\sum_{\io \in (\CI_H)_0^F}
   \sum_{\substack{ t \in z_1Z_{L_1}^{0F} \\
                    Z_G(t) = H}}
        a_{\io,H}(t)\xi\nat_g(t)q^{b_{H}(\io)}
             \lp\vG^{H}_{N_g,1}, \CX^{H}_{\io}\rp, 
\end{equation*}
where $L_1 = L^{w_L'}$ and $z_1 = \dz^{w_L'}$.  
Here $a_{\io,H}$ is given as follows; 
take $\a \in H$ such that $\a\iv F(\a) = \dw_L' \in N_H(L)$.  Then 
we have $L_1 = \a L\a\iv$ and $z_1 = \a \dz_{w_L'}\a\iv$. 
By putting 
$t_1 = \a\iv t \a \in (\dz Z_L^0)^{Fw_L'}$ we have
\begin{equation*}
\tag{9.6.2}
a_{\io,H}(t) = |\CW_{\th_1}|\iv\sum_{y \in \CW_{H,\th_1}}
        \th_0^{w}(t_1\iv\dz_{w})\p(\hat z_{w}\iv)
                \Tr((\g w_Hy)\iv, \wt E)\Tr(w, E_{\io}) 
\end{equation*}
with $w = w_L'y = w_Lw_Hy$.
\par
In the formula (9.6.1), we shall replace the inner product
$\lp\vG^{H}_{N_g,1},\CX_{\io}^{H}\rp$ by 
$\lp\vG^{H}_{N_g}, \CX_{\io}^{H}\rp$, which is easier to 
handle with.
Now it is easy to see, by making use of Lemma 9.5, that 
\begin{equation*}
\lp \vG_c'', \CX_{\io}\rp_{G^F} = q^{-b_0 - b(\io) +1}n_0\iv\b_{\io}
\end{equation*}
for an appropriate integer $n_0$ and $\b_{\io} \in \CA$.
We also note the relation 
\begin{equation*}
\lp\vG_c, \CX_{\io}\rp_{G^F} = f_{N_c}\lp\vG_{c,1},\CX_{\io}\rp_{G^F}
                 + \lp\vG_c'', \CX_{\io}\rp_{G^F}.
\end{equation*}
Applying these formulas to the case of $H$, we have 
\begin{equation*}
\tag{9.6.3}
\vD^{(g)}_{H} = \sum_{\io \in (\CI_H)_0^F}
   \sum_{\substack{ t \in z_1Z_{L_1}^{0F} \\
                    Z_G(t) = H}}
        a_{\io,H}(t)\xi\nat_g(t)f_{N_g}\iv\bigl\{q^{b_{H}(\io)+b_0}
             \lp\vG^{H}_{N_g}, \CX^{H}_{\io}\rp + qn_H\iv\b_{H,\io}\big\}, 
\end{equation*}
where $N_g = \ad(g)N_c$ is a nilpotent element in $\Lie H$, 
$n_H$ is an integer independent of $q$, and $\b_{H,\io} \in \CA$. 
Moreover, $f_{N_g}$ is the number of linear characters of 
$Z_{{}^gM \cap H}({}^g\la)^F$.
\par
Next we shall compute $\lp\vG^{H}_{N_g}, \CX^{H}_{\io}\rp$.
Recall that $Z_G(g\iv tg) \supset L$, and 
$N_g = \Ad(g)N_c \in \Lie H$. It follows that 
$\Ad(g)N \in \Lie H$.  Since $g \in G^F$, we may replace
$N$ by $N_1 = \Ad(g)N$ in parameterizing the nilpotent orbits $\CO_N^F$ in 
$\Fg$. Then we can identify $A_{H,\la} =
Z_{H}(N_1)/Z^0_{H}(N_1)$ as a subgroup of $A_G(N)$. 
As in 5.8, there exists $c_0 \in A_{\la}$ such that $-N_c^*$ is 
$G^F$-conjugate to $N_{cc_0}$.  Then $N_g = \Ad(g)N_c$ is
$H^{F}$-conjugate to $(N_1)_c$, and so $-N^*_g = \Ad(g)(-N_c^*)$ is
$H^{F}$-conjugate to $(N_1)_{cc_0}$, with $c, c_0 \in A_{H, \la}$.
In particular, for $\io'' \in \CI_H$ such that $\supp \io'' = \CO_{N_g}$, 
we see that $\CY_{\io''}(-N_g^*) = \e(cc_0)$, which is
independent from $H$. 
Let $\CO_1$ be the nilpotent orbit in $\Lie H$ such that 
$\supp \io' = \CO_1$ with $E_{\io'} = E_{\io} \otimes \ve$.
Then we have 
\begin{equation*}
\tag{9.6.4}
q^{b_{H}(\io) + b_0}\lp\vG^{H}_{N_g}, \CX^{H}_{\io}\rp = 
    \begin{cases}
        \z_{\CI_0}\iv\e(cc_0)\iv 
               &\quad\text{ if }  \CO_{N_g} = \CO_1, \\
        q\b'_{H,\io}
               &\quad\text{ if }   \CO_{N_g} \ne \CO_1,
    \end{cases}        
\end{equation*}
for some $\b'_{H,\io} \in \CA$.
\par
In fact, by Lemma 5.12 (i) applying to $H$, we have
\begin{equation*}
q^{b_{H}(\io) + b_0}\lp\vG^{H}_{N_g}, \CX^{H}_{\io}\rp = 
    \z_{\CI_0}\iv q^{(\dim\supp(\io') -\dim\supp(\io''))/2}
         \BP_{\io'',\io'}(q\iv)\ol{\CY_{\io''}(-N_g^*)}.
\end{equation*}
(Note that $\z_{(\CI_H)_0}$ for $\CI_H$ coincides with $\z_{\CI_0}$
for $\CI_G$ 
since both of them correspond to $(L, C_0,\CE_0)$.) 
If $\supp(\io') = \supp(\io'')$, then we have $\io' = \io''$ and 
$\BP_{\io'',\io'} = 1$.  This implies the first equality.  
Now assume that $\io' \ne \io''$.
By Lemma 5.12 (ii), we have 
\begin{equation*}
q^{(\dim\supp(\io') -\dim\supp(\io''))/2}
         \BP_{\io'',\io'}(q\iv) \in q\BZ.
\end{equation*} 
Since we may assume that 
$\z_{\CI_0}\iv, \ol{\CY_{\io''}(-N_g^*)} \in \CA$, we obtain 
the second equality.
\para{9.7.}
We consider the sum of $a_{\io,H}(t)\xi\nat_g(t)f_{N_g}\iv$ 
in (9.6.3) for a fixed $\io$.
Recall that $L_1 = \a L\a\iv$ and that $t$ is an element of 
$z_1Z_{L_1}^{0F}$ such that $Z_G(t) = H$ and
that $g\iv tg \in Z_M(\la_c)$.  Put $M_{H} = gMg\iv \cap H$.
Then $M_{H}$ is the Levi subgroup in  
$H$ associated to the nilpotent element $N_g$, and 
$t \in (Z_{M_{H}}(N_g) \cap z_1Z_{L_1}^0)^F$.    
Put $M'_H = \a\iv M_{H}\a = {}^{\a\iv g}M \cap H$, which 
is stable by $F''_H = Fw_L'$, and put 
\begin{equation*}
X_{M_H} = (Z_{M'_H}(N_{\a\iv g})\cap \dz Z_{L}^{0})/
           (Z_{M'_H}^1(N_{\a\iv g})\cap Z_{L}^{0}).
\end{equation*}
Then  
\begin{equation*}
\tag{9.7.1}
(Z_{M_{H}}(N_g) \cap z_1Z_{L_1}^{0})/
      (Z^1_{M_{H}}(N_g) \cap Z_{L_1}^{0}) \simeq X_{M_H}
\end{equation*} 
and the action of $F$ on the left hand side corresponds to the
action of $Fw_L'$ on $X_{M_H}$.
\par
We apply Lemma 8.5 (i) to $H$.  Then 
one can write, for $w = w_L'y$, 
\begin{equation*}
\th_0^{w}(t_1\iv\dz_w)\p(\hat z_w\iv) = 
        \th_0^{w'_L}(t_1\iv\dz_{w'_L})
          \p(\hat z_{w_L'}\iv)\Psi_{t_1\iv}(y),
\end{equation*} 
where $\Psi_{t_1\iv}$ is a
linear character of $\CW_{H,\th_1}$.
Then (9.6.2) can be rewritten as 
\begin{equation*}
\tag{9.7.2}
a_{\io,H}(t) = \frac{|\CW_{H,\th_1}|}{|\CW_{\th_1}|}
    \th_0^{w'_L}(t_1\iv\dz_{w_L'})\p(\hat z_{w_L'}\iv)
       \lp \wt E\otimes \Psi_{t_1\iv}, E_{\io}\rp_{w'_L\CW_{H,\th_1}}, 
\end{equation*}
where $\lp \ , \ \rp_{w'_L\CW_{H,\th_1}}$ is the inner product 
on the $\CW_{H,\th_1}$-invariant functions on $w'_L\CW_{H,\th_1}$.
Put $Y = Z_{M_{H}}(N_g) \cap z_1Z_{L_1}^0$ and 
$Y_1 = Z_{M_{H}}^1(N_g)\cap Z_{L_1}^0$.
Then by Lemma 8.5 (ii), together with the isomorphism in (9.7.1), 
$\Psi_{t_1}$ depends only on $t \mod Y_1^F$.
\par
On the other hand, we have an isomorphism 
\begin{equation*}
Y/Y_1 \simeq 
   (Z_{M \cap {}^{g\iv}H}(\la_c)\cap \dz Z^0_{{}^{g\iv }L_1})/
   (Z^1_{M \cap {}^{g\iv}H}(\la_c)\cap Z^0_{{}^{g\iv}L_1})
\end{equation*}
and a natural map 
\begin{equation*}
  (Z_{M \cap {}^{g\iv}H}(\la_c)\cap \dz Z^0_{{}^{g\iv}L_1})/
   (Z^1_{M \cap {}^{g\iv}H}(\la_c)\cap Z^0_{{}^{g\iv}L_1}) \to 
         Z_M(\la_c)/Z_M^1(\la_c) \simeq \bar A_{\la}.
\end{equation*}
Hence we have a map
\begin{equation*}
\tag{9.7.3}
Y^F/Y_1^F \to (Y/Y_1)^F \to \bar A_{\la}^F
\end{equation*}
satisfying the property that $\xi\nat(g\iv tg) = \xi(a)$, where 
$a \in \bar A_{\la}^F$ is the image of $t \in Y^F$ under this map.
In particular, we see that $\xi\nat_g(t)$ also depends only on 
$t \mod Y_1^F$.
\par
Put $\bar t = t \mod Y_1^F \in Y^F/Y_1^F$.
It follows from the above argument, together with (9.7.2), that 
\begin{align*}
\tag{9.7.4}
\sum_{ t \in \bar t}&a_{\io,H}(t)\xi\nat_g(t) \\ 
   &= \frac{|\CW_{H,\th_1}|}{|\CW_{\th_1}|}\xi\nat_g(\bar t)
         \lp \wt E\otimes \Psi_{\bar t\iv}, E_{\io}\rp_{w'_L\CW_{H,\th_1}}
      \p(\hat z_{w_L'}) \sum_{t \in \bar t}\th_0^{w'_L}(t_1\iv z_{w_L'}),
\end{align*}
where $\Psi_{\bar t\iv}$ (resp. $\xi\nat_g(\bar t)$) denotes 
$\Psi_{t_1\iv}$ (resp. $\xi\nat_g(t)$) for $t \in \bar t$.
Since $\th_0^{w_L'}$ is a linear character on $Z_L^{0Fw_L'}$, we have
\begin{equation*}
\tag{9.7.5}
\sum_{ t \in \bar t}\th_0^{w_L'}(t_1\iv z_{w_L'}) = 
       \begin{cases}
          |Y_1^F|\bar\th_0^{w_L'}(\bar t\iv z_1)
        &\quad\text{ if $\bar\th_0^{w_L'}|_{Y_1^F}$ is
                  trivial},  \\
          0       &\quad\text{ otherwise},
       \end{cases}
\end{equation*}
where $\bar\th_0^{w_L'}(\bar t\iv z_1)$ denotes a common value of 
$\th_0^{w_L'}(t_1\iv z_{w_L'})$ for $t \in \bar t$. 
 
\par
By definition, $f_{N_g}$ is the number of linear characters of
$Z_{M_{H}}(N_g)^F$, which coincides with the order of 
$Z_{M_{H}}(N_g)^F/(Z_{M_{H}}(N_g)^F)\der$.
Hence one can write $f_{N_g}$ as 
\begin{equation*}
f_{N_g} = f_H|Z_1^{0F}|,
\end{equation*}
where $Z_1 = Z(Z_{M_{H}}(N_g))$ is the center of $Z_{M_{H}}(N_g)$, 
and $f_H$ is contained in a finite subset of $\BQ$ independent of $q$. 
On the other hand, the description of $Y_1$ in 8.4 (applied to the case
$M_{H}$) implies that $Z_1^0 \subseteq Y_1$.
Thus $Y_1^F$ is divisible by $Z_1^{0F}$ and $|Y_1^F|/|Z_1^{0F}|$ is
a polynomial in $q$ (in fact, it is the order of the group of
$F$-fixed points of the group $Y_1/Z_1^0$ whose connected component 
is a torus). Let $g_H$ be the constant term of the polynomial 
$|Y_1^F|/|Z_1^{0F}|$.
\par
Note that we may assume that $\th_0^{w}(t) \in \CA$ for any 
$w \in \CW$ since $\th_0^w$ comes from a linear character $\th_1$
of $(L/L\der)^{F^m}$, and $L/L\der$ is a finite group. 
Thus summing up the above argument, we have the following formula.
\begin{equation*}
\tag{9.7.6}
\begin{split}
\sum_{t \in Y^F}&a_{\io,H}(t)\xi\nat_g(t)f_{N_g}\iv 
  = n_1\iv q\b  \\
      &+ \frac{|\CW_{H,\th_1}|}{|\CW_{\th_1}|}f_H\iv g_H\p(\hat z_{w_L'}\iv)
       \sum_{ \bar t \in Y^F/Y_1^F}
        \bar\th_0^{w_L'}(\bar t\iv z_1)\xi\nat_g(\bar t)
         \lp \wt E\otimes \Psi_{\bar t\iv}, E_{\io}\rp_{w'_L\CW_{H,\th_1}}, 
\end{split}
\end{equation*}
where $n_1 \in \BZ$ is independent of $q$, and $\b \in \CA$. 
(We understand that $\bar\th_0^{w_L'}(\bar t\iv z_1) = 0$ if 
$\bar\th_0^{w_L'}|_{Y_1^F}$ is non-trivial).
\para{9.8.}
We return to the setup in 9.6. 
In view of (9.6.4) and (9.7.6), one can rewrite the equation (9.6.3)
in the form
\begin{equation*}
\tag{9.8.1}
\begin{split}
\vD^{(g)}_H = n_1\iv q\b&   
  + \frac{|\CW_{H,\th_1}|}
 {|\CW_{\th_1}|}f_H\iv g_H\z_{\CI_0}\iv\e(cc_0)\iv \p(\hat z_{w_L'}\iv) \times \\
       &\times \sum_{ \bar t \in Y^F/Y_1^F}
      \bar\th_0^{w_L'}(\bar t\iv z_1)\xi\nat_g(\bar t)
         \lp \wt E\otimes \Psi_{\bar t\iv}, E_{\io}\rp_{w'_L\CW_{H,\th_1}}, 
\end{split}
\end{equation*}
where $n_1 \in \BZ$ is independent of $q$ and $\b \in \CA$, and $\io$
is the unique element in $(\CI_{H})_0$ such that 
$\supp (\io') = \CO_{N_g}$ with $E_{\io'} = E_{\io}\otimes\ve$.
\par
We now compute the inner product 
$\lp \wt E\otimes \Psi_{\bar t\iv}, E_{\io}\rp_{w'_L\CW_{H,\th_1}}$
in the right hand side of (9.8.1).
Note that $E \in \CW_{\th_1}\wg$ is  
an extension of the $\Om_{\th_1}$-stable character  
$E_1$ of $\CW_{\th_1}^0$ to $\CW_{\th_1}\rtimes \Om_{\th_1}$,
where $E_1 \in (\CW_{\th_1}^0)\wg$ is of the form  
$E_1 = E^{\mu}\boxtimes\cdots\boxtimes E^{\mu}$ 
($t/d$-times) with 
$E^{\mu} \in \FS_{n/t}\wg$ corresponding to a partition 
$\mu$ of $n/t$. 
Let $\mu^*$ be the partition dual to $\mu$.
Then $N = N_{\la}$ with $\la = t\mu^*$ by 8.2.
(In general, 
for a partition $\r = (\r_1, \dots, \r_k)$ and 
$a \in \BZ_{>0}$, we put $a\r = (a\r_1, \dots, a\r_k)$). 
By our assumption, $\supp(\io') = \CO_{N_g}$. 
Then $\CW_H$ is of the form 
\begin{equation*}
\tag{9.8.2}
\CW_H \simeq \FS_{\nu_1}\times\cdots\times \FS_{\nu_k},
\end{equation*} 
where $\nu = (\nu_1, \dots, \nu_k)$ is a partition of $n/d$
such that $(t/d)\mu^*$ is a refinement of $\nu$, i.e., 
$\nu_i$ is a sum of parts of $(t/d)\mu^*$.
Moreover 
$\CW_{H,\th_1}$ is given as 
\begin{equation*}
\tag{9.8.3}
\CW_{H,\th_1} \simeq \CW^0_{H,\th_1}\rtimes \Om_{\th_1}
\end{equation*}
with
\begin{equation*}
\CW^0_{H,\th_1} \simeq \FS_{\nu'} \times \cdots \times \FS_{\nu'}
\quad (t/d\text{-times}),
\end{equation*}
where $\nu'$ is a partition of $n/t$ such that $(t/d)\nu' = \nu$.
Since $\nu' = (\nu_1', \dots, \nu_k') $ is 
a partition of $n/t$ whose parts are sums of parts of $\mu^*$, 
$\mu^*$ determines a partition $(\mu^{(i)})^*$ of $\nu_i'$ for 
$i = 1, \dots, k$.  We denote by $\mu^{(i)}$ its dual partition.
Hence $\mu = \mu^{(1)} + \cdots  + \mu^{(k)}$.
We define an irreducible character $E^0$ of $\CW^0_{H,\th_1}$ by 
$E^0 = E_1^0\boxtimes\cdots\boxtimes E_1^0$ with 
$E^0_1 \in \FS_{\nu'}\wg$ such that 
\begin{equation*}
E^0_1 = E^{\mu^{(1)}}\boxtimes\cdots\boxtimes E^{\mu^{(k)}}.
\end{equation*}
We remark that 
\par\medskip\noindent
(9.8.4) \ $\lp E, E_{\io}\rp_{\CW^0_{H,\th_1}} = 1$,  and 
$E^0$ is the unique irreducible character of $\CW^0_{H,\th_1}$ which 
appears in the restrictions to $\CW^0_{H,\th_1}$ 
of both of $E$ and $E_{\io}$.
\par\medskip
In fact, one can write $E|_{\CW^0_{H,\th_1}} = E^0 + \sum E'$
with $a(E') > a(E^0)$.  
On the other hand, $E_{\io}$ is given as 
$E_{\io} = E^{\r^{(1)}}\boxtimes\cdots\boxtimes E^{\r^{(k)}}$ with 
$E^{\r^{(i)}} \in \FS_{\nu_i}\wg$, where $\r^{(i)}$ is a partition 
of $\nu_i$ such that $(\r^{(i)})^* = (t/d)(\mu^{(i)})^*$.  It follows 
that $\r^{(i)} = \mu^{(i)} \cup \cdots \cup \mu^{(i)}$ ($t/d$-times).
(For a partition $\la$ and $\mu$, we denote by $\la \cup \mu$ the
partition obtained by rearranging the parts of $\la$ and $\mu$ in 
decreasing order.) Then one can write as 
$E_{\io}|_{\CW^0_{H,\th_1}} = E^0 + \sum E''$ with 
$a(E'') < a(E^0)$.  Hence (9.8.4) holds.
\par
Since $E^0$ is $\Om_{\th_1}$-stable, (9.8.4) implies that there 
exists a unique extension  $\wt E^0$ of $E^0$ to $\CW_{H,\th_1}$ 
which appears in the decomposition of $E_{\io}|_{\CW_{H,\th_1}}$ 
with multiplicity one.  On the other hand, again by (9.8.4), the 
restriction of $E$ to $\CW_{H,\th_1}$ also contains a certain 
extesnion $(\wt E^0)'$ of $E^0$ with multiplicity one.  Hence 
we can write $(\wt E^0)' = \wt E^0\otimes \w$ with 
some $\w \in \Om_{\th_1}\wg$.
Since
$(\CI_{H})^F_0 = (\CI_H)_0^{F''_H}$,  
$E_{\io}$ is stable by $F''_H$. Since $E$ is also $F''_H$-stable, 
we see that $\wt E^0$ and $(\wt E^0)'$ are  $F''_H$-stable.  
It follows that $\w$ is $F''_H$-stable.
\par
Now by applying Lemma 8.5 (ii) to the group $H$, we see that
\par\medskip\noindent
(9.8.5) \ There exists a unique class 
$x_H \in  X_{M_H}^{F''_H}$ 
satisfying the following. 
\begin{equation*}
\lp E\otimes\Psi_{t\iv}, E_{\io}\rp_{\CW_{H,\th_1}} = 
     \begin{cases}
        1 &\quad\text{ if } t \in x_H, \\
        0 &\quad\text{ otherwise.}
     \end{cases}
\end{equation*}
\par\medskip
We pass to the extension $\wt E$.  Then the above 
arguments show that the restrictions of 
$\wt E\otimes \Psi_t\iv$ and $E_{\io}$ to $\lp w_L'\rp \CW_{H,\th_1}$ 
contain a unique irreducible character which is an extension of 
$\wt E^0$ to $\lp w_L'\rp\CW_{H,\th_1}$ for $t \in x_H$.   
\par
Under the notation in 9.7, $Y^F/Y_1^F$ is regarded as a subset of
$M_H^{F''_H}$.
If $x_H$ is contained in $Y^F/Y_1^F$, $x_H$ determines an element in 
$\bar A_{\la}^F$ by (9.7.3), which we denote by $z_H \in \bar A_{\la}^F$.
Then summing up the above arguments, we have 
\begin{equation*}
\tag{9.8.6}
\begin{split}
\sum_{\bar t \in Y^F/Y_1^F}&\bar\th_0^{w_L'}(\bar t\iv z_1)\xi\nat_g(\bar t)
       \lp\wt E\otimes\Psi_{\bar t\iv}, E_{\io}\rp_{w_L'W_{H,\th_1}} \\
   &= \begin{cases}
          \bar\th_0^{w_L'}(x_H\iv z_1)\xi(z_H)\a_H
                &\quad\text{ if } x_H \in Y^F/Y_1^F, \\
           0    &\quad\text{ otherwise, }
      \end{cases}
\end{split}
\end{equation*}
where $\a_H$ is a root of unity determined by the extension
$\wt E$, which is independent of $q$.   
Substituting this into (9.8.1), we have
\begin{equation*}
\tag{9.8.7}
\vD^{(g)}_H \in \d_H\frac{|\CW_{H,\th_1}|}{|\CW_{\th_1}|}
             f_H\iv g_H\z_{\CI_0}\iv\e(cc_0)\iv
                  \p(\hat z_{w_L'}\iv)\bar\th_0^{w_L'}(x_H\iv z_1)
                    \xi(z_H)\a_H + n_1\iv q\CA,
\end{equation*}
where $\d_H = 1$ if $x_H \in Y^F/Y_1^F$ and 
$\bar\th_0^{w_L'}|_{Y_1^F}$ is trivial, and 
$\d_H = 0$ otherwise.
\para{9.9.}
We shall compare $x_H$ for various $H$ appearing in 9.8. 
Since 
\begin{equation*}
Z_{M'_H}(N_{\a\iv g}) \cap \dz Z_L^0 = 
    Z_{{}^{\a\iv g}M}(N_{\a\iv g}) \cap \dz Z_L^0
\end{equation*}
and similarly for $Z^1_{M'_H}(N_{z\iv g}) \cap Z_L^0$, 
a similar argument as in 8.4 implies that 
\begin{equation*}
X_{M_H} \simeq \BZ/(t/d)\BZ \simeq X_M.
\end{equation*}
We have a natural isomorphism $f: X_{M_H} \to X_M$ which is 
compatible with the action of $F''_H$ and of $F''$ (the both
coincide with the action of $F$ on $\BZ/(t/d)\BZ$). 
Assume that $H$ satisfies the property (9.8.2).
Then it is easy to check that $w_L \in \CW_H$, and so we have
\begin{equation*}
\tag{9.9.1}
w_L = w_L'.
\end{equation*}
We note that 
\begin{equation*}
\tag{9.9.2}
f(x_H) = x_G.
\end{equation*}
In fact, let $E_{\io}\in \CW_H\wg, E^0 \in (\CW^0_{H,\th_1})\wg$ 
be as in 9.8 for $H$, and $E_{\io_G} \in \CW\wg$, 
$E_G^0 \in (\CW^0_{\th_1})\wg$ the corresponding objects for $G$.
  Here $\io$ (resp. $\io_G$) is the unique element 
in $(\CI_H)_0$ (resp. $\CI_0$ ) such that $\supp \io' = \CO_{N_g}$ 
(resp. $\supp \io'_G = \CO_N$).  
Then $E_{\io}$ is contained in $E_{\io_G}|_{\CW_H}$ with 
multiplicity one.
On the other hand, 
we have $E|_{\CW^0_{\th_1}} = E_G^0$, and 
the restriction of $E_G^0$ on $\CW_{H,\th_1}^0$ contains $E^0$ 
with multiplicity one by the property (9.8.4).  Let $\wt E^0$  be the
extension of $E^0$ to $\CW_{H,\th_1}$ appearing in 
$E_{\io}|_{\CW_{H,\th_1}}$, and let $\wt E_G^0$ a similar object 
as $\wt E^0$ for $E^0_G$.  Then the above fact shows that 
$\wt E^0$ occurs in the restriction of $\wt E_G^0$ with multiplicity 
one.  This implies that  $E$ coincides with $\wt E^0_G \otimes \w$, 
where $\w \in \Om_{\th_1}\wg$ is given as in 9.8 for $H$. 
Thus (9.9.2) is proved.
\par
It follows from (9.9.2) that $x_H$ in (9.8.5) depends only on $E$, 
and does not depend on the choice of $H$, which we denote by $x_E$.  
We also denote by $z_E$ the element $z_H \in \bar A_{\la}^F$ determined 
from $x_H$ (see (9.8.6)).   
A similar argument as above shows that the root of unity 
$\a_H$ in (9.8.6) also depends only on $E$, and not on $H$, which 
we denote by $\a_E$.
\par
Summing up the above argument, (9.8.7) can be written as 
\begin{equation*}
\tag{9.9.3}
\vD^{(g)}_H \in \d_H\frac{|\CW_{H,\th_1}|}{|\CW_{\th_1}|}
             f_H\iv g_H\z_{\CI_0}\iv\e(cc_0)\iv
       \p(\hat z_{w_L'}\iv)\bar\th_0^{w_L}(x_E\iv z_1)\xi(z_E)\a_E + qn_H\iv\CA,
\end{equation*}
where $\d_H = 1$ if $x_E \in Y^F/Y_1^F$ and 
$\bar\th_0^{w_L}$ is trivial on $Y_1^F$, and
$\d_H = 0$ otherwise, and
$n_H$ is an integer independent of $q$.  
\para{9.10}
We are now ready to prove Proposition 9.2, by completing   
the computation of $\lp\vG_{c,\xi}, \x_{A_E}\rp$.
We now look at the formula (9.4.2).  Take $H$ such that
$H = Z_G(t)$ for $t \in Z_{L_1}^{0F}$.  Then the set of $g$ 
in the second sum in (9.4.2) corresponds to the set of 
semisimple conjugacy classes in $Z_M(\la_c)^F$ which are 
conjugate to a fixed $t$ in $G^F$.   
Let $e_H$ be the number of $g$ occurring in the second sum 
in (9.4.2).  Then the above observation implies that 
$e_H$ is bounded by a positive integer 
independent of $q$.
\par
Now by substituting (9.9.3) into (9.4.2), together with the above
remark, we have 
\begin{equation*}
\tag{9.10.1}
\lp \vG_{c,\xi}, \x_{A_E}\rp_{G^F} = 
   Q\xi(z_E)\e(cc_0)\iv + qm\iv\b
\end{equation*}
with some integer $m$ independent of $q$ and $\b \in \CA$.  Here
\begin{equation*}
\tag{9.10.2}
Q = \a_E\z_{\CI_0}\iv\p(\hat z_{w_L'}\iv)\bar\th_0^{w_L}(x\iv_E z_1) 
          |\CW_{\th_1}|\iv\sum_H |\CW_{H,\th_1}|e_Hf_H\iv g_H,
\end{equation*}
where $H$ runs over the subgroups such that $\d_H = 1$.
$Q$ is independent of $c, \xi$, and also is contained in a finite subset 
of $\CA$ independent of $q$.
\par
By Theorem 7.9, 
$\x_{A_E}$ coincides with $\nu_ER_x$ for a certain
$x = (c', \xi') \in \bar A_{\la}^F \times (\bar A_{\la}\wg)^F$, 
where $\nu_E \in \Ql^*$ is a certain root 
of unity.  By (4.5.1) and (4.5.2), together with Theorem 2.6, 
we see that 
\begin{equation*}
\tag{9.10.3}
\lp \vG_{c,\xi}, \nu_x R_x\rp_{G^F} = 
         \nu_E|\bar A^F_{\la}|\iv \xi(c')\xi'(c).
\end{equation*}
On the other hand, suppose that $\b \ne 0$ in (9.10.1).  
Then the absolute value of $qm\iv\b$ turns out to be very large 
if we choose $q$ large enough since $\b$ is contained in the ring 
$\CA$ of algebraic integers of the fixed cyclotomic field. 
This implies that the absolute value of 
$\lp \vG_{c,\xi}, \x_{A_E}\rp_{G^F}$ becomes very large  
since $Q\xi(z_E)\e(cc_0)\iv$ is contained in a finite 
subset of $\CA$ independent of $q$.  
This contradicts the formula (9.10.3), and we conclude that
$\b = 0$, and we have 
\begin{equation*}
\tag{9.10.4}
\lp \vG_{c,\xi}, \x_{A_E}\rp_{G^F} = Q\xi(z_E)\e(cc_0)\iv.
\end{equation*}
Comparing (9.10.3) and (9.10.4), we see that 
\begin{equation*}
c' = z_E, \quad \xi' = \e\iv, \quad \nu_E = Q\e(c_0)\iv|\bar A^F_{\la}|.
\end{equation*}
Note that the last equality implies that 
\begin{equation*}
\tag{9.10.5}
  |\CW_{\th_1}|\iv\sum_H|\CW_{H,\th_1}|e_Hf_H\iv g_H 
      = \pm |\bar A^F_{\la}|\iv,
\end{equation*}
where $H$ runs over the $F$-stable reductive subgroups of $G$ 
containing $L$ such that $\d_H = 1$, which is of the form 
$H = Z_G(t)$ for $t \in Z_L^0$ satisfying the property that 
there exists $g \in G^F$ such that $g\iv tg \in Z_M(\la_c)^F$. 
Then we have   
\begin{equation*}
\tag{9.10.6}
\nu_E = \pm \z_{\CI_0}\iv\e(c_0)\iv\a_E\p(\hat z_{w_L'}\iv)\bar\th_0^{w_L}(x\iv_E),
\end{equation*}
and the signature $\pm 1$ is determined by (9.10.5).
Thus Proposition 9.2 is proved.
\para{9.11.}
Returning to the setting in Theorem 8.6, we shall show that
the statement (i) in Theorem 8.6 holds without the restriction 
on $q$. The argument is divided into two steps.  First we show 
that Lusztig's conjecture holds without restriction on $q$, 
and that the scalar constants are  determined explicitly, 
by applying the specialization argument based on the Shintani 
descent identities of character sheaves
in [S1, Corollary 2.12].  In the second step, we show that the 
parametrization of almost characters as given in Theorem 8.6
holds without the restriction on $q$.
\para{9.12.}
For a positive integer $c$, we put 
$\CP_c = \{ r \in \BZ_{\ge 1} \mid r \equiv 1 \pmod c\}$.
Let $s \in G^*$ be as in (8.1.1).  We choose a positive
integer $c$ such that $F^c$ acts trivially on 
$A_{\la}$, and that $F^c$ stabilizes $\ds$ and $z_L$.
Then for any $r \in \CP_c$, $\ol\CM^{(r)}_{s, N}, \CM^{(r)}_{s, N}$ 
(objects for $G^{F^r}$) are naturally identified with 
$\ol\CM_{s,N}, \CM_{s, N}$ (objects for $G^F$).  
We denote by $\r_x^{(r)}$ (resp. 
$R_y^{(r)}$ the irreducible character (resp. the almost character) of 
$G^{F^r}$ corresponding to $x \in \ol\CM^{(r)}_{s,N}$ 
(resp. $y \in \CM^{(r)}_{s, N}$).
In particular, the set $\CE(G^{F^{mr}}, \{ s\})^{F^r}$ is 
naturally identified with the set $\CE(G^{F^m},\{ s\})^F$ for a
sufficiently divisible $m$.  
\par
Similarly, the set 
$\CM^{L, (r)}_{s_L, N_0}$ (for $L^{F^r}$) is identified with the set
$\CM^L_{s_L,N_0}$ (for $L^F$).
Let $\d^{(mr)} = \d^{(mr)}_{z,\e}$ be a cuspidal 
irreducible character of $L^{F^{mr}}$ corresponding to 
$(z,\e) \in \CM^L_{s_L, N_0}$.  Then $\CW_{\d^{(mr)}}$
and $\CZ_{\d^{(mr)}}$ are independent of $r$, which we denote by
$\CW_{\d}, \CZ_{\d}$.  
By (6.7.2) and (7.9.1), we have $\CW_{\d} = \CW_{\th_1}$.
For each $E \in (\CW_{\th_1}\wg)^{F''}$, we consider the character
sheaf $A_E$, and its characteristic function $\x^{(r)}_{A_E}$ with 
respect to $F^r$.  We also consider the almost character 
$R_{E}^{(r)}$ of $G^{F^r}$.  It follows from the proof of Theorem 7.9
(see (7.9.4)) that $\x_{A_E}^{(r)}$ coincides with $R^{(r)}_{E}$ 
up to scalar.  Then by Proposition 9.2 that there exists a 
positive integer $r_0$ such that this scalar $\nu^{(r)}_E$ is 
described by (9.10.6) and (9.10.5) if $r > r_0$. We put 
$\CP_c' = \{ r \in \CP_c \mid r > r_0\}$.
One can check that $\nu_E^{(r)}$ is independent of the choice of 
$r \in \CP'_c$, by replacing $c$ by its appropriate multiple, if
necessary. We denote by $\nu_E$ this common value $\nu_E^{(r)}$. 
We have the following proposition.    
\par
\begin{prop} 
Let the notations be as in 9.12.
Assume that $q$ is arbitrary.  Then for  
each $E \in (\CW_{\th_1}\wg)^{F''}$, we have 
\begin{equation*}
\x_{A_{E}} = \nu_{E}R_{E}.
\end{equation*}
\end{prop}
\begin{proof}
By Theorem 4.7, we have 
\begin{equation*}
\tag{9.13.1}
Sh_{F^{mr}/F^rw}(\tilde\d^{(mr)}|_{L^{F^{mr}}\s^rw})
   = \mu_0(R_{z,\e}^{L_w})^{(r)}
\end{equation*} 
for $w \in \CZ_{\d}$, where $(R^{L_w}_{z,\e})^{(r)}$ is the almost 
character of $L_w^{F^r}$ corresponding to $(z,\e)$.
Note that by Theorem 4.7, $\mu_0$ is independent of the choice 
of $r \in \CP_c$ under the appropriate choice of the extension 
$\wt\d^{(mr)}$ of $\d^{(mr)}$.
\par
Let $A_0 = \CE\otimes A_{z, \e}$ be the cuspidal character sheaf on $L$ 
as in Theorem 7.9.  We have 
$\CW_{\CE_1} = \CW_{\d}$ and $\CZ_{\CE_1} = \CZ_{\d}$. 
Then $A_0$ gives rise to an $F^r$-stable character 
sheaf $A^w_0$ of $L_w$ for $w \in \CZ_{\d}$.  We denote by 
$\x^{(r)}_{A^w_0}$ the characteristic function on $L_w^{F^r}$
induced from the isomorphism $\f_0^w: (F^r)^*(A_0^w) \isom A_0^w$.
Then by Theorem 7.9, we have
\begin{equation*}
\tag{9.13.2}
\x^{(r)}_{A^w_0} = \nu_0(R_{z, \e\iv}^{L_w})^{(r)},
\end{equation*}
where $\nu_0 = \z_{\CI_0}\iv\e(c_0)\iv$ is independent of the choice
of $r \in \CP_c$.
\par
Now the map $a_{F^rw}: C(L^{F^r}\ssim_{F^rw}) \to 
      C(G^{F^{mr}}\ssim_{F^r})$ is defined as in [S2, 3.3].
By [S2, (3.6.3)], for $w = w_{\d}y \in \CZ_{\d}$, we have 
\begin{equation*}
\begin{split}
a_{F^rw}(\wt\d^{(mr)}&|_{L^{F^{mr}}\s^rw}) \\
  = &\sum_{E \in (\CW_{\d}\wg)^{F^r}}q_y^{-mr/2}q^{-\tilde l(w)mr/2}
      \Tr(T_{\g_{\d}y}, \wt E(q^{mr}))(\wt\r^{(mr)}_E|_{G^{F^{mr}}\s^r}),
\end{split}
\end{equation*} 
where $q_y$ is a power of $q$.  (For the notation, see [S2, 3.5]). 
By applying the 
Shintani descent operator to the above equality, and by using 
Theorem 4.7, we have 
\begin{equation*}
\tag{9.13.3}
\begin{split}
Sh_{F^{mr}/F^r}\circ &a_{F^rw}(\wt\d^{(mr)}|_{L^{F^{mr}}\s^rw}) \\
  = &\sum_{E \in (\CW_{\d}\wg)^{F^r}}q_y^{-mr/2}q^{-\tilde l(w)mr/2}
      \Tr(T_{\g_{\d}y}, \wt E(q^{mr}))\mu_ER^{(r)}_E,
\end{split}
\end{equation*}
where $\mu_E$ is also independent of the choice of $r \in \CP_c$.
\par
On the other hand, by the Shintani descent identity for character
sheaves ([S1, Corollary 2.12]), the following formula holds for 
each $r \in \CP_c$;    
\begin{equation*}
\begin{split}
(-1)^{d_{\Bw}}Sh_{F^{mr}/F^r}\circ a_{F^rw}\circ 
       N^*_{F^{mr}/F^rw}(\x_{A^w_0}^{(r)})
 = \sum_A\biggl(\sum_{i \in I_A }c_i\xi_i^r\biggr)\x_A^{(r)},
\end{split}
\end{equation*}
where $A$ runs over all the elements in $\wh G^F$, 
and $\{ \xi_i \mid i \in I_A\}$ and 
$\{ c_i \mid i \in I_A\}$ are certain finite subsets of $\Ql$
associated to $A$. 
Then combining (9.13.1), (9.13.2) and (9.13.3), we see that
\begin{equation*}
\begin{split}
&(-1)^{d_{\Bw}}\nu_0\mu_0\iv\sum_{E \in (\CW_{\d}\wg)^{F^r}}
  q_y^{-mr/2}q^{-\tilde l(w)mr/2}
      \Tr(T_{\g_{\d}y}, \wt E(q^{mr}))\mu_ER^{(r)}_E \\
  = &\sum_A\biggl(\sum_{i \in I_A }c_i\xi_i^r\biggr)\x_A^{(r)}.
\end{split}
\end{equation*}
By using the orthogonality relations for Hecke algebras (see 
[S2, (3.6.4)]), we can deduce a formula
\begin{equation*}
\tag{9.13.4}
P_{\CW_{\d}}(q^{mr})\mu_{E}R_{E}^{(r)} = 
     \sum_A\biggl(\sum_{\i \in I_A}d_i\z_i^r\biggr)\x_A^{(r)}
\end{equation*}
for certain subsets $\{ d_i \mid i \in I_A \}$ 
$\{ \z_i \mid i \in I_A\}$ of $\Ql$ (depending on $E$), where 
$P_{\CW_{\d}}(t)$ is a polynomial in $t$ (a generalization of 
Poincar\'e polynomial).
\par
As discussed in 9.12 we have 
\begin{equation*}
\tag{9.13.5}
\x^{(r)}_{A_{E}} = \nu_{E}R^{(r)}_{E}
\end{equation*}
for $r \in \CP_c'$. 
Substituting (9.13.5) into (9.13.4), 
we see that 
\begin{equation*}
\tag{9.13.6}
\sum_{i \in I_A}d_i\z_i^r = \begin{cases}
       \mu_{E}\nu_E\iv P_{\CW_{\d}}(q^{mr}) 
               &\quad\text{ if } A = A_{E}, \\
       0       &\quad\text{ otherwise}
                                \end{cases} 
\end{equation*}
for any $r \in \CP_c'$.
By applying a variant of Dedekind's theorem ([S1, (3.7.6)]), 
we see that (9.13.6) holds for any $r \in \CP_c$. 
In particular, substituting (9.13.6) into (9.13.4), we see that 
(9.13.5) holds for any $r \in \CP_c$.  By putting $r = 1$, 
we obtain the proposition. 
\end{proof}
\par
Next we show that 
\begin{lem}
For each $E \in (\CW_{\th_1}\wg)^{F''}$, we have 
\begin{equation*}
R_{E} = R_{z_E, \mu_1\iv}.
\end{equation*}
\end{lem}
\begin{proof}
Let $\CP_c'$ be as in the proof of Proposition 9.13.  We choose 
$r \in \CP_c'$ and take $m$ large enough so that $F^{mr}$ is 
a sufficiently divisible extension of $F^r$.  Then $F^{mr}$ is 
also sufficiently divisible for $F$.  Let $\d = \d^{(mr)}$ be 
a cuspidal irreducible character of $L^{F^{mr}}$.
Then by the definition of $\CP_c'$, $\CZ^{(r)}_{\d}$ is identified with
$\CZ^{(1)}_{\d}$.  It follows that for any $E \in (\CW_{\d}\wg)^{F''}$, 
we have
\begin{align*}
Sh_{F^{mr}/F^r}(\wt\r^{(mr)}_E|_{G^{F^{mr}}\s^r}) &=
\mu_{E}^{(r)}R_{E}^{(r)}, \\ 
   \quad Sh_{F^{mr}/F}(\wt\r^{(mr)}_E|_{G^{F^{mr}}\s})&= 
        \mu_{E}^{(1)}R_{E}^{(1)}
\end{align*}
with some constants $\mu_{E}^{(r)}, \mu_{E}^{(1)}$.
But by Proposition 9.2, we know that for $r \in \CP_c'$,
$R_{E}^{(r)}$ coincides with $R^{(r)}_{z_E,\e_1\iv}$ with
$(z_E, \e_1\iv) \in \CM_{s, N}$.  It follows that 
$\r^{(mr)}_E$ coincides with $\r^{(mr)}_{z_E, \e_1\iv}$, and 
so we have $R_{E}^{(1)} = R_{z_E, \e_1\iv}$.  This proves the lemma.
\end{proof}
\para{9.15.}
Combining Proposition 9.13 and Lemma 9.14, we see that 
(8.6.1) holds for any $q$.  This completes the proof of 
Theorem 8.6.

\par\medskip 
 


\bigskip
\begin{thebibliography}{[DLM2]}
\par  
\bibitem [C] {C:C} R. Carter, Finite groups of Lie type: Conjugacy 
classes and complex characters, A Wiley-Interscience Publication, 
1985.
\par
\bibitem [DLM1] {DLM1: DLM} F. Digne, G.I. Lehrer and J. Michel,
  On Gel'fand-Graev characters of reductive groups with disconnected
  centre, J. Reine angew. Math. {\bf 491} (1997), 131--147.
\par
\bibitem [DLM2] {DLM2} F.Digne, G.I. Lehrer and J. Michel,
The space of unipotently supported class functions on a finie 
reductive group, J. Alg.  
\par 
\bibitem [G] {G:G} M. Geck, A note on Harish-Chandra induction,
  Manuscripta Math. {\bf 80} (1993), 393--401.
\par 
\bibitem [HL] {HL:HL} R. Howlett and G. I. Lehrer, Induced cuspidal
  representations and generalized Hecke rings, Invent. Math.  
  {\bf 58} (1980), 37--64. 
\par
\bibitem [K1] {K1:K} N. Kawanaka, Generalized Gelfand-Graev
  representations and Ennola duality, in ``Algebraic groups and
  related topics,'' Advanced Studies in Pure Math., Vol {\bf 6},   
  Kinokuniya, Tokyo and North Holland, Amsterdam, 1985, 179--206.
\par
\bibitem [K2] {K2:K} N. Kawanaka, Generalized Gelfand-Graev
  representations of exceptional simple algebraic groups over a finite 
  field, I, Invent. Math. {\bf 84} (1986), 575--616.
\par
\bibitem [K3] {K3:K} N. Kawanaka, Shintani lifting and Gelfand-Graev
  representations, in ``The Arcata conference on representations  of
  finite groups,'' Proceedings of Symposia in Pure Math., 
  Vol.{\bf 47--1}, Amer. Math. Soc., Providence, RI, 1987, 147--163
\par
\bibitem [Le] {Le:Le} G. I. Lehrer, The characters of the finite
  special linear groups, J. of Algebra. {\bf 26} (1973), 564--583.
\par
\bibitem [L1] {L1} G. Lusztig,  ``Characters of reductive groups over
  a finite field'', Ann.\ of Math.\ Studies, Vol.{\bf 107}, Princeton
  Univ Press, Princeton , 1984.
\par
\bibitem [L2]{L2}G. Lusztig, Intersection cohomology complexes on 
a reductive group,
Invent. Math. {\bf 75} (1984), 205--272.
\par
\bibitem[L3] {L3}G. Lusztig, Character sheaves, I Adv. in
Math. {\bf 56} (1985),
193--237, II Adv. in Math. {\bf 57} (1985), 226--265, III, Adv. in
Math. {\bf 57}
(1985), 266--315, IV, Adv. in Math. {\bf 59} (1986), 1--63, V, Adv. in Math.
{\bf 61} (1986), 103--155.
\par
\bibitem[L4]{L4} G. Lusztig, On the character values of finite
Chevalley groups at unipotent elements, J. of Algebra, {\bf 104} (1986),
146--194.
\par
\bibitem [L5] {L5} G. Lusztig, On the representations of reductive
  groups with disconnected centre, Ast\'erisque {\bf 168} (1988),
  157--166.
\par
\bibitem[L6] {L6}G. Lusztig, Green functions and character sheaves, 
Annals of Mathematics, {\bf 131} (1990), 355--408.
\par
\bibitem [L7] {L7} G. Lusztig, A unipotent support for irreducible
  representations , Adv. in Math., {\bf 94} (1992), 139--179.
\par
\bibitem[LS]{LS} G. Lusztig and N. Spaltenstein, On the generalized
Springer correspondence for classical groups, Advanced Studies in Pure
Math. {\bf 6} (1985), pp.289--316.
\par
\bibitem[M]{M} I.G. Macdonald, ``Symmetric functions and Hall
  polynomials'', second edition.  Clarendon Press. Oxford 1995.
\par
\bibitem [S1] {S1:S} T. Shoji, Character sheaves and almost characters
of reductive groups, II, Adv. in Math. {\bf 111} (1995), 314--354.
\par  
\bibitem [S2] {S2:S} T. Shoji, Shintani descent for special linear 
groups, J. Algebra {\bf 199} (1998), 175--228.
\par 
\bibitem [S3] {S3} T. Shoji, Generalized Green functions and 
unipotent classes for finite reductive groups, in preparation.
\par
\bibitem [ShS] {ShS} T. Shoji and K. Sorlin, Subfield symmetric 
spaces for finite special linear groups, to appear in Rep.
Theory. 
\bibitem [SpS] {SS:SS} T.A. Springer and R. Steinberg, Conjugacy
  classes, in \lq\lq Seminar on Algebraic groups and related topics'', 
  Lecture Note in Math., Vol. {\bf 131},  Part E, Springer-Verlag,
  1970.
\end{thebibliography}
\end{document}